\documentclass[10pt, reqno]{amsart}
\usepackage{latexsym}
\usepackage{amsfonts, amsmath, amscd}
\usepackage{amssymb}
\usepackage[all]{xypic}
\usepackage{color} 
\usepackage{enumerate}
\usepackage{verbatim}

\usepackage{wrapfig}

\usepackage{amsthm}    

\usepackage{graphicx, psfrag}
\usepackage{pstool}
 \usepackage{caption}
\usepackage{subcaption}

\numberwithin{equation}{section}
\renewcommand{\theequation}{\arabic{section}.\arabic{equation}}

\def\be{\begin{equation}}
\def\ee{\end{equation}}
\def\bear{\begin{eqnarray}}
\def\eear{\end{eqnarray}}
\def\best{\begin{eqnarray*}}
\def\eest{\end{eqnarray*}}

\renewcommand{\theequation}{\arabic{section}.\arabic{equation}}
\renewcommand{\thesubsection}{\arabic{section}.\arabic{subsection}}

\newtheorem{theorem}{Theorem}[section]
\newtheorem{prop}[theorem]{Proposition}

\newtheorem{lemma}[theorem]{Lemma}
\newtheorem{cor}[theorem]{Corollary}

\newtheorem{defn}[theorem]{Definition}

\newtheorem{remark}[theorem]{Remark}
\newenvironment{rem}{\begin{remark}\rm}{\end{remark}}

\newtheorem{example}[theorem]{Example}
\newenvironment{ex}{\begin{example}\rm}{\end{example}}

\def\non{\noindent}

\def\al{\alpha}

\def\om{\omega}
\def\bd{\partial}

\def\dim{\mathrm{dim\;}}
\def\de{\delta}

\def\ga{\gamma}
\def\la{\lambda}

\def\ep{\varepsilon}
\def\Hom{\mathrm{Hom}}
\def\J{\mathcal{J}}

 \def\ker{\mathrm{Ker\; }}
 \def\cok{\mathrm{Coker \;}}

\def\ra{\rightarrow}

\def\r#1{\right#1}
\def\l#1{\left#1}

\def\ma#1{\mathop {#1} \limits}

\def\Si{\Sigma}
\def\De{\Delta}

\def\ti{\times}

\def\N{{ \mathbb N}}
\def\Z{{ \mathbb Z}}
\def\R{{  \mathbb R}}

\def\P{{ \mathbb P}}
\def\Q{{ \mathbb Q}}
\def\cx{{ \mathbb C}}

\def\La{\Lambda}

\def\wt#1{\widetilde{#1}}
\def\wh#1{\widehat{#1}}
\def\ov#1{\overline{#1}}

\def\M{{\mathcal M}}
\def\U{{\mathcal U}}
\def\oM{\ov\M}

\def\longra{\longrightarrow}
\def\H1{{\mathcal H}^1}
\def\si{\sigma}
\def\Si{\Sigma}
\def\Aut{\mbox{Aut}}

\def\F{\mathbb F}

\def\phi{\varphi}
\def\longra{\longrightarrow}

\def\x{{\mathcal X}}

\def\L{{\mathcal L}}

\def\op#1{\overset{\circ}{#1} \vphantom{V}}

\def\Ev{\mathrm{Ev}}
\def\ev{\mathrm{ev}}
\def\st{\mathrm{st}}
\def\ct{\mathrm{ct}}
\def\pr{\mathrm{pr}}
\def\pr{\mathrm{pr}}
\def\ft{\mathrm{ft}}
\def\De{\Delta}

\def\raG{ \ma\longra^G}
\def\J{\mathcal{J}}
\def\JV{\mathcal{JV}}

\def\fle{\preccurlyeq}
\def\fge{\succcurlyeq}

\def\fg{\succ}
\def\s{\mathcal{S}}

\begin{document}

\title[Relative GW]{GW Invariants Relative normal crossing Divisors} 

\vskip.2in
\author{Eleny-Nicoleta Ionel}
\address{Department of  Mathematics,  Stanford University}
\email{ionel@math.stanford.edu}

\thanks{Research supported in part by the NSF grants DMS-0605003 and DMS-0905738} 

 \begin{abstract}
In this paper we introduce a notion of symplectic normal crossing divisor $V$ and define the  GW invariant of a symplectic manifold $X$ relative such a divisor. Our definition includes normal crossing divisors from algebraic geometry. The invariants we define in this paper are key ingredients in symplectic sum type formulas for GW invariants, and  extend those defined in our previous joint work with T.H.  Parker \cite{ip1}, which covered the case $V$ was smooth. The main step is the  construction of a compact moduli space of relatively stable maps into the pair $(X, V)$ in the case $V$ is a symplectic normal crossing divisor in $X$. 
\end{abstract}

\maketitle

\setcounter{section}{-1}
\setcounter{equation}{0}
\section{Introduction}

In previous work with Thomas  H. Parker \cite{ip1} we constructed the relative Gromov-Witten invariant $GW(X,V)$ of a closed symplectic manifold $X$ relative a smooth ``divisor'' $V$, that is a (real)  codimension 2 symplectic submanifold. These relative invariants  are defined by choosing  an almost complex structure  $J$ on $X$ that is compatible with both  $V$ and the symplectic form, and counting $J$-holomorphic maps that intersect $V$ with specified multiplicities.  An important application is the  symplectic sum formula that relates the GW invariant of a symplectic sum $X\#_V Y$ to the relative GW invariants of  $(X, V)$ and $(Y, V)$ (see  \cite{ip2}  and the independent approaches \cite{lr}, \cite{li} and  \cite{sft}).

In this paper we introduce a notion of symplectic normal crossing divisor $V$, and define the  GW invariant of a symplectic manifold $X$ relative such a divisor. Roughly speaking, a set $V\subset X$ is a symplectic normal crossing divisor if it is locally the transverse intersection of codimension 2 symplectic submanifolds compatible with $J$ (the precise definition is given in Section \ref{s1}).
 
There are many reasons why one would want to extend the definition of relative GW invariants to include normal crossing divisors, and we already have several interesting applications in mind. One is a Mayer-Vietoris type formula for the GW invariants: a formula describing how the GW invariants behave when $X$ degenerates into several  components and that allows one to  recover the invariants of $X$ from those of the components of the limit.   The simplest such degenerations come from the symplectic sum
along a smooth divisor. But if one wants to iterate this degeneration, one is immediately confronted with several pieces whose intersection is no longer smooth, but instead are normal crossing divisors.  Normal crossing divisors appear frequently in algebraic geometry, not only as the central fiber of a stable degeneration but also for example as the toric divisor in a toric manifold which then appears in the context of mirror symmetry.   We also have some purely symplectic applications in mind in which normal crossing divisors arise from Donaldson's theorem; these appear in a separate paper \cite{ip-don}.

The general approach in this paper is to appropriately adapt the  ideas in \cite{ip1}  but now allow the divisor to have a simple type of singularity, which we call symplectic normal crossing. This is defined in Section \ref{s1}, where we also present many of the motivating examples. The notion of simple singularity is of course relative: the main issue here is to be able to control the analysis of the problem; the topology of the problem, though perhaps much more complicated is essentially of a combinatorial nature so it is much easier controlled.

There are several new features and problems that appear when the divisor $V$ has such singular locus.  First, one must include in the moduli space holomorphic curves that intersect the singular locus, and one must properly record the contact information about such intersections. In Section \ref{s2} we describe how to do this and construct the corresponding moduli space $\M_s(X,V)$ of stable maps into $X$ whose contract intersection with $V$ is described by the sequence $s$.  There is a lot of combinatorics lurking in the background that keeps track of the necessary  topological information along the  singular locus, which could make the paper unnecessarily longer. We have decided to keep the notation throughout the paper to a minimum, and expand its layers only as needed for accuracy in each section.  We give simple examples of why certain situations have to be considered, explain in that simple example what needs to be done, and only after that proceed to describe  how such situations  can be handled in general. In the Appendix we describe various needed stratifications associated to a normal crossing divisor, and  topological data associated to it.

\smallskip
The other more serious problem concerns the construction of a compactification $\ov \M_s(X,V)$ of the relative moduli space. In the usual  Gromov compactification of stable maps into $X$, a  sequence of holomorphic maps that have a prescribed contact to $V$ may limit to a map that has components in $V$ or even worse in the singular locus of $V$; then not only the contact information is lost in the limit, but the formal dimension of the corresponding boundary stratum of the stable map compactification is greater than the dimension of the moduli space. This problem already appeared for the moduli space relative a smooth divisor, where the solution was to rescale the target normal to $V$ to prevent components from sinking into $V$; but now the problem is further compounded by the presence of the singular locus of $V$.  So the main issue now is to how to precisely refine the Gromov compactness and construct an appropriate relatively stable map compactification $\ov \M_s(X,V)$ in such a way that its boundary strata are not larger dimensional than the interior.

In his unpublished Ph.D. thesis, Joshua Davis \cite{da} described how one can construct a relatively stable map compactification for the space of genus zero maps relative a normal crossing divisor, by recursively blowing up the singular locus of the divisor. As components sunk into this singular locus, he recursively blew it up to  prevent this from happening. This works for genus zero, but unfortunately not in higher genus. The main reason for this is that in genus zero a dimension count shows that components sinking into $V$ cause no problem, only those sinking into the singular locus of $V$ do. However, that is not the case in higher genus, so then one would also need to rescale around $V$ to prevent this type of behavior. But then the process never terminates: Josh had a simple example in higher genus where a component would sink into the singular locus. Blowing up the singular locus forced the component to fall into the exceptional divisor. Rescaling around the exceptional divisor then forced the component to fall back into the next singular locus, etc. 
\smallskip

In this paper we present a different way to construct a relatively stable map compactification $\ov \M_s(X,V)$, by instead  rescaling $X$ simultaneously normal to all the branches of $V$, a procedure we describe in Section \ref{s3}. When done carefully, this is essentially a souped up version of the rescaling procedure described in  \cite{ip1} in the case $V$ was smooth. Unfortunately, the naive compactification that one would get by simply importing the description of that in \cite{ip1} simply does not work when the singular locus of $V$ is nonempty! There are two main reasons for its failure: the first problem is that the ``boundary stratum'' containing curves with components over the singular locus is again larger dimensional than the ``interior'' so it is in some sense too big;  the second problem is that it still does not capture all the limits of curves sinking into the singular locus, so it is too small! This seems to lead into a dead end, but upon further analysis in Sections \ref{s5} and \ref{s6}  of the limiting process near the singular locus two new features appear that allows us to still proceed. 
\smallskip

The first new feature is the refined matching condition that the limit curves must satisfy along the singular locus of $V$. It turns out that not all the curves which satisfy the naive matching conditions can appear as limits of maps in $\M_s(X,V)$. The naive matching conditions require that the curves intersect $V$ in the same points with matching order of contact, as was the case in \cite{ip1}, while the refined ones along the singular locus essentially require that their slopes in the normal directions  to $V$ also match. So the refined matching conditions also involve the leading coefficients of the maps in these normal directions, and they give conditions in a certain weighted projectivization of the normal bundle to the singular locus, a simple form of which is described in  Section \ref{s4}. Luckily, this is  enough to cut back down the dimensions of the boundary to what should be expected. In retrospect, these refined matching conditions already appeared in one of the  key Lemmas in our second joint paper \cite{ip2} with 
Thomas H.  Parker on the symplectic sum formula, but they do not play any role in the first paper \cite{ip1} because they are automatically satisfied when $V$ is smooth. 
\smallskip

The second  new feature that appears when $V$ is singular is that unfortunately one cannot avoid trivial components stuck in the neck (over the singular  locus of $V$), as we show in some simple examples at the end of Section \ref{s4}. This makes the refined matching conditions much more tricky  to state, essentially because these  trivial components do not have the right type of leading coefficients. The solution to this problem is to realize that the trivial components are there only to make the maps  converge in Hausdorff distance to their limit, and in fact they do not play any essential role in the compactification, so one can simply collapse them in the domain, at the expense of allowing a node of the collapsed domain to be between not necessarily consecutive levels. The refined matching condition then occurs only at nodes between two nontrivial components, but needs to take into account this possible jump across levels. It is described more precisely in Section \ref{s5}.
\smallskip

This finally allows us to define in Section \ref{s6} the compactified moduli space $\ov\M_s(X, V)$ of relatively $V$-stable maps into 
$X$,  which comes together with a continuous map
\bear\label{st.ti.Ev}
\st\ti \Ev: \ov \M_s (X, V) \ra \ov \M_{\chi_s, \ell(s)} \ti    \prod _{ x}  \P_{s(x)} (N V_{I(x)})
\eear
The first factor is the usual stabilization map recording the domain of $f$, which may be disconnected,  but the new feature is the second factor $\Ev$. It is a refinement of the usual (naive) evaluation map $\ev$ at  the points $x$ that are mapped into the singular 
locus of $V$, and it also records  the weighted projectivization of the leading coefficients of  $f$  at $x$ in all the normal directions to $V$ at $f(x)$. This is  precisely the map that appears in the refined matching conditions. 
 
In Section \ref{s7} we then show that for generic $V$-compatible perturbation $(J,\nu)$ (under the assumption of Remark \ref{nu.univ}) the  image of $\ov \M_s(X,V)$ 
under the map (\ref{st.ti.Ev}) defines a homology class $GW_s(X, V)$ in dimension 
\best
\dim \ov \M_s(X,V)= 2 c_1(TX)A_s + (\dim X-6)\frac {\chi_s} 2 + 2\ell(s) - 2A_s\cdot V
\eest
called the GW invariant of $X$ relative the normal crossing divisor $V$. The class $GW_s(X, V)$ is independent of the perturbation $\nu$ and is invariant under smooth deformations of the pair $(X,V)$ and of $(\omega, J)$ though $V$-compatible structures.  When $V$ is smooth these invariants agree with the usual relative GW invariants as constructed in \cite{ip1}.

\medskip
There is a string of recent preprints that have some overlap with the situation considered in our paper, in that they all generalize in some way the normal crossing situation from algebraic geometry. First of all, there is certainly an overlap between what we call a symplectic normal crossing divisor in this paper and what fits into the exploded manifold setup considered by Brett Parker \cite{expl}.  There is also some overlap with the logarithmic Gromov-Witten invariants \cite{gs} considered by Gross and Siebert in the context of algebraic geometry (see also the Abramovich-Chen paper \cite{ac} on a related topic). However, the precise local structure near the divisor is very different: log geometry  vs symplectic normal crossing vs  exploded structures. Furthermore, the moduli spaces constructed in these papers and in particular their compactifications are completely different, even when applied to the  common case  when $V$ is a smooth divisor in a smooth projective variety, see Remarks \ref{b-p} and \ref{gs-e} for more details. This means that a priori even in this common case each one of these other approaches many lead to different invariants, some even different from the usual relative GW invariants. 
\medskip

This paper is based on notes from a talk the author gave in Sept  2006 in the Moduli Space Program at Mittag-Leffler Institute, during a month long stay there.  The notes were expanded in the fall of 2009 during the Symplectic and Contact Geometry and Topology program at MSRI. We thank both research institutes for their hospitality.  The author would also like to thank the referees for their tireless requests to add more details to the paper which we hope helped improved the exposition. 

\setcounter{equation}{0}
\section{Symplectic normal crossing divisors}\label{s1}
\medskip

In this section we define a notion of symplectic normal crossing divisors, generalizing that from  algebraic geometry, and encoding the geometrical information required for the analysis of \cite{ip1} and \cite{ip2} to extend after appropriate modifications. Clearly the local model of such divisor $V$ should be the union of $k$ coordinate planes in $\cx^n$, where the number of planes may vary from point to point. But we also need a local model for the symplectic form  $\omega$ and the tamed almost complex structure $J$ near such divisor, and we require that each branch of  $V$ is  both $\omega$-symplectic and $J$-holomorphic. This will allow us to define the order of contact of $J$-holomorphic curves to $V$. We also need a good description of the normal directions to the divisor, as these are the directions in which the manifold $X$ will be rescaled when components of the $J$-holomorphic curves fall into $V$. So we keep track of both the normal bundle to each branch of $V$ and its inclusion into $X$ describing the neighborhood of that branch. 
\begin{defn}[Local model]\label{def.loc.mod} In $\cx^n$, consider the union $V$ of $k\ge 0$ (distinct) coordinate hyperplanes 
$H_i=\{x |x_i=0\} $, together with their normal direction $N_i$ defined by the projection $\pi_i:\cx^n \longra H_i$, $\pi_i(x)=x_i$,  and the inclusion $\iota: (N_i, 0)  \ra (\cx^n, H_i)$. We say that these form a {\em model for a normal crossing divisor}  in $\cx^n$ with respect to a pair  
$(\omega, J)$ if  all the divisors $H_i$ are both $\omega$-symplectic and $J$-holomorphic. 
\end{defn}
\begin{rem} There is a natural action of  $\cx^*$  on the model induced by scaling by a factor of $t^{-1}$ in the normal direction 
to each $H_i$,  for $i=1, \dots, k$. This defines a rescaling map $R_t: \cx^n\ra \cx^n$ for  $t\in \cx^*$. By construction, the $R_t$ leaves the divisors $H_i$ invariant, but not pointwise, and may not preserve $J$.  However, as $t\ra 0$,  $R_t^* J$ converges uniformly on compacts to a $\cx^*$ invariant limit  $J_0$ which depends on the 1-jet of $J$ along the divisor. 
\end{rem}

\begin{defn}\label{D.ncd} Assume $(X,\omega,J)$ is a symplectic manifold with a tamed almost complex structure.  $V$ is called {\em a normal 
crossing divisor} in $(X,\omega, J)$ with normal bundle $N$ if there exists a smooth manifold $\wt V$ with a complex line bundle 
$\pi: N\ra \wt V$ and an immersion $\iota: U_V\ra X$ of some disk bundle $U_V$ of  $N_V$ into $X$ satisfying the following properties: 
\begin{itemize}
\item[(i)] $V$ is the image of the zero section $\wt V$ of $N$
\item[(ii)] the restriction of $\iota^* J$ to the fiber of $N$  along the zero section induces  the complex multiplication in the bundle $N$.  
\item[(iii)] at each point $p\in X$ we can find local coordinates on $X$ in which the configuration $(X, \pi, \iota, V)$ becomes identified with one of the local models in  Definition \ref{def.loc.mod}.   
\end{itemize}
Such a pair $(J,\omega)$ is called  {\em compatible} with the divisor $V$. $N$ is called the {\em normal bundle} of $V$ and $\wt V$ the (smooth) {\em resolution} of $V$.  A connected component of $\wt V$ is called a component or {\em (global) branch} of $V$. 
  \end{defn} 
Note that $\iota$ induces by pullback from $X$ both a symplectic structure $\omega$ and an almost complex structure $J$ on the 
total space of the disk bundle in $N$ over on $\wt V$, which serves as a global model of $X$ near $V$.  Its zero section $\wt V$ is both symplectic and $J$-holomorphic and serves as a smooth model of the divisor $V$. $N$ is also a complex line bundle whose complex structure comes from the restriction of $J $ along the zero section.  Thus $N$ also comes with a $\cx^*$ action which will be used to rescale $X$ normal to $V$. 

\begin{rem}  We are not requiring $J$ to be locally invariant under this $\cx^*$ action. We also are not imposing the condition that the branches are perpendicular with respect to $\omega$  or that the projections $\pi_i$ are $J$-holomorphic.  We also allow transverse self intersections of various components of $V$. When each component of $V$ is a  submanifold of $X$  the divisor is said to have {\em  simple} normal crossing singularities.  Any of these  assumptions would simplify some of the arguments, but are not needed. 
\end{rem}
In this paper we only work with  $J$'s compatible with  $V$ in the sense of Definition 3.2 of \cite{ip1}. This is  a condition on the normal 1-jet of $J$ along $V$: 
 \begin{enumerate}
 \item[(b)] \hskip.2in $ [(\nabla_\xi J + J\nabla_{J\xi}J)(v)]^N = [(\nabla_vJ)\xi + J (\nabla_{Jv}J)\xi]^N$ for all $v\in TV$, $\xi \in NV$;
 \end{enumerate}
discussed in the Appendix. This extra condition is needed to ensure that the stable map compactification has codimension 2 boundary strata, so it gives an invariant, independent of  parameters. A priori, even when $V$ is smooth the relatively stable map compactification may have real codimension 1 boundary strata without this extra assumption.  

\smallskip

A symplectic normal crossing divisor could be defined locally in terms of an atlas of charts on 
$(X, \om, J)$ compatible with  $V$. The local models suffice to construct both a smooth resolution $\wt V$ of $V$ by separating its local branches as well as the complex normal bundle $N$, in effect proving a tubular neighborhood theorem in this context. For simplicity of exposition, we decided instead to include the global existence of $\wt V$ and $N$ as part of the definition of a normal crossing divisor.  

One could also define a notion of normal crossing in the smooth or even orbifold category. In this paper we insist that the normal bundle $N$ carry a complex structure, which induces a local $\cx^*$-action normal to $V$. Otherwise, one only has an $\R_+$ action, which is the typical situation in SFT, leading to further complications. 
\begin{ex}
A large class of examples is provided by  algebraic geometry.  Assume $X$ is a smooth  projective variety and $V$ a smooth normal crossing divisor in $X$ in this category (i.e. the normalization of 
$V$ is a smooth projective variety).  
Then $V$ is a symplectic normal crossing divisor for $(X, J_0, \omega_0)$ where $J_0$ is the integrable complex structure and $\omega_0$ the Kahler form. For example  (a) $X$ could be a Hirzebruch surface and $V$ the union of the zero section, the  infinity section and several fibers  or  (b) $V$ could be the union of a section and a nodal fiber in an elliptic surface $X$. 

An important example of this type is when $X$ is a toric manifold and $V$ is its toric divisor, which is a case considered in mirror symmetry, see for example \cite{au2}. 
\end{ex} 
\begin{ex} \label{ex.cp2} A particular example to keep in mind is  $X=\cx\P ^2$ with a degree 3 normal crossing divisor $V$. For example $V$ could be a smooth elliptic curve, or $V$ could be a nodal sphere, or finally $V$ could be a union of 3 distinct lines.  In the second case the resolution $\wt V$ is  $\cx\P^1$ with normal bundle $O(7)$ while in the last case it is $\cx\P^1\sqcup \cx\P^1\sqcup\cx \P^1$, each component with normal bundle $O(1)$. Of course, in a complex 1-parameter family, a smooth degree three curve can degenerate into either one of the other two cases. 

Another motivating example of this type comes from a smooth quintic 3-fold degenerating to a union of 5 hyperplanes in $\cx \P^4$.  
\end{ex}
\begin{rem}\label{DM} Another special case is $X=\ov \M_{0,n}$ the Deligne-Mumford moduli space of stable genus  $0$ curves and $V$ the union of all its boundary strata (i.e. the stratum of nodal curves).  The usual description of each boundary stratum and of its normal bundle  provides the required local models for a symplectic normal crossing divisor. This discussion can also be extended to  the orbifold setting to cover the higher genus case $\ov\M_{g, n}$ and includes its smooth finite  covers, the moduli space of Prym curves \cite{lo} or the moduli space of twisted $G$-covers \cite{acv}. 
\end{rem} 
Of course, there are many more symplectic examples besides those coming from algebraic geometry. 
\begin{ex}\label{V.sm} Assume $V$ is a symplectic codimension two submanifold of $(X, \omega)$. The symplectic neighborhood theorem allows us to find a $J$ and a model for the normal direction to $V$, so $V$ is  normal crossing divisor in $(X, \omega, J)$. Of course in this case  $V$ is a smooth divisor, so it has empty singular locus. 
\end{ex}
One may have hoped that the union of  several transversely intersecting codimension two symplectic submanifolds would  similarly be a normal crossing divisor.  Unfortunately, if the singular locus is not empty,  that may not be the case: 
\begin{ex}\label{ex1} Let  $V_1$ be an exceptional divisor in a symplectic 4-manifold  and $V_2$ a sufficiently small generic perturbation of it, thus still a symplectic submanifold, intersecting transversely $V_1$. This configuration cannot be given the structure of a normal crossing divisor, simply because one cannot find a $J$ which preserves both. If such a $J$ existed, then all the intersections between $V_1$ and $V_2$ would be positive, contradicting the fact that exceptional divisors have negative self intersection.  
\end{ex}

This example illustrates the fact that a normal crossing divisor is not a purely symplectic notion, but rather one also needs the existence of an almost  complex structure  $J$ compatible with  the crossings.  The positivity of intersections of all branches is a necessary condition for such a $J$ to exist in general. 

\begin{rem}  One could ask  what are the necessary and sufficient conditions for  $V$ inside a symplectic manifold $(X,\omega)$  to be a  normal crossing divisor with respect to some $J$ on  $X$. Clearly $V$ should be locally the transverse intersection of symplectic submanifolds, and the intersections should be  positive. If we assume that the branches of $V$ are orthogonal wrt $\omega$, the existence of an $\om$-compatible $J$  which is compatible with  $V$ is straightforward (see Appendix). In general, one might be able to use a homotopy argument to prove that positivity of intersections is the only obstruction to the existence of a $J$ compatible with  $V$ and tamed by $\om$.  We do not pursue this issue further in this paper. 
\end{rem}

\begin{ex}\label{ex2} Symplectic Lefschetz pencils or fibrations  provide another source of  symplectic normal crossing divisors.  Assume $X$ is a symplectic manifold which has a symplectic Lefschetz fibration with a symplectic section, for example one coming from Donaldson Theorem \cite{do2} where the section comes from blowing up the base locus.  Gomph \cite{g2} showed that in this case there is an almost complex structure $J$ compatible with this fibration. We could then take $V$ the union of the section with a bunch of fibers, including possibly some singular fibers. 
\end{ex}

\begin{ex}\label{Don.ex} (Donaldson divisors) Assume $V$ is a normal crossing divisor in 
$(X,\omega, J)$, $J$ is $\om$-compatible and $[\omega]$ has rational coefficients.  We can use Donaldson theorem \cite{do} to obtain a smooth divisor $D$ representing the Poincare dual of $k\omega$ for $k\gg0$ sufficiently large, such that $D$  is $\ep$-$J$-holomorphic and 
$\eta$-transverse to $V$ (see also \cite{au}). Choosing carefully the parameters $\eta$ and $\ep$, one can find a sufficiently small, $\om$-tamed deformation of $J$ such that $V\cup D$ is also a normal crossing divisor (cf. \cite{ip-don}). 
\end{ex}

\begin{rem}  The definition of a normal crossing divisor  works well under taking products of symplectic manifolds with divisors in them. If $V_i$ is a normal crossing divisor in $X_i$  for $i=1,2$ then 
$\pi_1^{-1}(V_1)\cup \pi_2^{-1} (V_2)= V_1\ti X_2\cup X_1\ti V_2$ is a normal crossing divisor in $X_1\ti X_2$, with normal model 
$\pi_1^*N_1 \sqcup \pi_2^* N_2$.   Note that even if $V_i$ were smooth divisors, the induced divisor in the product $X_1\ti X_2$ is singular along $V_1 \ti V_2$. 
\end{rem}

\begin{rem}\label{sym.sum.ex} The definition of a normal crossing divisor  also behaves well under symplectic sums. Assume $U_i\cup V$ is a symplectic divisor in $X_i$ for $i=1, 2$ such that the normal bundles of $V$ in $X_i$ are dual. If  $U_i$ intersect $V$ in the same divisor $W$ then Gomph's argument \cite{g} shows that the divisors $U_i$ glue  to give  a normal crossing divisor $U_1 \#_W U_2$  in the symplectic sum  $X_1\#_V X_2$. 
\end{rem}

\begin{rem} \label{N.strata}{\bf (Stratifications associated to a normal crossing divisor)} Any symplectic normal crossing divisor $V$ in $(X,\omega, J)$ induces a stratification of $X$, whose closed stratum $V^k$ consists of those points $X$ where at least $k$ local branches of $V$ meet.  Each closed stratum has a {\em smooth resolution}  $\wt {V^k}\ra V^k$ which comes with an induced $(\om, J)$ and an intrinsic symplectic normal crossing divisor $\wt V^{k+1}$  over the lower depth stratum $V^{k+1}$, as described in  \S \ref{S.A.nc.strat}. 

More precisely, for each point $x\in X$, its  {\em depth } $k(x)$ is the largest  $k$ such that $x\in V^k$, or equivalently the cardinality of $\iota^{-1}(x)$, where $\iota: \wt V  \looparrowright  X$ is the immersion parameterizing $V$. So a point in $X\setminus V$ has depth 0 while points in the singular locus of $V$ have depth at least 2. This defines an upper semi-continuous function $\textsf {depth}:X\ra \N $ whose level sets are the open strata where precisely $k$ local branches of $V$ meet. The fiber of $\wt {V^k}\ra V^k$ over a depth $k$ point  $x$ is intrinsically $\iota^{-1}(x)$ and keeps track of the $k$ local branches of $V$ meeting at $x$. In fact, for any finite set $I$ of order $k$, we get a resolution $V_I\ra V^k$; its fiber at a depth $k$ point $x$ consists of bijections $I \ra \iota^{-1}(x)$, with a symmetric group $S_I$ action reordering $I$, see \S \ref{S.A.nc.strat} for more details.  \end{rem}

\begin{rem}\label{b-p} A special case of symplectic normal crossing divisor $V$ (with simple crossings) is the union of codimension 2 symplectic submanifolds which intersect orthogonally wrt $\omega$, and whose local model matches that of toric divisors in a toric manifold. This is a case that fits in the exploded manifold set-up of  Brett Parker (see Example 5.3 in the recent preprint  
\cite {expl}), so in principle one should be able to compare the relative invariants we construct in this paper with the exploded ones of  \cite{expl}.  It is unclear to us what is exactly the information that the exploded structure records in this case, and what is the precise relation between the two moduli spaces. But certainly the relatively stable map compactification we define in this paper seems to be quite different from the exploded one, so it is unclear whether they give the same invariants, even in the case when $V$ is smooth. \end{rem}

\begin{rem}  \label{gs-e} In a related paper, Gross and Siebert define log GW invariants in the algebraic geometry setting \cite{gs}. If $V$ is a normal crossing divisor in a smooth projective variety $X$, then it induces a log structure on $X$. However, even when $V$ is a smooth divisor, Gross and Siebert explain that the stable log compactification they construct is quite different from the relatively stable map compactification constructed earlier in that context by J. Li \cite{li} (and which agrees with that of \cite{ip1} in this case). So a priori, even when $V$ is smooth, the usual relative GW invariants may be different from the log GW invariants of \cite{gs}. The authors mention that in that case at least there is a  map from the moduli space of stable relative maps to that of stable log maps, which they claim could be used to prove that the invariants are the same. Presumably there is also a map from the relatively stable map compactification that we construct in this paper to the appropriate stable log compactification in the  more general case when $V$ is a normal crossing divisor in a smooth projective variety.  

In another paper \cite{ac}  Abramovich and Chen explain how, in the context of algebraic geometry, the construction of a  log moduli space when $V$ is a normal crossing divisor (with simple crossings) follows from the case when $V$ is smooth by essentially functorial reasons. Again,  it is unclear to us how exactly the two notions of log stable maps of \cite{gs} and \cite{ac} are related in this case.  
\end{rem} 
\begin{rem} One note about simple normal crossing vs general normal crossing: they do complicate the topology/combinatorics of the situation, but if set up carefully the analysis is unaffected.  If  the local model of $X$ is holomorphic near $V$ (as is the case in last two examples above),  even if  $V$  did not have simple crossings, one could always  blow up the singular locus $W$ of $V$ to get a total divisor $ \pi^{-1}(V)=Bl (V)\cup E$ with simple normal crossing in $Bl(X)$, where $E$ is the exceptional divisor.  Blowing up in the symplectic category is a more delicate issue, but when using 
the appropriate local model, one can  always  express (a symplectic deformation) of the original manifold $(X, V)$ as a symplectic sum of its blow up $(Bl (X), Bl(V))$ along the exceptional divisor $E$ with a standard piece $(\P, V_0)$ involving the normal bundle of the blowup locus. Since we are blowing up the singular locus of $V$, the proper transform $Bl(V)$ intersects nontrivially the exceptional divisor $E$; the symplectic sum $Bl (X)\#_E \P=X$  then also glues $Bl(V)$ to the standard piece on the other side to produce $V$, as in Remark \ref{sym.sum.ex}. So a posteriori, after proving a symplectic sum formula for the relative GW of normal crossing divisors passing through the neck of a symplectic sum, one could also express the relative GW invariants of the original pair  $(X, V)$ as universal expressions in the relative GW invariants of its the blow up and those of the piece obtained from  the normal bundle of the blow-up locus.  
\end{rem}
For the rest of the paper, unless specifically mentioned otherwise,  a (normal crossing) divisor $V$ in $X$ will always mean a symplectic normal crossing divisor for some tamed pair $(\om, J)$ on $X$ which is compatible with  $V$ to first order. We denote by $\J(X, V)$ the space of such tamed pairs, described in  \S \ref{A.spaces.param}. The restriction of $(\om, J)$ to $V$ pulls back to a tamed pair on the smooth resolution of each depth $k$ stratum of $V$, and each stratum is itself a symplectic normal crossing divisor in the smooth resolution of the previous one, cf. \S \ref{S.A.nc.strat}.  Basic facts about stratifications are reviewed in \S \ref{S.A.stratif}.


\setcounter{equation}{0}
\section{Outline of the construction of $\ov\M(X, V)$}\label{s-out}

\medskip

The construction of the moduli spaces $\ov\M(X, V)$ relative normal crossing divisors $V$ is modeled on the construction of \cite{ip1} and \cite{ip2} for the case when $V$ is smooth, which in 
turn is modeled on the classical construction of the absolute moduli space $\ov\M(X)$. 

This section reviews some of the main steps involved in the construction, focusing on the global aspects of the theory, while the later sections focus on the local considerations required in making these statements precise and proving them. Roughly speaking, the construction of the moduli space and its GW invariants involves three main ingredients: Gromov compactness theorem, transversality and gluing.  

\subsection{ Brief review of the absolute moduli space.} Fix a closed, finite dimensional smooth symplectic manifold $(X,\om, J)$ with a tamed almost complex structure and let $\ov\M^J(X) $ denote the moduli space of stable $J$-holomorphic maps into $X$. To construct GW invariants one needs to consider deformations of the structure, so let $\J(X)$ be the space of smooth tamed pairs $(\om, J)$ on the target $X$, and $\ov\M(X)\ra \J(X)$ the corresponding universal moduli space of stable, pseudo-holomorphic maps  $f:C\ra X$.  For transversality reasons, one often needs to consider further extensions $\J(X)\hookrightarrow \JV(X)$ of the space of parameters  by either taking various Sobolev completions or turning on perturbations $(J, \nu)$ of the $J$-holomorphic map equation, see Remark  \ref{transv}.

As a topological space $\ov\M(X)$ is constructed in a functorial fashion starting with the moduli space $\M_{A, g, n}(X)$ of $J$-holomorphic maps $f:C\ra X$ with smooth stable  genus $g$ domains with $n$ marked points and representing $A=f_*[C]\in H_2(X)$. A nodal curve $C$ is called stable if $\Aut C$ is finite, and a smooth map $f:C\ra X$ is stable if $\Aut(f, C)$ is finite. 

One uses Gromov compactness Theorem to define a compactification $\ov\M_{A, g, n}(X) $ consisting of stable maps with nodal domains, together with a stabilization map $\st:  \ov\M_{A, g, n}(X) \ra \oM_{g, n}$ to the Deligne-Mumford moduli space. Points of $\oM(X)$ are potential limits of sequences $f_n:C_n\ra X$ of $J_n$-holomorphic maps  (where $J_n\ra J_0$ in $\J(X)$ or at least in $C^0$) obtained by (a) first passing to the limit $C$ of the domains in $\ov\M_{g, n}$ and (b) inductively rescaling  the domains to prevent the energy of the maps from concentrating. The process terminates in finitely many steps (by energy considerations) extracting a subsequence still denoted $f_n$ with a stable map limit $f_0:C_0\ra X$. Here $f_0 $ is an equivalence class of stable maps up to reparametrizations of the domain, $C$ is the stable model of $C_0$, and the convergence $f_n\ra f_0$ is in Hausdorff distance, energy density and uniform convergence on compacts away from the nodes.  

 This perspective uses the functorial description of the Deligne-Mumford moduli space $\ov\M_{g, n}$  as the classifying space of deformations of the domains  including its universal curve $\ov\U\ra \ov\M$ and its semi-local models (i.e. local models around the nodes and local trivializations on away from them). 

With the Gromov topology, all the standard maps entering in GW theory are continuous, including the projection $\pr:\ov\M(X)\ra \J(X)$, the stabilization-evaluation map
\bear\label{M-abs-se}
\st \ti \ev: \ov\M_{A, g, n}(X)\ra \ov\M_{g, n} \ti X^n
\eear
and the map $\pi: \ov\M_{A, g, n+1}(X)\ra \ov\M_{A, g, n}(X)$ that forgets one of the marked points. 
\smallskip

The universal moduli space also comes with a stratification $\oM(X)=\cup_t\M_t(X)$
 by the topological type $t$ of the domains, with the top stratum $\M(X)$ consisting of maps $f:C \ra X$ with smooth  domains. When $C$ is nodal, separating its nodes gives a lift $\wt f:\wt C\ra X$ of $f$ to the smooth resolution $\wt C\ra C$. Ordering the extra marked points of $\wt C$  defines a finite resolution  $\wt \M_t(X)\ra \M_t(X)$ of each open stratum of the moduli space, with an action of the symmetric group reordering the extra data. It also identifies the resolution $\wt \M_t(X)$ with a subset of the top stratum $\M(X)$ of another universal moduli space of stable maps $\wt f:\wt C\ra X$ with smooth, but possibly disconnected domains. The maps which descend to $C$ are precisely those which match at the nodes, i.e. those in the inverse image $\ev^{-1}(\De)$ of the diagonal under the evaluation map at the pairs of marked points corresponding to the nodes of $C$. 

Transversality, assuming it can be achieved stratawise over the parameter space $\J(X)$, implies that each open stratum of the universal moduli space has a smooth resolution, and $\ev$  is a submersion. The gluing theorem describes families of deformations encoding how the strata fit together, showing that for generic parameter the moduli space $\ov\M(X)$ has a stratified smooth resolution, and that in fact the boundary divisor is a normal crossing divisor in the smooth/orbifold/stratified spaces category. The pushforward by (\ref{M-abs-se}) of the fundamental cycle of this resolution is the GW-cycle, which is invariant by the usual cobordism argument. One can avoid gluing altogether in defining the GW-cycle whenever all  boundary strata are codimension two, which is the approach we take here. A separate paper  \cite{i-nor-sum} in preparation focuses on gluing in the context of the generalized symplectic sum formula. 
\begin{rem} \label{R.top.T.strat} The moduli space $\oM(X)$ has {\em many} types of stratifications that enter in both Gromov compactness (via a limiting argument) and gluing (via a deformation argument). Each one is described by a upper semi-continuous (USC) map 
\bear\label{D.stra.tau}
\tau: \oM(X)\ra {\mathcal T}
\eear
to some partially ordered set $\mathcal T$,  whose level sets $\M_t(X)=\tau^{-1}(t)$ are by definition the open strata where the topological information $t$ is fixed, see  \S \ref{S.A.stratif}.  Examples include (i) $\tau(f)=\#$ nodes of the domain (ii) $\tau(f)=(A, g, n)$  or (iii)  $\tau(f)=$ dual graph associated to $f$.  Note that (ii) is continuous, while the other two are only USC, and that (iii) is a refinement of the other two, see \S \ref{S.A.stratif}. 

In each instance the meaning of $\tau$ and $\mathcal T$ depends on what topological information is relevant at a given moment in the analysis arguments, which we expand or contract as needed.
\end{rem} 

\subsection{The construction of the relative moduli spaces.} Assume next that $V$ is a symplectic normal crossing divisor $\iota:\wt V \looparrowright X$ in $(X,\omega)$ and restrict to the space $\J(X, V)$ of tamed $V$-compatible pairs $(\om, J)$. Then the absolute moduli space $\ov\M(X)$ has a further stratification indexed by the topological type of the intersection with  $V$. The level zero stratum $\M(X, V)$ of the relative moduli space consists of stable maps into $X$ without any components or nodes in $V$, and further decorated by their contact information to $V$. 

Next comes the task of proving a refined version of Gromov compactness for $\M(X, V)$, obtained by rescaling to prevent energy from concentrating both in the domain (giving rise to unstable domain components) and also in the target along $V$ (giving rise to a level $m$ building $X_m$).  As before, this refined notion of convergence {\em defines} $\ov\M(X, V)$ as a topological space (with a Hausdorff topology) in a functorial fashion, together with the corresponding notion of stability,  refined evaluation map (involved in the matching conditions), and refined stratifications $\si:\oM(X, V) \ra \s$ by topological type. Its proof occupies most of the remaining part of the paper. 

To outline some of the ingredients, recall that a nodal curve $C$ comes with marked points and singular points (all together special points), smooth resolutions and smooth deformations, which enter in the description of the absolute moduli space $\oM(X)$.  The resolutions of $C$ give rise to resolutions of open strata of $\oM(X)$, while deformations of $C$ smoothing the nodes encode the normal directions to the strata. When $C$ has unstable rational components, there is also a map collapsing them. As we shall see, this  resolution/deformation/collapse picture also extends to the target for a level $m$ building and in fact to each stratum of the relative moduli space, becoming responsible for its structure.  


Level $m$ buildings $X_m$ associated to a normal crossing divisor $V$ in $X$ are described in Section \ref{s3} and here are some of their key features. They come with a zero divisor $V_m$ and a singular divisor $W_m$ which combine into a total divisor $D_m$. The building is defined in terms of a smooth resolution 
$(\wt X_m, \wt D_m)$ via an attaching map $\xi$ that attaches together components of the symplectic normal crossing divisor $\wt D_m$ to produce the singular divisor $W_m$ of the building $X_m$. The  resolution $\wt X_m$  has several levels, with $(X, V)$ on the level zero and a fibration in each positive level. For each $l\ge 1$, there is a map $p_l: (X_m, D_m) \ra (X_{m-1}, D_{m-1})$
that collapses down that level, and their composition $p:X_m\ra X$ collapses all positive levels: 
\bear\label{O.attach}
(\wt X_m, \wt D_m) \ma\longra^\xi (X_m, D_m)\ma\longra^{p}
 (X, V)
\eear

A level one building $X_1$ also comes with families of deformations ${\mathcal X}\ra B$ over a small disk $B$ whose central fiber is $X_1$ and the rest of the fibers $X_\la$ are smooth, and  can be regarded as a symplectic sum of $X_1$ along the singular locus with gluing parameter $\la$. More generally, 
a  level $m$ building $X_m$ comes with families of symplectic deformations ${\mathcal X}\ra B$ that separately smooth out the singular locus in level $l$ which combine into families of deformations
\bear\label{O.defm.x} 
\xymatrix{  
X_\la\ar@{^(->}[r]^{\iota_\la} &{\mathcal X} \ar[r]^\pi & B 
}
\eear
over a product of small disks. The center fiber is $X_m$ while the fiber $X_\la$ over a point $\la=(\la_1,\dots ,\la_m)$ can be regarded as the sum of $X_m$ along the singular divisor using gluing parameters $\la$. The fiber over a point with $k$ nonzero coordinates is a level $m-k$ building.  

Implicit in this description is the fact that all diagrams involve not just smooth objects (resolutions, quotients and deformations), but also symplectic ones, each one of them coming  with various spaces of tamed pairs of parameters $(\om, J)$ compatible with the extra structure. 

\medskip

The relative moduli space $\ov\M(X, V)\ra \J(X, V)$ consists of equivalence classes of $J$-holomorphic maps $f:C\ra X_m$ into a level $m$ building satisfying certain matching conditions along the singular locus. To fully describe the moduli spaces $\ov\M(X, V)$ we need to consider spaces of $J$-holomorphic maps $f:C\ra X_m$ into a level $m$ building, with their resolutions $\wt f:\wt C\ra \wt X_m$ and deformations $f_{\mu, \la}: C_\mu \ra X_\la$ as well as the refined evaluation map. 
\begin{rem} {\bf (Stability and automorphisms)} \label{nu.univ} To simplify the discussion in this paper, in all local analysis arguments we assume that all nontrivial components of the domains $C$ are already stable, and were further decorated to kill their automorphism groups. For Gromov compactness-type arguments it suffices to add extra marked points. For transversality arguments, if the domains are already stable, one can pass to a finite cover of the Deligne-Mumford moduli space $\ov\M$, using instead the space $\ov\M^G\ra \ov\M$ of twisted $G$-covers for a suitably chosen finite group $G$ (cf. \cite{ip-don}).

At first the assumption that all (nontrivial) domain components  have trivial automorphism group seems to be a very restrictive assumption, but it can always be achieved  whenever one of the branches of $V$ is a sufficiently positive Donaldson divisor. This not only provides a global way to stabilize all domains, but simultaneously simplifies the analysis required in proving transversality and gluing (at the expense of complicating the topology and combinatorics). Moreover, the general case follows  by functorial considerations from this seemingly very special case (cf. \cite{ip-don}). 
 \end{rem}
\begin{rem} \label{transv} In this paper, we only work with Gromov-type perturbations $(J, \nu)\in\JV$ of the $J$-holomorphic map equation:
\bear\label{eq.f=nu}
 \ov\bd_{jJ} f (z) =\nu(z, f(z)) \quad  \iff \quad \ov\bd_{j J_\nu} F=0 \quad \mbox{ where }\hskip.4in
\\
\label{graph}
F:C\ra\ov \U\ti X\quad \mbox { is the graph $F(z)=(z, f(z))$ of  $f:C\ra X$} 
\eear
These are global and functorial perturbations induced by deformations $J_\nu$ of the product almost complex structure $j\ti J$ on $\ov \U \ti X\subset  \P^N\ti X$, where  $ \ov \U\hookrightarrow \P^N$ is a fixed embedding of the universal curve $\ov \U\ra \ov\M$. The projection $(J, \nu)\mapsto J$ defines a fibration with fiber $\nu\in {\mathcal V}=\Hom^{0,1}(T\ov\U, TX)$ and a zero section:
\bear\label{forget.nu}
\xymatrix{
{\mathcal V}\ar@{^(->}[r]&  \JV\ar[r]&\ar@ /_1pc/ [l]   \J.
}
\eear
When the domain $C$ has trivial automorphisms the graph $F$ is an embedding, thus  somewhere injective, which is often sufficient for the standard transversality argument in \cite{ms2} to extend.  On the other hand, Gromov-type perturbations vanish on all unstable components of $C$ (thus on trivial components), so these are always $J$-holomorphic. Under the assumption in Remark \ref{nu.univ}, Gromov-type perturbations provide enough transversality at nontrivial components to be able to deal with the potential lack of transversality at the trivial ones. 


In the context of relative moduli spaces, we work only with $V$-compatible parameters $(J, \nu)\in \JV(X, V) $ for the equation (\ref{eq.f=nu}), see  \S \ref{A.spaces.param}. If we let $(X', V')=(\ov\U\ti X,\; \ov\U\ti V)$, this condition is equivalent to requiring that the almost complex structure $J_\nu$ on $X'$ be compatible with  the divisor $V'$ to first order. The situation extends to a level $m$ building, where the corresponding almost complex structure $J_\nu$ on $X'_m$ is required to be compatible with  the total divisor, and to its deformations ${\mathcal X}\ra B$ where the almost complex structure is required to be compatible with the fibration. So in most of the arguments below we switch back and forth between $f$ and its graph $F$, converting statements about a $(J, \nu)$-holomorphic map $f:C\ra X$  into statements about its $J_\nu$-holomorphic graph 
$F:C\ra \ov\U\ti X$. 
\end{rem} 

\setcounter{equation}{0}
\section{The Relative Moduli Space $\M(X,V)$}\label{s2}

\medskip

Assume $(X,\omega, J)$ is a smooth symplectic manifold with a normal crossing divisor $V$, and restrict to the subspace $\J(X,V)$ of $V$-compatible almost complex structures on $X$. For transversality purposes we also need to turn on Gromov-type perturbations $(J, \nu)\in \JV(X, V)$. 

The definition of $\oM(X, V)$ takes several stages. In this section we describe the main piece
\bear\label{D>M.S}
\M(X,V)\ra \JV(X, V)
\eear
consisting of stable $(J, \nu)$-holomorphic maps $f:C\ra X$ into $X$ without any components or nodes in $V$, such that all the points in $f^{-1}(V)$ are marked, and each marked point $x$ of $C$ comes decorated by a sequence of multiplicities recording the order of contact of $f$ at $x$ to each local branch of $V$. There is a forgetful map 
\bear\label{ft.V.1}
\ft_V: \M(X, V) \ra \oM(X)
\eear 
that forgets the contact information to $V$ and whose image is the open subset 
\bear\label{def.U.mod}
U = \{ f\in  \oM(X)\;|\; \text{ $f$ has no components and no nodes in $V$ }\}.
\eear
Here the domain $C$ of $f$ could be nodal, as long as all its nodes are off $V$, and could be disconnected, as long as none of its irreducible components is mapped entirely in $V$.  We also define the leading order section (\ref{sect.lead}) and use it to construct a refined evaluation map on $\M(X, V)$. 
\smallskip

Assume $f:C\ra X$ is any map in \eqref{def.U.mod}. Then $f^{-1}(V)$ consists of finitely many smooth  points of $C$ which we can mark and decorate by their local contact information to $V$ as follows. 

Fix a point $x\in f^{-1}(V)$ such that  $p=f(x)$ is a depth $k\ge 1$ point of $V$. Choose local coordinates $z$ about the point $x$ in the domain, and index the $k$ branches of $V$ meeting at $f(x)$ by the set $I(x)$ as in  Remark  \ref{N.strata}.  For each $i\in I(x)$, choose a local coordinate at $f(x)$ in the normal bundle to the branch labeled by $i$, see (\ref{nor.vk}). Lemma 3.4 of \cite{ip1}  implies that the normal component $f_i$ of the  map $f$  around $z=0$ has an expansion: 
\bear\label{f.exp}
f_i(z)=a_i z^{s_i}+O(|z|^{s_i})
\eear
where $s_i$ is a positive integer and  the leading coefficient $a_i \ne 0$. The multiplicity $s_i$ is independent of the local coordinates used, and records the order of contact of $f$ at $x$ to the $i$'th  local branch of $V$. 

Thus each point $x\in f^{-1}(V)$ comes with the following contact information: 
\begin{enumerate}[(a)]
\item a depth $k(x)\ge 1$ that records the codimension in $X$ of the open stratum of $V$ containing  $f(x)$; 
\item an indexing set $I(x)$ of length $k(x)$ that keeps track of the local branches of $V$ meeting at $f(x)$; intrinsically $I(x)$ is the fiber $\iota^{-1}(f(x))$ of the immersion $\iota:\wt V \looparrowright X$ at the point  $f(x)$; 
\item  a sequence of positive multiplicities $s(x)=(s_i(x))_{i\in I(x)}$, recording the order of contact of $f$ at $x$ to each local branch of $V$ at $f(x)$. 
\end{enumerate}
\begin{rem}\label{s.empty} There are several ways to encode such discrete topological information. One way is via a map $s(x): I(x) \ra \N_+$ defined by $i\mapsto  s_i(x)$, in which case $I(x)$ is its domain and $k(x)=|I(x)|$. Either way, $s(x)$ keeps track of how the local intersection number of $f$ at $x$ with $V$ is {\em partitioned} into intersection numbers to each local branch of $V$ cf. \eqref{deg.s}.  The sequence $s$ is obtained by putting together the multiplicities $s(x)$ for all $x\in f^{-1}(V)$ and keeps track of all the contact information of $f$ to $V$. 

By convention, we include in $s$ the information about ordinary marked points (i.e. for which $f(x) \notin V$). Such points have depth $k(x)=0$, $I(x)=\emptyset$ and $s(x)=\emptyset$. 
\end{rem}
With this, let $R= f^{-1}(V)$ be the contact points and $P\supseteq R$ all  the marked points. Each  $x\in P$ has an associated sequence $s(x)$ which records the order of contact of the map $f$  at $x$ to each local branch of $V$, including the indexing set $I(x)$ of the branches.  The cardinality $k(x)=|I(x)|$ is the {\bf depth} of  $x$  while  the   {\bf degree} 
\bear\label{deg.s}
\deg s(x)\; \ma =^{\rm def}\; \sum_{i\in I(x)} s_i(x) 
\eear  
is the local intersection number of $f$  at $x$  with  $V$. For each {\em positive} depth point $x$, its {\bf isotropy} $\Gamma_x$ is the group of roots of unity of {\bf order}
\bear\label{multipl.s}
\al(x) \; \ma =^{\rm def}\; \mathrm{gcd} \{ s_i(x)\;|\;  i\in I(x)\}, 
\eear
the greatest common divisor of its contact multiplicities. The total degree of $s$ is purely  topological: 
\bear\label{s=av}
 \deg s=\sum_{x\in P} \deg s(x)=A\cdot V
\eear
where $A\in H_2(X)$ is the homology class of the image of $f$. Denote by $\ell(s)=|P|$ the {\bf length} of $s$, i.e. the total number of marked points, out of which $\ell_+(s)=|R|$ are contact points (i.e. with depth $k(x)>0$). Intrinsically, the depth is a map $k:P\ra \N$ which partitions the marked points into ordinary marked points $P_0$ (of depth zero) and $R^k$ the depth $k\ge 1$ contact points. 

\begin{rem}\label{R.aut.s} Combinatorial information can be described by first {\em ordering} all the data involved: (i) the $\ell(s)$ marked points $P$ of $C$, which include all contact points $R=f^{-1}(V)$  and (ii) the $k(x)$ branches $I(x)=\iota^{-1}(f(x))$ of $V$ for each $x\in f^{-1}(V)$. 
There is a free action of a wreath product of symmetric groups reordering this data, with the subgroup 
\bear\label{D.s.aut} 
G_s= S_{R}  \ltimes \ma\prod_{x\in R}  S_{I(x)}  \quad \text{ of order }\quad  
|G_s|= \ell_+(s)!\ma \prod  _{x\in R} k(x)!
\eear 
fixing the depth zero marked points. The group $\Gamma_x$ of roots of unity acts on the local coordinate $z$ on the domain, fixing all leading coefficients \eqref{f.exp} of $f$ at $x$. Note that data  \eqref{deg.s}-\eqref{D.s.aut} is intrinsically associated to the map $f:C\ra X$, independent of choices. So decorating each element of $U$ by a choice of an ordering of its contact information to $V$  {\em defines} the moduli space $\M(X, V)$ as a disjoint union of pieces $\M_s(X, V)$ indexed by the (ordered) data $s$. If $[s]$ is the equivalence class of $s$ (up to reordering), and $\M^{[s]}(X, V)$ the union of the pieces associated to the same equivalence class $[s]$, then the group $G_s$ acts on $\M^{[s]}(X, V)$. The quotient 
$\M^{[s]}(X, V)/G_s=\M_{[s]}(X, V)$ is the subset of $U$ consisting of maps whose (unordered) contact information to $V$ is  $[s]$.\end{rem}
In effect, we described a stratification $\mathsf{s}$ of both $\M(X, V)$ and $U$, with open strata indexed by the topological type of the intersection of $f$ with $V$!  But the moduli space $\oM(X)$ has other  natural stratifications (cf. \S \ref{S.A.stratif}), for example $\tau:\oM(X)\ra {\mathcal T}$ by the dual graph. When restricted to $U\subseteq \oM(X)$ or its resolution $\M(X, V)$, the dual graph stratification can be refined
\bear
\tau: \M(X, V) \ra  {\mathcal T}_{V}
\eear
to include the topological type of the intersection of $f$ with $V$ described above.  The (refined) dual graph $\tau(f)$ associated to a point  $f\in \M(X, V)$ has its half-edges (corresponding to marked points $x$ of $C$) decorated with their contact information $s(x)$ to $V$. Its vertices $v$ correspond to irreducible components $C_v$ of the domain $C$ (decorated by their genus $g_v$ and the homology class $A_v=f_*[C_v]$ of the image), while its edges correspond to the nodes $D$ of $C$. So far the contact information to $V$ decorates only the half edges of $\tau$, not its vertices or edges, because we assumed $f$ had no components or nodes in $V$. 

\begin{rem}\label{R.dual.defm}  Any nodal curve $C$ comes with a dual graph $\tau=\tau(C)$; the dual graph  associated to a map $f:C\ra X$ refines it by adding decorations to $\tau$. For each  node $x$ of $C$, we can consider either the resolution $\wt C\ra C$ resolving the node, or a family of deformations $C_\mu$ smoothing the node. The dual graph of $\wt C$ and respectively $C_\mu$ is obtained from $\tau$ by cutting the edge corresponding to $x$ in half or respectively contracting it. Forgetting a marked point $C$ removes a half-edge of $\tau$, collapsing a component of $C$ removes a vertex of $\tau$. Intrinsically, this encodes how the  stratification interacts with other natural structures (resolutions/deformations/projections) on the moduli spaces, see \S \ref{S.A.stratif}.\end{rem}
\begin{rem}\label{top.info} Fixing any (refined) dual graph $s$ therefore defines an open stratum $\M_s(X, V)=\tau^{-1}(s)$ of the moduli space consisting of maps with that dual graph, with the 'top stratum' consisting of smooth domains. From the dual graph $s$ we can read off not just the marked points $P_s$ of the domain (half edges $x$ of the graph) and their contact information $s(x)$ to $V$, including the number of depth zero (ordinary) marked points $n_s$ but also the topological type $\Si_s$  of the domain of $f$, including the indexing set $D_s$ of its nodes and its {\em Euler characteristic} $\chi_s=\chi(\Si_s\setminus D_s)$, as well as the {\em total homology class} $A_s=f_*[\Si_s]$. The domain of $f$ could be disconnected, and the dual graph $s$ includes this information. We could even record in the dual graph $s$ the relative homology class of each domain component, i.e. the information about the rim tori, see Section 5 of \cite{ip1}. The construction there is purely topological, so extends easily to the case when $V$ is a normal crossing divisor. We will not explicitly describe it in this paper. \end{rem}

For each refined dual graph $s$ we can forget the decorations associated to $V$ to get the  usual dual graph $\tau=\mathsf{ft_V}(s)$, or forget the graph structure and only record the contact information $(s(x))_{x\in P}= \mathsf{ev}(s)$ associated to its half edges $P=P(s)$. This is encoded in the following commutative diagram:
 \bear\label{M.ft.V.strat}
\xymatrix{
\M(X, V) \ar[r]^{\tau}\ar[d]_{\ft_V}&{\mathcal T}_{V}\ar[d]^{\mathsf{ft_V}} \ar[r]^{\mathsf{ev}}& \s_{V}
\\
\oM(X) \ar[r]^{\tau}&{\mathcal T}}
\eear
which we eventually extend to the compactification $\oM(X,V)$ so that the composition 
$\mathsf{s}=\ev\circ \tau$ is continuous (locally constant). 
For the rest of the paper, $s$ denotes either a refined dual graph, or only the contact information $(s(x))_{x\in P}$ associated to the marked points, depending on the context. Correspondingly, $\M_s(X, V)$ denotes the subset of the relative moduli space consisting of maps with that associated topological information, intrinsically a fiber of either $\tau$ or $\mathsf{s}$. 
 
\begin{rem} ({\bf Automorphisms and Stability II}) In this paper, by a stable object of a moduli space we always mean one whose automorphism group is finite. As we already saw, often there are several natural ways to describe the objects of the moduli space e.g. by first adding geometric decorations, and then forgetting about them; intrinsically these  describe resolutions of the moduli spaces, which affect the automorphism groups (e.g. whether an automorphism fixes the extra decorations or can permute them). 

For example, an element of $\oM(X)$ is an equivalence class $[f]$ or more precisely $[f, C, J, \nu, \om]$ of stable $(J, \nu)$-holomorphic maps $f:C\ra X$ up to reparametrizations of their domain. For each fixed domain $C$, the automorphism group $\Aut(f)$ of $f:C\ra X$ is a subgroup of $\Aut C$ and $[f]$ is stable iff $\Aut(f)$ is finite.  When the domain $C$ is disconnected, we have a {\em choice} of whether we order or not its components which may affect $\Aut(C)$, but does not affect the stability of $C$.  Moreover, the stability of $[f]\in \oM(X)$ is equivalent to a topological condition which can be read off its dual graph $\tau(f)$ (see \S \ref{S.maps.resol}), so we may restrict to 'stable dual graphs'.  Any automorphism of $f$ induces an automorphism of its dual graph. 
\end{rem} 
Similarly, an element of $\M(X, V)$ is an equivalence class up to reparameterizations of the domain, but now of {\em decorated} maps $f:C\ra X$: here $f$ is a stable $(J, \nu)$-holomorphic map without any components or nodes in $V$, such that $R=f^{-1}(V)$ is a subset of the marked points $P$ of the domain, and for each $x\in P$ the $k(x)$ branches $I(x)$ of $V$ over $f(x)$ are {\em ordered}, and decorated by the contact multiplicities $s(x)=(s_i(x))_{i\in I(x)}$ of $f$ to $V$ at $x$. Because both $P$ and $I(x)$ are ordered, we can talk about the $k$'th marked point or the $i$'th local branch of $V$ at it,  despite the fact that in $X$ the local branches of $V$ may globally intertwine. 

There is an action of the product of symmetric groups $G=\prod_{x\in P} S_{I(x)}$ reordering this information, inducing a natural action on a subset of the moduli space $\M(X, V)$. The quotient $\M(X, V)/G$ consists of equivalence classes whose marked  points remain ordered, but the local branches are {\em unordered}. The further quotient by $S_R$ consists of equivalence classes whose depth zero marked points $P_0$ remain ordered, but the contact multiplicities to $V$ are unordered, etc.
\begin{rem} \label{R.s.intrinsic} Equivalently, for each map $f:C\ra X$ we can regard the marked points of the domain as indexed (parametrized) by a fixed ordered set $P$ (via an injection $P \hookrightarrow C\setminus\{ nodes\}$), and the local branches $\iota^{-1}(f(x))$ of $V$ at $f(x)$ indexed  by fixed ordered set $I(x)$, one for each $x\in P$. Each element $f\in \M(X, V)$ then comes decorated by a {\em choice} of a bijection 
\bear\label{D.rho}
\rho_f(x):I(x) \ra \iota^{-1}(f(x)),
\eear 
one for each contact point $x\in f^{-1}(V)$, where $\iota:\wt V \looparrowright X$ is the immersion parameterizing $V$.  Forgetting $\rho_f(x)$ corresponds to forgetting the order of the branches, or equivalently taking the quotient by $S_{I(x)}$. Evaluating $\rho_f(x)$ at $i\in I(x)$ gives a point in the resolution $\wt V \ra V$ over $f(x)$.   
\end{rem}


\smallskip

The {\em ordinary evaluation map} $\ev_x:\M_s(X,V)\ra V_{k(x)}$ at a marked point $x\in P(s)$ is defined by $\ev_x(f)=f(x)$, and comes with a natural lift 
\bear\label{ev.naive}
\ev_x:\M_s(X,V)\ra V_{I(x)}
\eear
to the resolution $V_{I(x)}$ keeping track of the indexing \eqref{D.rho} of the local branches of $V$ at $f(x)$.  Here $V_{k}$, perhaps more appropriately denoted $(X,V)_k$, is the depth $k$ stratum in the stratification of $X$ induced by the normal crossing divisor $V$, and $V_I$ is its resolution \eqref{A.D.V.I}.  As before, this evaluation map includes the ordinary marked points with the convention $V_{\emptyset}=X$. Evaluating simultaneously at all the marked points $x\in P(s)$ associated to $s$ gives
\bear\label{or.ev}
\ev:\M_s(X,V)&\ra& V_{s} \ma=^{\mathrm {def}}  \ma \prod_{x\in P(s)}  V_{I(x)}. 
\eear
It is important to note the target of this evaluation map: there are several other choices that may seem possible (e.g.  $X^{\ell(s)}$), but this is the only choice for which the evaluation map can be a submersion, without loosing important information.  

When the depth $k(x)\ge 2$, the evaluation map (\ref{ev.naive}) does not record enough information to state the full matching conditions appearing in the relatively stable map compactification. We also need to record the leading coefficient of the expansion (\ref{f.exp}). For each $f:C\ra X$ in $\M_s(X,V)$ and each $x\in f^{-1}(V)$ let  $a_i(x)$ the   leading  coefficient (\ref{f.exp}) of $f$ at $x$ in the normal direction $N_i$ to the branch labeled by $i\in I(x)$. As explained in Section 5 of \cite{ip2},  $a_i(x)$  is naturally an element of $(N_i)_{f(x)}\otimes (T_x^*C)^{s_i(x)}$ so it  defines a section of the bundle
\bear
\ev_x^*N_i\otimes \L_x^{s_i(x)}  
\eear
where $\L_x$ is the relative cotangent bundle to the domain at the marked point $x$. If we denote by $E_{s,x}\ma= \ev_x^* N V_{I(x)}\otimes_{s(x)}  \L_x$ the bundle over the moduli space whose fiber at $f$ is  
\bear
\ma\bigoplus_{i\in I(x)}  \ev_x^* N_i\otimes \L_x^{s_i(x)}
\eear
then the {\em leading order section} at $x$ is defined by  
\bear\label{sect.lead}
\si_x:\M_s(X,V)&\ra& E_{s,x}= \ev_x^* N V\otimes_{s(x)}  \L_x
\\
\nonumber
\si_x(f) &=& (a_i(x))_{i\in I(x)}
\eear
It records the leading coefficients in all normal directions, a crucial information needed in later arguments. This section was already used in \cite{i2} to essentially get an isomorphism between the relative cotangent bundle of  the domain and that of the target for the moduli space of branched covers of $\P^1$. 


The {\em refined evaluation map} at the marked point $x$ is defined by 
\bear\label{enh.ev} 
\Ev_x:\M_s(X,V)&\ra&\P_{s(x)} (N {V_{I(x)}})
\\
\nonumber
\Ev_x(f) &=&[\si_x(f)]
\eear
where the $\P_{s(x)} (N{V_{I(x)}}) $ is the weighted projectivization with weight $s(x)$ of  the normal bundle $N {V_{I(x)}}$ of the depth $k(x)$ stratum $V_{I(x)}$.  More precisely, 
$\P_{s(x)} (N {V_{I(x)}})$ is  a bundle over $V_{I(x)}$ whose fiber is the weighted projective space obtained as the quotient  by the $\cx^*$ action with weight $s_i(x)$ in the normal direction $N_i$ to the branch labeled by $i\in I(x)$.

Of course, if $\pi:\P_{s(x)} (N{V_{I(x)}}) \ra V_{I(x)}$ is the projection then 
\best
\pi\circ \Ev_x=\ev_x
\eest 
which explains the name; $\Ev_x$  is a nontrivial refinement of $\ev_x$ only when the depth $k(x)\ge 2$. \begin{rem}
By construction the leading order terms are nonzero, so the image of $\si_x$ is away from the zero  sections of each term in (\ref{sect.lead}). The image of the refined evaluation map \eqref{enh.ev} similarly lands in the complement of all the coordinate hyperplanes, where $\P_{s(x)}(NV_{I(x)})$ is smooth, though perhaps not reduced: the isotropy subgroup is $\Gamma_x\le \cx^*$, the roots of unity of order \eqref{multipl.s}. \end{rem}

\begin{rem}\label{R.form.Ev} Note that \eqref{sect.lead} and \eqref{enh.ev} only depend on the restriction of $f$ to an infinitesimal neighborhood of $V$ in $X$, and not on the rest of the dual graph of $f$, while \eqref{or.ev} also includes the value of $f$ at the depth zero marked points. We formally include depth zero marked points $x$ by setting $\Ev_x=\si_x=\ev_x$, $E_{s, x}= \L_x$ and $(\P^1)^{\emptyset}=\P_{\emptyset}(0)$= pt. 
\end{rem}
The combined refined evaluation map at all the marked points (including depth zero points) is 
\bear\label{enh.ev.full} 
\Ev:\M_s(X,V)&\ra& \P_s(NV) \ma=^{\mathrm{def}}\prod_{x\in P(s)} \P_{s(x)} (N V_{I(x)}) 
\eear
We next show that the contact information $s$ to $V$ is locally constant on $\M(X, V)$, see \eqref{M.ft.V.strat}. So the strata of $\M(X, V)$ corresponding to different $s$'s are {\em disjoint}, which is not true on \eqref{def.U.mod}.  

\begin{lemma} Consider a sequence $\{ f_n\}$ of maps in $\M_s(X,V)$  and assume its  limit $f$ in the usual stable map compactification has no components in $V$. Then $f \in \M_s(X,V)$.
\end{lemma}
\begin{proof}  A priori, there are two reasons why $f$ would fail to be in $\M_s(X, V)$: 
\begin{enumerate}
\item[(a)] $f$ has a node in  $V$ or
\item[(b)] the contact information of $f$ to $V$ is not given by $s$.
\end{enumerate} 
Note that case (b) includes the cases when in the limit  the multiplicity of intersection jumps up or when a depth $k$ marked point falls into a higher depth stratum of $V$. 

Since ALL points in $f_n^{-1}(V)$ are already marked, indexed by the same set $P$, they persist as marked points for the limit $f$, in particular they are distinct from each other and from the nodes of $f$.  On the other hand, let  $\wt f:\wt C\ra X$ denote the lift of $f$ to the smooth resolution of its domain. Since  $f$ has no components in $V$, then each point in  
$\wt f^{-1}(V)$  has a well defined sequence $s_0$ that records the local multiplicity of intersection of $\wt f$ at that point with each local branch of $V$. At those points of $\wt f^{-1}(V)$  which were limits of the marked points in $f_n^{-1}(V)$, the multiplicity $s_0(x) \ge s(x)$, as the multiplicity could go up when  the leading coefficients converge to 0.  But then
\best
[f]\cdot V =\sum_{x\in \wt f^{-1}(V)} s_0(x) \ge \sum_{x\in P} s_0(x) \ge \sum_{x\in P} s(x) = [f_n]\cdot V
\eest
Since  $[f_n]=[f]$ then  $\wt f^{-1}(V)=P$, which means that $f$ has no nodes in $V$, ruling out (a), and that $s_0(x)= s(x)$ for all $x\in P$, which  rules out (b). \end{proof}  

Next assume $s$ is a  dual graph as in Remark \ref{top.info}. Lemma 4.2 in \cite{ip1} easily extends to this context: 
\begin{lemma} Each open stratum $\M_s(X, V)\ra \JV(X, V)$ of the universal moduli space  is cut transversally and the refined evaluation map \eqref{enh.ev.full} is a submersion at all points $(f, C, J,\nu)$ with $\Aut \;C=1$ (or more generally whose graph $F:C\ra \ov\U \ti X$ is somewhere injective). 

In particular, for generic $V$-compatible parameter $(J, \nu)$, each stratum is a smooth manifold of dimension 
\bear\label{dim}
\dim  \M_s(X, V)=  2c_1(TX)A_s+({\rm dim} X-6)\frac{\chi_s}2 +2\ell(s)-2A_s\cdot V-2 |D_s|
\eear
near such points. 
\end{lemma}
\begin{proof} As usual, one proves this first for the stratum with smooth domains and then extends it to the nodal stratum by passing to the resolution and using transversality of the evaluation map to the diagonal, as in the  proof of Lemma \ref{L.d=d-2}. Note that any $f\in \M(X, V)$ has all  nodes off $V$ (depth zero), and at  such points $\Ev=\ev$. 

For the dimension count, when the domain is smooth, $\dim \M(X)$ is given by the first three terms in \eqref{dim}. By Lemma 4.2 of \cite{ip1}, the fact that $f\in \M(X)$ has a contact of order $s_i(x)$ at a point $x$ to the $i$'th branch of $V$ imposes a codimension $2s_i(x)$ condition on $f$, all together adding to $2\deg s=2A_s\cdot V$ so 
\best
\dim\M_s(X, V)=\dim \M(X)-2\deg s  
\eest
giving \eqref{dim} (when the domain has no nodes). 

When the domain of $f\in \M_s(X, V)$ has nodes, pass first to the resolution $\wt f:C \ra X$, whose dual graph $\wt s$ is obtained from $s$ by cutting each edge $x\in D_s$ in half. Since $f$ is continuous at the nodes, the resolution $\wt f$ is in the inverse image of the diagonal $\De$ under the evaluation map $\ev_{D_s}: \M_{\wt s}(X, V)\ra (X\ti X)^{|D_s|}$ at the pairs of marked points of $\wt C$ corresponding to the nodes $D$ of $C$. This imposes $|D_s|\cdot \dim X$ conditions on $\wt f\in \M_{\wt s}(X, V)$ so 
\best
\dim \M_s(X,V)=\dim \M_{\wt s}(X,V)- |D_s|\cdot \dim X
\eest 
Since $\ell(\wt s)= \ell (s)+ 2|D_s|$, $\chi_{\wt s}= \chi_s+2|D_s|$ and $A_{\wt s}=A_s$, and  we already proved \eqref{dim} for $\wt s$, this proves it for $s$ as well. 
\end{proof}
\begin{rem} Note that the dimension  \eqref{dim} of a stratum only depends on $A, \chi$, $\ell(s)$ and the number of nodes, but not on the rest of the dual graph, nor on the way the intersection $A\cdot V$ is partitioned into local intersection multiplicities $s(x)$ at each marked point $x$. The dimension of a top stratum is 
\bear\label{dim.M.s}
\dim\M_{A, \chi, s}(X,V)= 2 c_1(TX)A + (\dim X-6)\frac {\chi} 2 +  2\ell(s) - 2A\cdot V
\eear
while the codimension of each boundary stratum is twice the number of (depth zero) nodes. 
\end{rem}
We conclude this section by describing how the stratification of $(X, V)$ and its resolutions induce natural ones on the entire absolute moduli space $\oM(X)$, by keeping track of how each map intersects each stratum of $(X, V)$.  

Fix a map $f:C\ra X$ in $\oM(X)$ (for a $V$-compatible parameter).  Any point $x\in C$ has a depth $k(x)$ associated to its image $f(x)\in V$. Each component $\Si$ of $C$  has a depth defined by
\bear\label{D.depth.comp.si}
k(\Si)= \max \{ \;k \;| \;f(\Si)\subseteq V^k \;\} = \min\{ \;k(x)  \;| \;x\in \Si \;\}
\eear
and let $C^k$ be union of all its depth $k$ components.  The {\em contact set} of a component $\Si$ is its  subset of points whose depth is strictly higher than $k(\Si)$. The contact set $R$ of $C$ is the union of all the contact sets of its components. 
Then  $f$ has no components in $V$ iff all its components have depth 0; in this case  $R=f^{-1}(V)$. If also $f$ has no positive depth nodes, then $C$ is smooth near $R$ so $f$ lifts to an element of $\M_s(X, V)$ obtained by decorating $f$ with its contact information to $V$ (the lift is unique up to reordering the decorations, see Remark \ref{R.aut.s}). 

If  $f:C\ra X$ has no components in $V$ but has nodes in $V$, these can also be decorated by their contact information to $V$, except that now we get two possibly different sequences of multiplicities at each node.  Specifically, let $\iota:\wt C\ra C$ denote the resolution of $C$ obtained by resolving all its positive depth nodes; each such node $x$ of $C$ lifts to two points $x_\pm\in \wt C$. The lift $\wt f:\wt C\ra X$ no longer has  nodes in $V$, thus has a well defined order of contact $s$ to $V$. So  each node $x$ of $C$ has an associated indexing set $I(x)$ of the branches of $V$ over $f(x)$ and two sequences of multiplicities $s(x_\pm):I(x) \ra \N$, one for each local branch $x_\pm$ of $C$ at $x$. Such a lift $\wt f\in \M_{\wt s}(X, V)$ is essentially unique (i.e. up to reordering). This discussion extends to give:
\begin{lemma}\label{f.dec.depth.k} Fix a  $V$-compatible parameter $(J, \nu)\in\JV(X, V)\subseteq \JV(X)$. Then any element  $f:C\ra X$ in the absolute moduli space  $\ov\M(X)$ can be canonically decomposed into depth 
$k\ge 0$ pieces $f^k: C^k\ra V^k$ which have a well defined order of contact to $V^{k+1}$ and: 
\bear\label{dec.f.k.V.0}
f^{-1}( V^k\setminus V^{k+1}) = (C^k \setminus R^{>k}) \cup R^k
\eear
where  $R^k$, $R^{>k}$ denotes the collection of depth $k$ and respectively $>k$ contact points of $C$. Moreover, 
each $f_k:C^k\ra V^k$ has a lift  $\wt f^k: \wt C^k \ra \wt {V^k}$  (unique up to reordering) with the following properties: 
\begin{enumerate}[(i)]
\item  $\wt C^k$ is a resolution of $C^k$ obtained by resolving its higher depth contact 
points $R^{>k}$;
\item $\wt f^k \in \M_{s^k}(\wt {V^k}, \wt V^{k+1}) $; 
\end{enumerate}
\end{lemma}
\begin{proof} The first part is immediate from the  definition of depth. Next, if  $f:C\ra X$ is any $(J, \nu)$-holomorphic map, then its restriction $f|_{\Si}$ to  each depth $k$ component $\Si$ has a lift $\wt f:\Si \ra \wt {V^k}$  to the smooth resolution $\wt{V^k}$ of $V^k$. This lift then has a well defined contact information \eqref{f.exp} to the higher depth strata, i.e. to the  normal crossing divisor $\wt V^{k+1}$ of $\wt {V^k}$.   \end{proof}
 Intrinsically, the first part of Lemma \ref{f.dec.depth.k} describes a stratification of the restriction of the absolute moduli space $\ov\M(X)$ to $\J(X, V)$ which keeps track of how the curves meet the divisor $V$. The second part describes a resolution of these strata as a subset of relative moduli spaces $\sqcup_k\; \M_{s^k}(\wt {V^k}, \wt V^{k+1})$ of maps into the smooth resolution $(\wt {V^k}, \wt V^{k+1})$ of each stratum of $(X, V)$. 

\setcounter{equation}{0}
\section{Rescaling the target} \label{s3}
\medskip

Assume next that  $f_n$ is a sequence of maps in $\M_s(X, V)$ whose limit $f$ in the stable maps compactification has some components in $V$. We use the methods of \cite{ip1} to rescale the target normal to $V$ to prevent this from happening. In this section we describe the effect of rescaling on the target $X$. 

\subsection{\bf Brief review of the rescaling procedure in \S6 of \cite{ip1}.} 
Assume $V$ is a smooth divisor in $X$.  In local coordinates, if $x$ is a fixed local coordinate normal to $V$, rescaling $X$ by a factor of $\la\ne 0$ means we make a change of coordinates in a neighborhood of $V$ in the normal direction to $V$:
\bear\label{resc.x}
x_\la=x/\la.
\eear 
Under rescaling by an appropriate amount $\la_n$, depending on the sequence $f_n$, in the limit we will get not just a curve in $X$ (equal to the part of $f$ that did not land in $V$),   but also a  curve in the compactification of $N_V$, i.e. in 
\best
\F=\P(N_V\oplus \cx).
\eest
Here $\F$ is a $\cx\P^1$ bundle over $V$, with a zero and infinity section $V_0$ and $V_\infty$.  Under the rescaling (\ref{resc.x}), $x_\la$ can be thought instead as a coordinate on $\F$ normal to $V_0$. Let  
$y=1/x_\la$ be the corresponding  coordinate normal to $V_\infty$ inside $\F$, so that (\ref{resc.x}) becomes
\bear\label{xy=la}
x y=\la
\eear
This procedure has  the infinity section $V_\infty$ of $\F$ naturally identified with $V$ in $X$ such that furthermore their normal bundles are dual to each other, i.e. 
\bear\label{n=n*}
N_{V/X}\otimes N_{V_\infty/\F} \cong \cx
\eear
is trivial. This identification globally encodes the local equations (\ref{xy=la}) because $x, y$ are local sections of $N_{V/X}$ and $N_{V_\infty/\F}$ respectively. 
\begin{rem} Once an 
identification (\ref{n=n*}) is fixed, then for any (small) gluing parameter $\la\in \cx^*$ equation (\ref{xy=la})  is exactly the local model of the symplectic sum $X_\la$  of $X$ and $\F$ along $V=V_\infty$. Of course, topologically 
$X\#_{V=V_\infty} \F=X$. This means that an  equivalent point of view to rescaling $X$ by a factor of $\la$ normal to $V$ is to regard $X$ as a symplectic sum $X_\la$ of 
$X$ and $ \F$ with gluing parameter $\la$, with the above choice of coordinates and identifications, including (\ref{n=n*}).  The advantage of this perspective is that the rescaled manifolds $X_\la$ now fit together as part of a smooth total space 
${\mathcal X}\ra B$ as its fibers over $\la$ where they converge (in Hausdorff distance) as $\la\ra 0$ to the normal crossing divisor 
\bear\label{x.0}
X_0=X\ma \cup_{V=V_\infty} \F, 
\eear
obtained by joining $X$ to $\F$ along $V=V_\infty$. 
\end{rem}

In coordinates, denote by $U_\la$ the tubular neighborhood of $V$ in $X$ described by $ |x|\le |\la|^{1/2}$ and by $O_\la$ the complement of the tubular neighborhood of $V_\infty$ in $\F$ described  by $|y|\ge |\la |^{1/2}$. Rescaling $X$ around $V$ by $\la\ne 0$ gives rise to a manifold $X_\la$ together with a diffeomorphism 
\bear\label{R.la}
R_\la : X \ra X_\la
\eear
which is the identity outside $U_\la$ and which identifies $U_\la$ with $O_\la$ by rescaling it by a factor of $\la$, or 
equivalently via the equation (\ref{xy=la}). 
As $\la\ra 0$, $U_\la$ shrinks to $V$ inside $X$, but it expands in the rescaled version $X_\la$ to  $\F\setminus V_\infty$. So in the limit as $\la \ra 0$, the rescaled manifolds $X_\la$, with the induced almost complex structures 
$J_\la=(R_\la^{-1})^* J$ converge to the singular space $X_0$  defined by (\ref{x.0})  with an almost complex structure $J_0$ which agrees with $J$ on $X$ and is $\cx^*$-invariant on the $\F$ piece. 
\begin{rem}\label{R.ep.on.F} If $\om$ is the symplectic form on $X$ taming $J$, then as 
$\la \ra 0$ its rescaled version $(R_\la^{-1})^* \om$ converges to a singular $\om^0$ on $X_0$. The restriction of $\om^0$ to $X$ is $\om$, but its restriction to $\F$ is equal to $\pi_V^*\om_V$,  which is degenerate along the fibers of $\pi_V:\F\ra V$. However,  there exists a family $\om^\al=\pi_V^*\om_V+\al^2 \tau$ of symplectic structures on $\F$ taming $J_0$, indexed by the size $0<\al\ll1$ of the $\P^1$ fiber and described in the four displayed equations after the statement of Proposition 6.6 of \cite{ip1}. 

Conversely, starting with any such fixed choice $(\om^\al, J_0)$ on $\F$, matching along the divisor $V$ with a fixed $(\om, J)$ on $X$, Section 2 of \cite{ip2} describes how to interpolate between them in a neighborhood of the divisor to construct a tamed pair $(\om, J)$ on $\x\ra B$   which restricts to each fiber $X_\la$ as a tamed pair $(\om_\la, J_\la)$, where $R_\la ^*\om_\la$ is a small symplectic deformation of the original $\om$ (as small as needed provided $\al$ is sufficiently small). 
\end{rem}
After rescaling the sequence $\{ f_n\}$ by an appropriate $\la_n$, and passing to a subsequence, $R_{\la_n}(f_n)$ has a new limit inside ${\mathcal X}\ra B$ which is a map into $X_0$ and satisfies certain matching condition along $V_\infty =V$. In general, different components may fall into $V$ at different rates, so we need to rescale several (but finitely many) times to catch all of them, and in the limit we get maps into a building with several levels. 

\subsection{Rescaling the manifold $X$ normal to $V$}\label{S.resc.V} 

Assume now $V$ is a normal crossing divisor. We next describe the effect on the manifold $(X, \om, J)$ of rescaling it around $V$ (in a  tubular neighborhood of $V$) to obtain its refinement, a level one building.  Using our local models, we could extend the discussion above independently in each normal direction to $V$, so normal to each open stratum of ${V^k}$,  we  could rescale in $k$ independent directions. However, globally these directions may intertwine, and not be independent, so one has to be careful how to  globally patch these local pictures.  Here  is where we use the fact we have a normal bundle $N$ defined over  the resolution $\wt V$ of $V$, and rescale normal to $V$ using the $\cx^*$ action on $N$. 

\begin{rem}\label{sev.act} The $\cx^*$ action in  the complex line bundle $N$  induces in fact several different actions. The one used in this paper is the diagonal $\cx^*$ action in the normal bundle to each stratum $V^k$ of $V$.  When  the resolution of $V$ has several components we have a separate $\cx^*$ action for each component. In particular, when the divisor $V$ has simple crossings, then we also have a local  $(\cx^*)^k$ action on $X$ normal to each $V^k$;  this essentially happens only when the crossings are simple. 
\end{rem}

When we rescale once the manifold $(X, \om, J)$ in the normal direction to  $V$, on level one we get several pieces, one for each piece $V^k$ of the stratification of $V$ according to how many local branches meet there. The level zero unrescaled  piece  is still $(X, V)$, with the original $(\om, J)$. But now level one 
\bear\label{def.f.v}
\F_V=\ma\sqcup_{k\ge 1} \F_k 
\eear 
consists of several pieces $\F_k$, one for each depth $k\ge 1$. The first  piece is 
\best
\F_1=\P(N_V\oplus \cx),
\eest 
a $\P^1$ bundle over the $\wt V$, the  resolution of $V$, obtained by compactifying the normal bundle $N_V\ra \wt V$ by adding an infinity section. Similarly,  the $k$'th piece 
\bear\label{D.f.k1}
\F_k \longra \wt {V^k}
\eear
is a  $(\P^1)^k$ bundle over the resolution $\wt{V^k}$ of the closed stratum  $V^{k}$, described in more details in the Appendix. What is important here is that $\F_k$ is a bundle over a smooth manifold $\wt{V^k}$ which is obtained by separately compactifying each of the $k$ normal directions to $V$  along $V^k$, see  (\ref{nor.vk}). This means that its fiber at a point $p\in \wt{V^k} $ is 
\bear\label{F.k.nor.J}
\F_k =\ma \ti_{i\in I} \P(N_i\oplus \cx)
\eear
where $N_i$ is the normal direction to the $i$'th branch, and $I$ is an indexing set of the $k$ local branches of $V$  meeting at 
$p$.   Globally,  these $\P^1$ factors intertwine as dictated by the global monodromy of the $k$ local branches of $V$. 

Each piece $\F_k$ comes with a natural normal crossing divisor
\bear\label{def.W}
D_k=D_{k,\infty}\cup D_{k,0} \cup F_k 
\eear
obtained by considering together its zero and infinity divisors plus the fiber $F_k$ over the (inverse image of the) higher depth strata of $V^k$. The construction of the divisor  $D^k$ is described in \S \ref{S.A.nc.strat}, but here let us just mention that $D_{k,0}$ is the zero divisor in $\F_k$ where at least one of the $\P^1$ coordinates is  equal to $0$,  while the fiber divisor $F_k$ is the restriction of $\F_k$ to the stratum of $\wt{V^k}$ coming from the higher depth stratum $V^{k+1}$. 

\begin{rem}\label{cx.act} The $\cx^*$ action in the normal bundle to $V$ induces a fiberwise, diagonal $\cx^*$ action on each piece $\F_k$ with $k\ge 1$, whose restriction to $S^1$ is a Hamiltonian action. The $\cx^*$ action preserves the divisor $D_k$, but not pointwise. 
\end{rem} 
To keep notation manageable, we make the following convention: 
\bear\label{conv.f.0}
(\F_0,D_0) \ma= ^{def} (X,V)
\eear
This is consistent with our previous conventions that $V^0=X$, so  $\F_0=\P(N_{V^0/X}\oplus \cx)=X$ and $D_0=F_0=V$ is the fiber divisor over the lower stratum, as the  zero and infinity divisors are empty in this case. However, one difference that this notation obscures is the fact that while $X=\F_0$ is on level zero (unrescaled), the rest of the pieces $\F_k$ for $k\ge 1$ are all on level one (all appeared as the result of rescaling once normal to $V$). 

\begin{figure}[h]
\centering
\includegraphics[height=4cm]{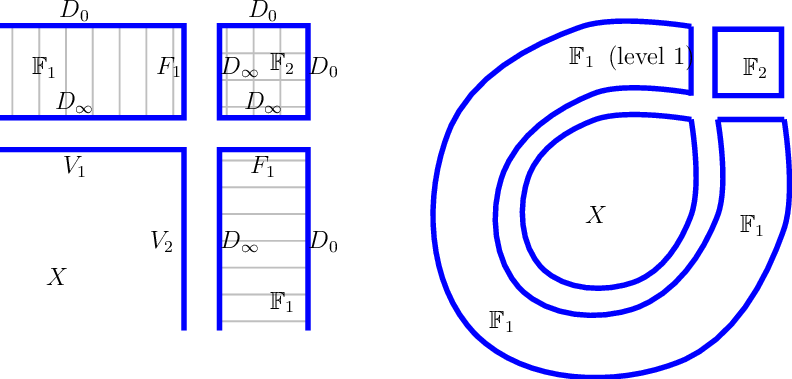}
\caption{The pieces of a level 1 building}\label{F.1}
\end{figure}
\begin{defn}\label{D.build.1} A {\em level one building} $X_1$ is obtained from  
\best
(\wt X_1, \wt D_1)= \ma\bigsqcup_{k\ge 0} \; (\F_k, \; F_k\cup D_{k, \infty} \cup D_{k,0}) 
\eest 
by  identifying the fiber divisor $F_k$ of $\F_k$ with the infinity divisor $D_{k+1,\infty}$ of  $\F_{k+1}$: 
 \bear\label{level1}
(X_1,D_1)=\ma\bigcup_{F_{k}=D_{k+1,\infty}}  (\F_k, \; F_k\cup D_{k, \infty} \cup D_{k,0})
\eear
Denote by  $W_1$ the {\em singular locus} of $X_1$ where all the pieces are attached to each other, by 
$V_1$ the {\em zero divisor} of $X_1$ obtained from the  zero divisors, and by $D_1= W_1\cup V_1$ {\em the total divisor}. A level one building comes with a {\em resolution}  $(\wt X_1, \wt D_1)$  and 
{\em an attaching map} 
\bear\label{at.map}
\xi: (\wt X_1, \wt V_1) \ra (X_1, V_1).
\eear 
It also comes with  a {\em collapsing map} to level 
zero $(X_0,V_0)=(X,V)$:
\bear\label{def.collapse}
p:(X_1, V_1)\ra (X,V)
\eear
which is identity on level 0, but collapses the fiber of each level one piece  $\F_k$, $k\ge 1$. 
\end{defn}
\begin{rem} The precise identifications required to construct this building are described in \S \ref{S.A.nc.strat}. It is easy to see that  both $F_{k}$ and $D_{k+1,\infty}$ are  normal crossing divisors, so in fact what we identify is their resolutions,  via a canonical map (\ref{F=D}). Their normal bundles are  canonically dual to each other, see (\ref{nor.op}), and the $\cx^*$ action on $N$ induces an anti-diagonal $\cx^*$ action in the normal bundle 
\best
N_{F_{k}} \oplus   N_{D_{k+1,\infty}}
\eest 
of each component $F_{k}= D_{k+1,\infty}$ of the singular divisor $W_1$ of $X_1$, where  $k\ge 0$. 
 \end{rem} 
\begin{ex}\label{ex.4dim} Assume $X$ has 4 real dimensions, and that the normal crossing divisor 
\best
V= V_1\ma\cup_{p_1=p_2} V_2
\eest
is the union of two submanifolds $V_1$ and $V_2$ intersecting only in  a point $p=p_1=p_2$. 

After rescaling once, we get 3 main pieces $X$, $\F_1$ and $\F_2$ together this an attaching map, see the left hand side of Figure \ref{F.1}.
Here  $\F_1$ is a $\P^1$ bundle over  $\wt {V^1}=V_1\sqcup V_2$, while  $\F_2$ is just $\P^1\ti \P^1$ (over the point $\wt V^2=p$). The divisor $D_{1,\infty}\subset \F_1$ is a copy of $\wt  V=V_1\sqcup V_2$ and it is attached to $F_0=V=V_1\ma\cup_{p_1=p_2} V_2 \subset X$. Similarly, $D_{2,\infty}\subset \F_2$ is $\P^1\ti \infty\cup \infty \ti \P^1$, and it is attached to $F_1$, which is the disjoint union of  two fibers of $\F_1$ over the points $p_1$ and $p_2$ in $V_1\sqcup V_2$. 

Note that   $\F_1$ does not descend as a bundle over $V$: the  two fibers of $\F_1$ over singular locus $p_1$ and $p_2$ are not identified with each other, but rather each gets identified with different of fibers of $\F_2$. 
\end{ex}


\begin{ex}\label{ex.4dim.b} Assume $X$ has 4 real dimensions, but the normal crossing divisor $V$ has only one component, self intersecting itself in just a point $p$. Locally, the situation looks just like the one in Example  \ref{ex.4dim}, with $V$ having two local branches meeting at $p$. The only difference is that globally $\wt V$ has only one connected component containing both points $p_1$ and $p_2$, see right hand side of Figure \ref{F.1}. 
\end{ex}

 \begin{rem} In the discussion above, we had a rescaling parameter $\la$ normal to $V$, which means that we considered the action of $\la\in \cx^*$ on the normal bundle $N$ over $\wt V$. If $\wt V$ has several connected components
 \bear
 \wt V=\ma\sqcup_{c\in C} \wt V_c
 \eear
then we could  independently  rescale normal to each one of them; this gives a $\la \in (\cx^*)^C$ action, rather than just the diagonal one we considered before. Rescaling in all these independent directions now gives a multi-building, where each floor has a level associated to each connected component of $\wt V$. We could talk about a room of the building which is on level one normal to some of the components, but level zero normal to other components. 
\end{rem}
By iterating the rescaling process,  we obtain level $m$ buildings where we rescale $m$ times normal to $V$, or more generally multi-buildings with  $m_c$ levels in the normal direction to each connected component  $\wt V_c$ of $\wt V$. 

 \begin{defn}\label{D.level.m} A {\em level $m$ building} is a singular space $X_m$ with a singular divisor $V_m$, called the {\em  zero divisor},  that is obtained recursively from  $(X_{m-1}, V_{m-1})$ by iterating the level one building procedure.  In particular, a level $m$ building  comes with  (i) a resolution $(\wt X_m, \wt D_m)$, (ii) an attaching map  $\xi_m$ and (iii) a collapsing map $p_m$.    The attaching map 
 \best
 \xi_m: (\wt X_m,\wt D_m )\ra (X_m, D_m)
 \eest  
attaches the floors together producing  the singular locus $W_m$ of $X_m$, where $D_m=W_m\cup V_m$ is the total divisor. The collapsing map  
 \bear\label{D.col.m}
 p_m:(X_m, V_m) \ra (X_{m-1}, V_{m-1})
 \eear
is the identity on $X_{m-1}\subset X_m$ and fiberwise collapses the top floor.
\end{defn} 
This inductively defines a tower of buildings 
\bear\label{D.tower.build} 
X =X_0\subseteq \dots  \subseteq X_{m-1}\subseteq X_m\subseteq \dots  
\eear
together with their resolutions $\wt X_m$ and natural maps between them. The connected components of the resolution $\wt X_m$ are called {\em rooms}, while the $i$'th floor of $X_m$ is  $X_i\setminus X_{i-1}$. The full projection 
 \bear\label{col.0}
 p:(X_m, V_m) \ra (X_0, V_0)=(X, V)
 \eear
collapses all positive floors, leaving the bottom one unaffected. Note that as we add floors, the building grows bigger in several (local) directions. Starting with  $(X_{m-1}, V_{m-1})$, on the top floor we add 
a new piece 
 \bear\label{D.piece}
 \F_{k, m}=\F_k(V_{m-1}).
 \eear
which is a $(\P^1)^k$-bundle  over the depth $k$ stratum of $V_{m-1}$, one for each $k\ge 1$. A point $\wt y$ of $\F_{k, m}$ has an associated level in each of the $k$ local directions,  keeping track of how many times the building was rescaled in that direction to get $\wt y$, with level zero unrescaled. Unlike the floors which partition the building, a point could be in several levels at once (in different local directions), so levels of the building overlap over depth $k\ge 2$ points, see \S \ref{rem.def.multilevel}.  

\begin{rem}\label{R.depth-level} We use {\bf depth}  to measure how many local branches of $V$ meet at a point. We now also have {\bf floors} and {\bf levels}, which measure slightly differently how the target was rescaled.    
\end{rem}
\begin{ex}\label{ex.rescale.D2}{\bf (Fundamental model)} When $(X, V)=(D^2, 0)$ is the disk, the level $m$ building is 
\bear\label{D2.rescaled}
X_m= D^2 \ma\cup_{0=\infty} \P^1\ma\cup_{0=\infty} \dots \ma\cup_{0=\infty} \P^1
\eear 
Its resolution $\wt X_m=  D^2\sqcup \P^1\sqcup \dots \sqcup \P^1$ comes with a locally constant {\em level map} 
\bear\label{D.level.disk}
l: \wt X_m\ra \N
\eear 
indexing its $m+1$ components in increasing consecutive order with $D^2$ on level zero. It uniquely descends to an upper semicontinuous level map $l: X_m\ra \N$, compatible with the inclusions \eqref{D.tower.build} but discontinuous along the total divisor. Note that the level zero of $X_m$ is $D^2\setminus 0$ while $0$ is already on level 1.  

Next, each point of $\wt X_m$ comes with a sign $\ep= \pm 1, 0$ keeping track of whether its coordinate is $\infty$, 0 or neither, inducing the {\em sign map}  
\bear\label{D.sign.disk}
\ep: \wt X_m\ra \{0, \pm 1\}
\eear 
Its restriction to $\wt D_m$ is locally constant, indexing the zero and respectively the infinity divisor; intrinsically $\ep$ keeps track of the weight of the $\cx^*$ action. Together, the restriction of $l$ and $\ep$ to $\wt D_m$ index the connected components of the total divisor $\wt D_m$.  

In coordinates, let $z=z_0$ be a complex coordinate on $D^2$ and $z_l$ be a homogenous coordinate on the $l$'th copy of $\P^1$ in \eqref{D2.rescaled}. These are coordinates on the resolution of \eqref{D2.rescaled}, with $z_l$ a coordinate on level $l$ and $z_l=0,\infty$ describing its zero and respectively infinity divisor.  The attaching map $\xi:\wt X_m\ra X_m$ identifies the point $z_{l-1}=0$ with the point $z_{l}=\infty$ to produce the $l$'th node $y_l$ and define the building  \eqref{D2.rescaled}. Conversely, the $l$'th node $y_l$ has two lifts $y_l^\pm$ to the resolution, on consecutive levels, with $y_l^+$ on the zero divisor in level $l-1$ and $y_l^-$ in the infinity divisor on level $l$.  The $\cx^*$ action $(\la, z)\ra \la z$ describes infinitesimal  coordinate changes on $D^2$ at 0. There is a similar $\cx^*$ action on each of the $m$ copies of $\P^1$, inducing all together a $(\cx^*)^m$ action on  \eqref{D2.rescaled}. 

The level and sign maps \eqref{D.level.disk}-\eqref{D.sign.disk} are an intrinsic combinatorial way to keep track of the topology of the tower of buildings \eqref{D.tower.build} obtained by rescaling the disk around 0. 
Note that for each positive level $l$, we also have a map $p_l:X_m \ra X_{m-1}$ that collapses that level, and more generally $p_J: X_m\ra X_{m-j}$ that collapses an order $j$ subset $J$ of the positive levels of \eqref{D2.rescaled}. 
\end{ex}
\par

\begin{wrapfigure}[12]{r}{0.40\textwidth}
\centering
\includegraphics[width=0.35\textwidth]{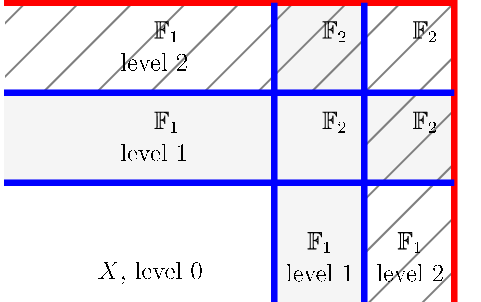}
\caption{A level 2 building}  \label{F.2} 
\end{wrapfigure}

\begin{ex}\label{ex4.rescale.1} In Example \ref{ex.4dim}, when we rescale again, a new level forms with five new rooms, see Figure \ref{F.2}. 
Two of the rooms are  
$\P^1$ bundles over $V_1$ and $V_2$ respectively, but we now have three extra copies of  $\F_2=\P^1\ti \P^1$. The way they come about is as follows: 

After the first rescaling, the zero divisor $V_{(1)}$ of the first floor consists of 4 pieces: $V_1$, $V_2$ but also two $\P^1$'s intersecting in  a point $p^{1}$ (coming from the zero section of depth two piece $\F_{2,1}$ on the first floor). 

When we rescale the second time, $\F_{2,2}$ is still a $\P^1\ti\P^1$ over the point $p^1$, but  $\F_{1,2}$ is now a  $\P^1$ bundle over 
 $(V_1\sqcup \P^1)\sqcup (V_2\sqcup \P^1)$,  so it has four pieces: a $\P^1$ bundle over  $V_1$ and respectively $V_2$, and two other $\P^1\ti \P^1$ pieces coming from rescaling over the two $\P^1$ fibers in $V_{(1)}$. 
\end{ex} 

\begingroup
\endgroup

 \par 
 
 All together this  level 2 building has $4$ copies of $\P^1\ti\P^1$, two of them simultaneously in both level one and two, see Figure \ref{F.2}, where level 1 is shaded and level 2 is hashed. 
 
 In general, such level $m$ building has $m^2$ copies of $\P^1\ti \P^1$, and $m$ copies of the $\P^1$ bundle over $V_1\sqcup V_2$. 
 
\begin{ex}\label{ex.4.rescale.c} In Example \ref{ex.4dim}, $\wt V=V_1\sqcup V_2$ has two connected components, so we can also independently rescale normal to each of them, getting instead a multi-building: a level (2,2) multi-building looks the same as the level $2$ building, except that we keep track separately of how many levels we have in each direction. For example, the  4 pieces $\P^1\ti \P^1$ described above now land one on each level $(i,\;j)$ for  $i,\; j=1$ or 2, while before the piece $(i, j)$ was in {\em both} level $i$ and level $j$ (see Figure \ref{F.2}).  More generally, the level $m$ building from the example above can be regarded as a level $(m,m)$ multibuilding when we independently rescale in the two directions, with exactly one copy of $\P^1\ti \P^1$ on each level $(i,\;j)$ for $i,\;j=1, \dots m$. But in this case a multi-building may have different number of levels in each direction, e.g. just one level normal to $V_1$ and three normal to $V_2$. 
 \end{ex}
 \par
 
\begin{ex} In Example \ref{ex.4dim.b}, locally everything looks the same as in Example  \ref{ex4.rescale.1} 
and even near $p$ we have two independent local directions in which we could rescale as in Example \ref{ex.4.rescale.c}. However, because $\wt V$ is now connected, globally there is only one scaling parameter $\la\in \cx^*$ normal to $V$, so the two local scaling parameters at $p$ are no longer independent (they are essentially equal).  
\end{ex} 

\subsection{Stratifications of buildings and semilocal models} \label{rem.def.multilevel} 
A level $m$ building $X_m$ comes with many stratifications, each recording relevant topological information about $X_m$. (See \S \ref{S.A.stratif} for a review of basic notions associated to stratifications.) First of all, the normal crossing divisor $\wt D_m$ induces a stratification $\textsf{depth}:\wt X_m \ra \N$ keeping track of how many branches of $D_m$ meet at a point (cf. Remark \ref{N.strata}). It pushes forward $\xi_*\textsf{depth}:X_m \ra \N$ to one on $X_m$ via the attaching map $\xi:\wt X_m \ra X_m$.  But $(X, V)$ is also stratified by depth, so we get a pullback stratification $p^*\textsf{depth}:X_m \ra \N$ induced by the collapsing map $p:X_m\ra X$, with  $p^*\textsf{depth}\ge \xi_* \textsf{depth}$. Therefore the stratification $\xi_* \textsf{depth}$ is finer, recording the fact that the depth of a point in $X_m$ is at most that of its projection in $X$. 
 
These enter in the description of the semilocal model of the maps $\wt X_m \ra X_m$ and  $X_m\ra X$ around one of their fibers. For each point $y\in X_m$, denote by $\wt y\in \wt X_m$ one of its resolutions and by $p\in X$ its projection. Denote by $k$ the depth of $p$ in $X$ and assume $k\ge 1$. Locally index the $k$ branches of $V$ near $p$ by $I$, and denote by $F_p$ the fiber of $X_m\ra X$ at $p$.


A neighborhood $U_p$ of $F_p$ in $X_m$ is a product of a small neighborhood $O_p$ of $p$ in $V^k$ and $k$ copies of the $m$-times rescaled disk $(D^2)_m$ of (\ref{D2.rescaled}),  one factor for each one of the $k$ branches of $V$ at $p$. This describes not only the tower of $(m+1)^k$ pieces of the resolution $\wt X_m$  with its total divisor $\wt D_m$, but also their attaching map, just as in Example \ref{ex.rescale.D2}, except that now we have $k$ directions to keep track of instead of one. For each $i\in I$ fixed, each point  $\wt y$ in the fiber of $\wt X_m \ra X_m \ra X$  over $p$ now comes decorated by both a level $l_i(\wt y)$ defined by \eqref{D.level.disk}, inducing the multilevel map
\bear\label{def.level.map}
&& l_{\wt y}: I\ra \N
\eear
and a sign $\ep_i(\wt y)$ defined by \eqref{D.sign.disk}, giving the multisign map 
\bear
\label{def.sign.map} 
&&\ep_{\wt y}: I \ra \{0, \pm 1\}. 
\eear
In general we can have a point which is (a) on different levels in different directions, which is the information recorded by $l$  and (b) on the zero divisor in some of the directions, on infinity divisor in other directions, and then in some other directions on neither, and this is precisely what $\ep$ records. The order of the domain of $l$ is the depth $k$ of $p$, while that of the support of $\ep$ is depth $\wt k$ of $\wt y$ (where $k \ge \wt k$). 

As long as the depth $k$ of the projection $p$ is constant, the multilevel map \eqref{def.level.map} is {\em locally constant} in $\wt y$, indexing the connected components of $\wt X_m$, and descends to an upper semicontinuous map in $y\in X_m$ (discontinuous along $D_m$).  Similarly, as long as the  depth $\wt k$ of $\wt y$ is constant, the restriction of the multisign map \eqref{def.sign.map}  to $\wt D_m$ is locally constant in $\wt y$, so together with the multilevel map $l$ it indexes the  connected components of each open stratum of $\wt D_m$  
(see Example \ref{A.ex.V.tower} for an intrinsic view point). 

If  $I$ is a set of order $k$, let $V_I\ra V^k$  denote the resolution \eqref{A.D.V.I} obtained by indexing  $k$ branches of $V$ by $I$. The discussion above uniquely describes a stratum $D_{J}$ of the resolution of $(\wt X_m, \wt D_m)$ in terms of 
\begin{enumerate}[(i)]
\item  the stratum $V_I$ of the resolution of $(X, V)$ that $D_J$ projects to under  $X_m\ra X$, and  
\item a multilevel $l:I \ra \N$ and multisign $\ep:I \ra \{ 0, \pm 1\}$ map.
\end{enumerate}
The strata $D_J$ of $\wt X_m$ are therefore indexed by pairs 
\bear\label{J=I.level}
J=(I,  \ep\ti l)
\eear 
with a symmetric group action $S_I$ reordering $I$. The strata $V_J$ of the zero divisor $\wt V_m$ corresponds to data \eqref{J=I.level}  for which 
\bear\label{D.ind.l.0} 
\text{ there exists at least one direction $i\in I$ with $l(i)=m$ and $ \ep(i)= +1$}
\eear while the other strata of the singular divisor come in dual pairs $W_{J_\pm}$ indexed by $J_\pm=(I, \ep_\pm \ti l_\pm)$ where 
\bear\label{l-l=ep}
\ep_+(i)=-\ep_-(i)= l_-(i)- l_+(i)\quad  \mbox {for all }i\in I 
\eear
and $\ep_\pm \not\equiv 0$, 
keeping track of the fact that a level $m$ building is obtained in each direction $i$ not only by joining together the zero and the infinity divisor on consecutive levels (when $\ep_\pm(i)=\pm 1$) but also joining together fibers of the pieces in the same level (when $\ep_\pm(i)=0$), see  (\ref{F=D}).  
 \begin{defn}\label{D.level.j} The {\em level $j$} of $X_m$ is  $L_j= \{ y\in X_m\; | \; l_y^{-1}(j) \ne \emptyset\}$, i.e. the collection of points which are on level $j$ in at least one direction. 
  \end{defn} 
A point in $X\setminus V$ is by convention in level zero. Note that $L_j$ is well defined, independent of the choice of indexing of the local branches of $V$ by $I$ (i.e. preserved by the symmetric group action $S_I$ reordering the branches). Intrinsically, the multi-level map $l_y$ describes a partition $I_0, \dots, I_m$ of the $k$ local branches $I$ of $V$ at $p(y)$  according to its value on each branch, where $I_j=l^{-1}_y(j)$ has order $\mu_j\ge 0$, see  Example \ref{E.S.nc}. Level $L_j$ consists of points with $\mu_j\ge 1$. 

\begin{rem}\label{cx.m.act} The $\cx^*$  action of  $N_V$ induces a $(\cx^*)^m$ action on a level $m$ building such that each factor $\al_l\in \cx^*$ rescales the level $l\ge 1$ piece of $X_m$  normal to its zero section, fixes $X\setminus V$ pointwise, and preserves the total divisor and its stratification, but not pointwise. It is modeled by the $(\cx^*)^m$ action on the rescaled disk  $(D^2)_m$ of Example \ref{ex.rescale.D2} in which the $l$'th factor $\al_l\in\cx^*$ acts nontrivially only on the $\P^1$ components of $(D^2)_m$ in level $l\ge 1$. The fiber $F_p$ of $X_m\ra X$ is a product of $k$ copies of $(D^2)_m$,  where the $i$'th $D^2$ is intrinsically the fiber $N_{i, p}$ of the normal bundle to the $i$'th branch of $V$ at $p$, thus $F_p$ comes with an induced diagonal action. 

The $l$'th factor  $\al_l\in\cx^*$ acts nontrivially only on the rooms of the building that are on level $l\ge 1$. 
 \end{rem}

\begin{rem} For each positive level $l$ there is a map collapsing that level 
\bear\label{D.p.L}
 p_l:(X_m, V_m) \ra (X_{m-1}, V_{m-1})
 \eear
defined on any building $X_m$ with at least $l$ floors, and more generally a collapsing map $p_J$  for any subset $J$ of the levels $\{1, \dots, m\}$.  Over a depth $k$ point $p\in V$ the map (\ref{D.p.L}) is the product of the  collapsing maps $p_l$ on each one of the $k$ copies of the fundamental model $(D^2)_m$.  In  Figure 2, $p_1$ collapses the first level (the shaded part) and $p_2$ the second level (the hatched part).  
\end{rem} 
\begin{rem}\label{R.ext.J.X.m}
The notion of a tamed pair $(\om, J)$ on a level $m$ building $(X_m, V_m)$ compatible with the total divisor $D_m$ is defined via one on the resolution $(\wt X_m, \wt D_m)$ that matches along the singular locus. Among these pairs, there is a subset which is also compatible with one of the collapsing maps (\ref{D.p.L}) or more generally with all of them, see \S \ref{A.spaces.param}.   

Moreover, for each pair $(\om, J)\in \J(X, V)$ on $X$, the rescaling procedure induces a $J$ on $X$, compatible with both the total divisor $D_m$ and with each of the collapsing maps.  Its restriction to the level zero is the original $J$ while the restriction to each other piece $\F_k$ is $\cx^*$ invariant. For each $\al>0$ sufficiently small, we can find an $\om^\al$ on $X_m$ such that $(\om^\al, J)$ is a tamed pair, compatible with the total divisor $D_m$, such that its restriction to level zero is the original structure $(\om, J)$ on $X$ while its restriction to each other piece $\F_k$ is compatible with the projection $\pi:\F_k \ra \wt V^k$ and the fiberwise $S^1$-action is Hamiltonian. In fact, we can arrange that the restriction of $\om^\al$ to $\F_k$ has the form $\pi^*\om_{V^k} + \al^2 \tau$, where $\al$ is the size of each $\P^1$ fiber, see Remark \ref{R.ep.on.F}. 

Observe that if we let  $(X', V')=(\ov\U\ti X,\; \ov\U\ti V)$ then the level $m$ building associated to $(X', V')$ is the product of $\ov\U$ with the level $m$ building associated to $(X, V)$. So the discussion above extends to  Gromov-type perturbations: any $V$-compatible parameter $(J, \nu)\in \JV(X, V)$ induces one on $X_m$ compatible with both $D_m$ and the collapsing maps.
\end{rem}

\subsection{Local model near the divisor}\label{S.local.mod}
Rescaling $X$ by a factor of $\la\ne 0$ along a normal crossing divisor $V$ similarly gives rise to manifold we denote 
$X_\la$ and an identification
\bear\label{R.la.2}
R_\la: X\ra X_\la
\eear 
extending (\ref{R.la}) and described in more details below. The manifold $X_\la$ again has two regions, one is the complement of the $|\sqrt \la|$-tubular neighborhood $U_{\la}$ of $V$ in $X$,  on which $R_\la$ is the identity, and the other one is identified with $O_\la$, the complement of the $|\sqrt \la|$-tubular neighborhood of the singular divisor  $W_1$ in $\F_V$. The only difference now is that we have several overlapping local models, coming from the stratification of $V$. 

This perspective allows us to think of the rescaled $X$ as a sequence of  manifolds $X_\la$ with varying $J_\la=R_\la^* J$, which as $\la\ra 0$ converge to a level one building $X_1$ with an induced $J$ on its resolution $\wt X_1$. In fact, $X_\la$ can be thought as an iterated symplectic sum: fix a level one building $X_1$ as in Definition \ref{D.build.1}, with appropriate identifications along corresponding divisors, including fixed isomorphisms (\ref{nor.op}). For each (small) gluing parameter $\la\in \cx^* $  we get the 'symplectic sum' $X_\la$ (diffeomorphic to $X$, but with a deformed symplectic form) and $X_\la$ converges to $X_1$ as $\la \ra 0$. This describes a local family of deformations 
\bear\label{cal.X}
\xymatrix{ 
X_\la\ar@{^(->}[r]&{\mathcal X}\ar[r]^\rho & B}
\eear
with fiber $X_\la$ over $\la\in B\setminus 0$ and central fiber $X_1$.  Moreover, (\ref{R.la.2}) can be regarded as  a local trivialization of (\ref{cal.X}) away from $\la=0$. This approach is mentioned in Remark 7.7 of \cite{ip1} and expanded on in \cite{ip2}. 
\begin{rem} When $(X, V)$ is a marked Riemann surface $(\Si, x)$ the level 1 building $\Si_1$ is $\Si$ together one spherical bubble (with two marked points, the zero and the infinity divisor) attached at each of the marked points $x$. This matches the picture appearing in Gromov compactness on the domain side. The family of deformations (\ref{cal.X})  is obtained by fixing an isomorphism between the two tangent planes at the node, for example by adding an extra marked point to make the bubble stable,  and then (uniquely) identifying the new family of curves with the fibers of the universal curve. The universal curve is a locally trivial fibration away from the nodal locus and has specific local models around the nodal locus.  
\end{rem}
\par

We next describe the rescaling procedure in more detail, working in regions which are obtained from neighborhoods of depth $k$ stratum of $V$ after removing neighborhoods of the higher depth stratum. 

\begin{wrapfigure}[13]{l}{0.28\textwidth}
\vskip-.1in
\centering
\includegraphics[width=0.22\textwidth]{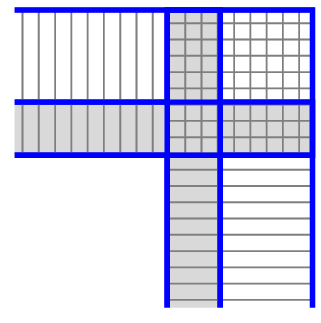}
\caption{Regions of $X$}\label{F.3}
\end{wrapfigure}
Recall that a normal crossing divisor $V\subset X$ comes with a smooth resolution $\wt V\ra V$, a normal bundle $\pi:N\ra \wt V$, and an immersion $\iota: \wt U \ra X$ defined on some disk bundle $\wt U$ of $N$. The restriction $\iota: \wt V\ra V$ to the zero section defines a stratification $\mathsf{depth}$ of $X$ (by the number of elements in the fiber of $\iota$) with closed strata $V^k$.

Denote by $\iota_\de$ the restriction of $\iota$ to the  $\de$-disk bundle $\wt U_\de$ of $N$, where $\de$ is sufficiently small.  Then $\iota_\de$ defines a stratification 
\bear\label{D.str.X.de}
\hskip1in \mathsf{depth}_\de:X\ra \N
\eear
of $X$ with closed strata $U_\de^k$ where $U_\de^1=U_\de$ is a neighborhood of $V$ and $U_\de^k$ consists of points $x\in X$ for which $\iota_\de^{-1}(x)$ has at least $k$ points. Note that $U^0_\de=X$ is independent of $\de$.

In Figure  \ref{F.3}, $U_\de$ is the hatched region while $U^2_\de$ is the doubly hatched one. The maps $\mathsf{depth}_\de$ are decreasing in $\de$ and limit as $\de\ra 0$ to $\mathsf{depth}$; the regions $\{ U_\de^k\}_\de$ provide a basis of neighborhoods of $V^k$ in $X$. For each $0<\ep<\de$, consider the annular  `neck region' 
\bear\label{defn.an.r}
A(\ep, \de) =\iota (\wt U_\de\setminus \wt U_\ep)
\eear
where the two stratifications differ from each other; it  corresponds to the shaded region in Figure \ref{F.3}. Note that $A(\ep, \de)$ intersects $V$ near depth $k\ge 2$ points, so it is  different from $U_\de\setminus U_\ep$. 
\par 

Let $S_\de$ be the radius $\de>0$ circle bundle of $N$. The complement of its image $A_\de=\iota(S_\de)$ separates $X$ into open pieces $X_\de^k= \mathrm{int} \;U_\de^k \setminus U_\de^{k+1}$ which lie in a neighborhood of $V^k$ but stay away from a neighborhood of the higher depth stratum $V^{k+1}$:
\bear\label{X.disjoint}
X\setminus A_\de= \ma \sqcup_{k\ge 0}  X^k_\de
\eear 
Lemma \ref{A.L.V-I.de} shows that  each piece $X_\de^k$ for $k\ge 1$ comes with 
\bear\label{proj.emb.X.k}
\text{ a projection $\pi_k:X_\de^k \ra V^k\setminus U_\de^{k+1}$ and an embedding $\eta_k:X^k_\de \hookrightarrow N_{V_k}\subset \F_k$}
\eear 
whose image is a product of $\de$-disk bundles in the fibers $(\P^1)^k$ of $\F_k$, but  with the $\de$-neighborhood of the fiber divisor $F_k$ removed. In other words,  $X^k_\de$ can be regarded, via $\eta_k$, as the subset of $\F_k$ whose complement is the union of the $1/\de$ neighborhood of the infinity divisor and the $\de$-neighborhood of the fiber divisor. But each piece $\F_k$ has a map $R_\la^k$ that fiberwise rescales $\F_k$ by the factor  $\la\in \cx^*$, see Remark \ref{cx.act}.  $R_\la^k$ identifies $X_\de^k$ in $\F_k$ with the complement $O_{\de}^k(\la)$ of the union of the $|\la|/\de$ neighborhood of the infinity divisor $D_{k, \infty}$ and the $\de$-neighborhood of the fiber divisor $F_k$. For $\de=|\sqrt \la|\ra 0$  the image of $X_\de^k$ under $R_\la^k$ converges to $\F_k \setminus (D_{k, \infty}\cup F_k)$. The rescaling map \eqref{R.la.2} is defined by first decomposing $X$ into pieces $X_{\de}^k$ where $\de=|\sqrt \la|$ and then identifying each piece $X_{\de}^k\cong O_{\de}^k(\la)$ using $R^k_\la$, with $R_\la$ being the identity on $X^0_\de=X\setminus U_{\de}$. 

In local coordinates, assume $p$ is a depth $k$ point of $V$ (but away from the higher depth strata);  locally index  the $k$ branches $I$ of $V$ coming together at $p$, and choose local coordinates $u_1, \dots, u_k$ normal to each one of these branches (so the $j$'th branch is given by $u_j=0$) such that the $\de$-neighborhood $U_\de$ of $V$ is given in these coordinates by $|u_i|\le \de$ for {\em some} $i\in I$. Then the region $X_\de^k$ is described by 
\bear\label{def.u.k}
|u_i|\le \de \mbox{ for all }i=1, \dots, k
\eear
with $p\in V^k\setminus U^{k+1}_\de$.  Under the change of coordinates $v_i= \la/u_i$ for each $i\in I$,  the region $X_\de^k$ is identified with the region $O_{\de}^k(\la)$ in $\F_k$ described by the equations 
\bear\label{def.v.k}
|v_i|\ge |\la| /\de \mbox{ for all }i=1, \dots, k
\eear
where $p\in V^k\setminus U^{k+1}_\de$. In particular, this change of coordinates provides an identification $R_\la$ of $X$, broken into pieces $X^k_{\de}$ with $\de=|\sqrt\la|$ as above, with  the complement of the $|\sqrt\la|$-tubular neighborhood of the singular divisor $W_1$ in a level one building $X_1$ (broken into pieces in $\F_k$), exactly as suggested by the iterated symplectic sum construction \eqref{cal.X}. 

In fact, the semi-local model of $X_\la$ over a point $p\in V^k\setminus U_\de^{k+1}$ is given by the locus of the equations 
\bear\label{eq.xy=l}
u_i v_i= \la\quad  \mbox{ where } \quad |u_i|\le \de \quad  \mbox{ for all }i\in I
\eear
where $u_i$ is a coordinate in $N_i$, the normal bundle to the $i$'th local branch of $V$ and $v_i$ is the dual coordinate in the dual bundle $N_i^*$ (which is allowed to equal infinity). Note that these equations are invariant under reordering of the branches, so they describe an intrinsic subset of $N_{V^k} \ti  \F_k$ where our semi-local analysis takes place. In particular, $(p, u, v)$ can be regarded as local coordinates on a thickening of the total space of the deformation \eqref{cal.X} in which the fibers $X_\la$ are cut out by \eqref{eq.xy=l}. 

The {\em $\de$-neck}  $N_\la(\de)$ of $X_\la$ is the region in the above coordinates where 
\bear\label{eq.de-neck}
\mbox{ $|u_i| < \de$ and $|v_i|<\de$ for some } i 
\eear
and it globally corresponds to the annular region  $A(|\la|/\de, \de)$ around $V$ in $X$. It decomposes into pieces depending how many of the $i\in I$ satisfy \eqref{eq.de-neck}, keeping track of the overlap of these semi-local models. 

The {\em upper hemisphere} $H_\la(\de)$ of $X_\la$ corresponds to the region 
$A(|\la|,\de)$ in $X$; in coordinates it is 
\bear
\mbox{ $|u_i| < \de$ and  $|v_i|< 1$ for some } i. 
\eear
Its intersection with $X_\de^k$ maps under $R_\la^k$ to the part of the upper hemisphere of $\F_k$ in $O_\de^k(\la)$.   

As $\la\ra 0$, $X_\la$ converges to a level one building $X_1$. For $\la$ sufficiently small ($|\la|<\de^2$), the $\de$-neck of $X_\la$ is described by the local models \eqref{eq.xy=l}-\eqref{eq.de-neck}, while its complement is canonically identified with the complement of a neighborhood of $W_1$ in $X_1$. As $\la\ra 0$ and then $\de\ra 0$, this neighborhood expands to $X_1\setminus W_1$, while the upper hemisphere region $H_\la(\de)$ converges to the upper hemisphere of $\F_V$.

So over a point $p\in V^k\setminus U^{k+1}_\de$, the level one building $X_1$ is described by the equations:
\bear\label{2-k-piece}
u_1 v_1=0, \quad \dots\quad,  u_k v_k=0
\eear 
where $|u_i|\le \de$ for all $i$, regarded as an intrinsic subset of $N_{V^k}\ti \F_k$ (or more precisely inside the pullback 
of $\F_k$ over $N_{V^k}$).

This describes the $2^k$ local pieces of $X_1$ coming together at $p$ along the singular locus $W_1$: each piece of $X_1$ is described by the vanishing of exactly $k$ coordinates, but some of them may be $u_i\in N_i$ in which case the rest are the complementarily indexed ones $v_j\in N_j^*$. The divisor $W_1$ has several local branches, each one described by the further vanishing of one of the remaining coordinates $u_i$ or $v_i$, matching the description in  \S \ref{rem.def.multilevel}. Note that the $u$'s are coordinates on level zero of $\wt X_1$, while the $v$'s are coordinates on level one (in that direction), with $u=0$ corresponding to the zero divisor and $v=0$ the infinity one (in that direction). The multilevel \eqref{def.level.map} and multisign map \eqref{def.sign.map} encode this information. 
\begin{rem}\label{rem.coord-iterate.m}
The rescaling procedure can be iterated finitely many times: start with $X$, rescale it  by $\la_1$, then rescale again the resulting manifold by $\la_2$,  etc. The limit as all $\la_a\ra 0$ is a level $m$ building $X_m$, with similar semi-local models described over a depth $k$ point $p$ as intrinsic subsets of $(\F_k\times \F_k)^m$, modeled in each direction $i$ by the process of $m$ times rescaling a disk at the origin, see \S \ref{rem.def.multilevel} and Figure 2. So one way 
to describe the iterated rescaling is to start by choosing coordinates $u_{i, 0}=u_i$ on $X$ normal to $V$ at $p$, together with $m$ other sets of dual normal coordinates: 
\bear\label{dual.nor.coord}
\mbox{$u_{i, l}\in N_i\cup \infty $  and $v_{i,l}\in N_i^*\cup \infty$  with $v_{i, l}=u_{i, l}^{-1}$ }
\eear
for all $ l=1, \dots, m$ and all $i\in I$. These provide semi-local coordinates on $X_m$ in a neighborhood of the fiber $F_p$ over $p$ of the collapsing map $X_m\ra X$. Rescaling $X$ means we 
\bear\label{resc.change.coord}
\mbox{make a change of coordinates $u_{i, l-1} v_{i, l} = \la_l$ at step $l$,} 
\eear 
for each $l=0, \dots, m$  and all $i\in I$. 

Intrinsically, this describes (i) a complete family ${\mathcal X}_m \ra B$ of global deformations (\ref{cal.X}) of $X_m$ indexed by the parameters $\la=(\la_1, \dots \la_m)\in B$, (ii) local trivializations of $\x$  away from the necks and (iii) local models of $\x$ in the necks. When restricted to subsets $B'\subset B$ we get partial families of deformations; in particular, if $B'$ is the subset where $i$ of the components of $\la$ are equal to zero, the fibers of ${\mathcal X}\ra B'$ are level $m-i$ buildings obtained by smoothing off only the singular divisors corresponding to nonzero gluing parameters. 
\end{rem} 
\begin{rem}\label{R.ext.J.cal.X}
For such deformation $\rho:{\mathcal X} \ra B$ of a level $m$ building $(X_m, V_m)$ we can use the semilocal models above to define a notion of a tamed pair $(\om, J)$ on the total space $\mathcal X$ compatible with the (singular) fibration $\rho$: it is an element $(\om, J)\in \J(\x)$ whose restriction to each fiber $X_\la$ is a tamed pair 
$(\om_\la, J_\la)$ compatible with the total divisor of $X_\la$ (cf. \S \ref{A.spaces.param}). Any $D_m$-compatible, tamed pair $(\om, J)$ on $X_m$ extends to a tamed pair on 
${\mathcal X} \ra B$ compatible with the fibration by the symplectic sum construction. 

This extends naturally  to Gromov-type perturbations using the fact that $ {\mathcal X}'=\ov\U\ti {\mathcal X} \ra B$ is a family of deformations of $X'_m= \ov\U\ti X_m$. In particular, any parameter $(J, \nu)$ on $X_m$ which is compatible with  the total divisor (for example one obtained from the rescaling process) extends to a parameter on  ${\mathcal X} \ra B$ compatible with the fibration. 
\end{rem}

\setcounter{equation}{0}
\section{Refined limits of a sequence of maps}\label{s4}
\medskip

Consider now a sequence $f_n$  of maps in $\M(X,V)$ which limits in the usual stable map compactification to a map $f_0$ which has components in $V$. In this section, after rescaling  $X$ once around $V$ we construct refined limits $f$ into a level one building $X_1$ such that $f$ no nontrivial components in the singular locus, projects to $f_0$ under the collapsing map $X_1\ra X$, but the depth of at least one nontrivial component has strictly decreased in the refined limit. This is the key step in the inductive rescaling procedure used in the next section to construct the relatively stable map compactification. 

We also need a detailed description of the collection of such possible refined limits $f$. In the case $V$ is a smooth divisor, we  proved in  \cite{ip2} that the limit satisfies a matching condition along the singular divisor. The arguments used there are semi-local (in a tubular neighborhood of $V$), and extend to the case when $V$ is a normal crossing divisor by similarly working in neighborhoods of each depth $k$ stratum of $V$, where we have $k$ different normal directions to $V$. When the limit $f$ has no trivial components in $W_1$, we also get a matching condition along $W_1$.  First of all, we shall see that similarly  $f^{-1}(W_1)$ consists of nodes of the domain.  The naive condition is that at each such node, $f$ has matching intersection points with  $W_1$, including the multiplicities, as was the case in \cite{ip1}. But it turns out that at points of depth $k\ge 2$,  this is not enough to define relative GW invariants, and needs to be further refined. 

\subsection{Maps into resolutions and deformations of buildings} \label{S.maps.resol}

In \S\ref{s3} we described the result of rescaling the target normal to $V$, which can be regarded as constructing a family ${\mathcal X}\ra B$ of deformations $X_\la$ of $X$ converging as $\la \ra 0$ to a singular space, a level $m$ building $X_m$; in turn, the level $m$ building was described in terms of its smooth resolution (via attaching and collapsing maps).  In order to state the result of rescaling on a sequence of maps, we similarly need to first define the notions of  resolutions/refinements/projections and deformations/limits in the context of maps. 

First of all, according to our conventions which match those of \cite{ip2}, every map $f:C \ra X$ is regarded as a map $(\phi, f):C \ra \ov \U\ti X$, where the first factor is the universal curve of the domains (with a fixed embedding into $\cx\P^K$ for large $K$).  The {\em energy} of $f$ in a region $N$ of the target is defined by 
 \bear\label{energy}
E(f; N)= \frac 12 \int _{f^{-1}(N)} |df|^2 +|d\phi|^2 
\eear
By \eqref{eq.f=nu} the total energy of a $(J, \nu)$-holomorphic map  $f$ is topological, equal to the symplectic area of $(\phi, f)$.  A $(J, \nu)$-holomorphic map  $f:C\ra X$ is {\em stable} iff $\Aut (f)$ is finite $\iff$ each unstable domain component represents a nontrivial homology class in $X$ $\iff$ each domain component has positive energy \eqref{energy}. Denote by $\al_V$ the minimum quanta of energy carried by a stable $(J, \nu)$-holomorphic map into $V$. Note that according to our conventions $\al_V>0$. 

The key idea of Chapter 3 of \cite{ip2} is that when $V$ is smooth and $\de>0$ small, the limits of maps which have energy at most $\al_V/2$ in the $\de$-neck cannot have any components in the singular locus $V=V_\infty$ 
and thus have well defined leading coefficients along $V=V_\infty$ which are furthermore uniformly bounded away from 0 and  infinity (the particular bound depends on $\de$ and the choice of metrics). 

When $V$ has normal crossings, it is no longer true that limits of maps with small energy in the neck have no components in the singular locus, rather they may have what we called  trivial components in Definition 11.1 of \cite{ip2} (when $V$ was smooth). The appropriate extension is:
\begin{defn}\label{def.triv.comp} A {\em trivial component} in $\F_k$ is a non-constant holomorphic map $f: C \ra (\F_k, D_0\cup D_\infty)$ 
decorated with its contact information to $D_0\cup D_\infty$  and such that (a)  its domain $C$ is an unstable sphere with only two contact points $x_\pm$ and  
 (b) its image lies in a fiber $(\P^1)^k$ of $\F_k$. In homogenous coordinates   
 \best
 f(z)= (a_1 z^{s_1}, \dots,  a_k z^{s_k} ) 
 \eest
 where some, but not all  of the $a_i$ could be zero or infinity. 
 \end{defn}
This is a stable map into $\F_k$ which does not have a stable model when projected into $V^k$: its domain is an unstable sphere with just two marked points and so its projection has energy zero. The trivial components are the only stable maps fixed by the $\cx^*$ action on the target. In general they are multiply covered and the Gromov-type perturbation vanishes on them, so they may cause trouble with transversality. Note that as a result of the rescaling process, we have uniform control only on the energy of the projection, as the area of the fiber of $\F_k$ might be arbitrarily small.

\begin{defn} Assume $g:C\ra X_m$ is any continuous map into a level $m$  building. A {\em resolution of $g$} is a continuous map  $\wt g:\wt C \ra \wt X_m$ into the resolution $\wt X_m$ of $X_m$ such that \begin{enumerate}[(i)]
\item  $\wt C$ is a resolution of $C$ obtained by resolving a collection $R$ of its nodes; so  each node $x\in R$ corresponds to two extra marked points $x_\pm$ of $\wt C$;
\item $\wt g$ has no components in the infinity divisor; 
\item $g$ is obtained from $\wt g$ via the attaching map $\xi$ of the domain that identifies $x_-=x_+$ for each $x\in R$ to produce the nodes $R$ of $C$ and also attaches the target $\wt X_m$ together  to obtain the building $X_m$:  
\vskip-.2in
\bear\label{g.lift.g} 
\xymatrix{  
\wt C \ar[r]^{\wt g}\ar[d]_{\xi}&{\wt X_m} \ar[d]^{\xi}
\\ 
C  \ar[r]^{g}&X_m
}
\eear\end{enumerate}
\end{defn} 
A map $g:C\ra X_m$ is called $J$-holomorphic if it has a resolution $\wt g$ which is a $J$ holomorphic map into $\wt X_m$, and $(J, \nu)$-holomorphic if it has a resolution whose graph is a $J_\nu$ holomorphic map into $\ov\U\ti \wt X_m$. It is easy to show that any  $(J, \nu)$-holomorphic map $g:C\ra X_m$ has an essentially unique resolution $\wt g:\wt C \ra \wt X_m$ (up to reordering the extra marked points).\smallskip

\begin{defn} Assume $g: C\ra X_m$ is a stable $(J,\nu)$-holomorphic map into a level $m\ge 0$ building. 
A {\em lift} or {\em refinement} of $g$ is a stable $(J, \nu)$-holomorphic map $g':C'\ra X_{m+1}$ such that $g$ is obtained from $p\circ g':C'\ra X_m$ after contracting trivial domain components:
\bear\label{g.ref.f0} 
\xymatrix{ C'\ar[r]^{g'}\ar[d]_{\ct} & {X_{m+1}}\ar[d]^{p} 
\\ 
C \ar[r]^{g}&X_m
}
\eear
where $p:X_{m+1}\ra X_m$ collapses one of the positive levels. The refinement $g'$ is called {\em strict} if there exists at least one nontrivial component $\Si$ such that 
its depth in $g'$ is strictly less than its depth in $g$. 
The map $g$ is called the {\em pushforward} of $g'$ and denoted $p_{*}(g')$.
\end{defn} 
The pushforward is obtained from $g'$ after collapsing the target $X_{m+1}\ra X_m$ as well as any trivial components landing in the collapsed fibers. For any other component $\Si$ of $C$, $g'|_\Si:\Si\ra X_{m+1}$ projects to $g|_\Si:\Si\ra X_{m}$ and its depth cannot increase under this projection, and  remains constant for all the components on which $p$ is the identity. So if $g$ has no nontrivial components in the total divisor $D_m$, then it does not  have any strict lifts. Even if it did have such a component, the existence of a strict $(J,\nu)$ holomorphic lift of it is not automatic.
\begin{rem}\label{R.lift.xi}  If  $g: \Si\ra V^k$ has depth $k$, then a lift $g'$ to $\F_k$ can be regarded as a section $\xi$ of $g^*N_{V^k}$ with poles along the infinity divisor. 
When $(J,\nu)\in\JV(X, V)$ the map $g'$ is $(J, \nu)$-holomorphic iff $\xi$ is in the kernel of the linearized normal operator $D^N_g$ (cf. (7.1) in \cite{ip1}) associated to $g$. A strict lift simply means that $\xi\not\equiv 0, 
\infty$ and so  $\ker D^N_g\ne 0$, a nontrivial condition whenever the index of $D^N_g$ is negative (cf Lemma 6.4 of \cite{ip1} and Lemma \ref{L.d=d-2}). 
\end{rem}

To define the notion of convergence of maps in this context, assume $g:C\ra X_m$ is a $J$-holomorphic map into a level $m$ building, where the parameter $(\om, J)\in \J(X_m, V_m)$. Denote by $\ov\U\ra \ov\M$ the universal curve classifying $C$, and by $\ov\U_k\ra \ov\M_k$ the universal curve classifying $C$ together with $k$ extra marked points.  Consider a family ${\mathcal X} \ra B$ of deformations \eqref{O.defm.x} of the level $m$ building with a fixed 
extension $(\om, J)\in \J({\mathcal X} \ra B)$ compatible with the fibration. The spaces  $\ov\U_k$, $\ov\M_k$, ${\mathcal X}$, $B$ all come with fixed background metrics. 
\vskip-.1in

\begin{defn}\label{D.Gr.top} 
With the notations above, consider $g_n:C_n\ra X_{\la_n}$ 
a sequence of $J_n$-holomorphic maps into the fibers of ${\mathcal X} \ra B$. We say $(g_n, J_n)$ {\em converges}  to $(g, J)$ denoted $g_n\raG g$ provided:
\begin{enumerate}[(i)]
\item there exist lifts denoted still $C_n$ and $C$ obtained by adding $k$ points such that 
$\Aut \;C=1$ and $C_n\ra C$ in $\ov\M_k$. 
\item the graphs of $g_n$ inside $\ov\U_{k}\ti \mathcal X$ converge in Hausdorff distance and energy density to the graph of $g$.
\item the distance between $J_n$ and $J|_{X_{\la_n}}$ wrt induced metric on $X_{\la_n}$ converges to 0.   
\item $g_n\ra g$ uniformly on compacts away from the nodes of $C$.
\end{enumerate}
\end{defn}
\non For $m=0$ this is equivalent with the usual definition of Gromov convergence in the context of the absolute moduli space. The definition is also consistent with both the  forgetful maps $\ov\U_{k+1} \ra \ov\U_k$ on the domain and the collapsing maps $X_{m+1}\ra X_m$ on the target. In fact, condition (ii) and elliptic regularity implies condition (iv); as we prove below it also implies a controlled convergence in local models around the nodes. 

Note that this notion of convergence is insensitive to inserting constant unstable components. To prevent that from happening one restricts to stable maps, which we are already implicitly doing when we assume that domains have only trivial automorphisms.

\begin{rem} The definition of convergence above has a natural extension to the space of Gromov-type perturbations $\JV$. In the most general case, one starts with any parameter $(J, \nu)\in \JV(\x\ra B)$,  and a sequence $(J_n, \nu_n)\in \JV(X_{\la_n})$. To define  $(g_n, J_n, \nu)\raG (g, J, \nu)$ simply extend condition $(iii)$ to $(iii)'$ the distance between the restrictions of $(J_n, \nu_n)$  and $(J, \nu)$ to $C_n\ti X_{\la_n}$ converge to zero wrt induced metric from $\ov\U_{k}\ti \mathcal X$. 
\end{rem}


\subsection{The limits of a sequence of maps after one rescaling} With these preliminary definitions, we are now ready to state the first rescaling result. For simplicity of notation, everywhere in this section we write $J$ instead of $(J, \nu)$ for any parameter $\JV(X, V)$. 
\begin{prop}\label{L.lim.level.1} Assume $J_n\ra J$ is a  
sequence of parameters in $\JV(X, V)$ and consider a sequence $f_n:C_n\ra X$  of stable $J_n$-holomorphic maps in $\M_s(X,V)$ converging to a stable $J$-holomorphic
limit $f_0:C_0\ra X$ which has at least one (nontrivial) component in the zero divisor $V$.

Then there exists a sequence $\la_n\ra 0$  of rescaling parameters such that, after passing to a subsequence, $R_{\la_{n}} f_{n}:C_{n} \ra X_{\la_{n}}$ have a limit $g:C \ra X_1$ with the following properties:
\begin{enumerate}[(a)]
\item $g$ is a stable $J$-holomorphic map into a level one building $X_1$ which is a strict refinement of  $f_0$ and without any nontrivial components in the singular divisor $W_1$.  
\item $g$ has a resolution $\wt g:\wt C\ra \wt X_1$ with no nontrivial components in the level zero part of the total divisor of $\wt X_1$. 
\end{enumerate}
A limit $g$ with these properties is called the {\em refined limit}, or more precisely the once-refined limit of the (sub)sequence $f_{n}$ using rescaling parameters $\la_{n}$. \end{prop}
\begin{proof} In  the case $V$ is a smooth divisor, we proved all the key analytical estimates that imply this result in Section 6 of \cite{ip1} (see also  Section 3 of \cite{ip2}). The arguments used there to construct the refined limit are semi-local, in a neighborhood of $V$, and if set up right, easily 
extend to the case when $V$ is a normal crossing divisor. The main rescaling argument consists of two parts, first the construction of the refined limit by rescaling the target near $V$, and later on a further analysis of the properties of this refined limit. For the first part of the argument, we work separately in neighborhoods of depth $k$ stratum $V^k$ but away from the higher depth strata, and where usual Gromov compactness arguments apply; for the second part of the argument, we work in the necks, where the transition between these local models happens.  As this is one of the crucial steps in the construction of the relatively stable map compactification, we  include below the complete details of both of these arguments. 
\smallskip

\non {\bf Step 0. Preliminary considerations}. Assume for simplicity that the domains $C_n$ are smooth, otherwise work separately on each of their components, and that the smooth resolution $\wt C_0$ of their original limit $C_0$ is a stable curve with trivial automorphism group, otherwise first add enough marked points. 

If $C_0$ were smooth, the universal curve $\ov\U\ra \ov\M$ can be locally trivialized in a neighborhood of $C_0$. Otherwise, we can find a local trivialization only away from the nodes $D$ of $C_0$ (in the complement of a disjoint union of small balls in $\ov \U$ around these nodes). Around each node $x$ of $C_0$ we can choose local coordinates on the universal curve such that the domains $C_n$ are described by
\bear\label{cood.nodes}
z w=\mu_{n} (x)
\eear 
where $\mu_0(x)=0$, and $z, w$ are the two local coordinates on $C_0$ at the node $x$ (one on each local branch).  

On the target $X$, consider $U_\de$ the $\de$-tubular neighborhood of $V$ in $X$, and $U_\de^k$ the corresponding neighborhood \eqref{def.u.k} of $V^k$. 
The region $U_\de^k \setminus U_{\de'}^{k+1}$ lies in a neighborhood of $V^k$ but stays away from a neighborhood of the higher depth stratum $V^{k+1}$.

Denote by $P_n\subset C_n$ the collection of marked points of $f_n$, which include all the points in $f_n^{-1}(V)$, with their contact information recorded by $s$. As $n\ra \infty$, they converge to the marked points $P_0$ of $C_0$, but the original limit $f_0$ may only have a partially defined contact information to $V$, as described by Lemma \ref{f.dec.depth.k} and the two paragraphs preceding it. In particular, the restriction $f^k_0$ of $f_0$ to $C_0^k$ lies in $V^{k}$ where it has a well defined order of contact to the higher depth stratum, and 
\bear\label{dec.f.k.V}
f_0^{-1}(V^k\setminus V^{k+1})=(C_0^k \setminus R)\sqcup R^k \quad \text{and} \quad f_0^{-1}(V^k)= C_0^{\ge k} \cup R^{\ge k}
\eear
Here $C_0^k$ and $C_0^{\ge k}$ denote the union of depth $k$ (and respectively at least $k$) components of $C_0$ and $R^k$ the depth $k$ contact points of $C_0$. Note that a depth $k$ component may contain finitely many higher depth points, and these by definition are the contact points of that component. A priori some of these contact points may not be part of the special points $D\cup P_0$ of $C_0$ (and some special points may not be contact points). We next construct a decomposition of the domains $C_0$ and $C_n$ into pieces close to the depth $k$ stratum but which stay away from the higher depth stratum. 

For each point $x\in C_0$ denote by $B(x,\ep)$ the ball about $x$ of radius $\ep$ in the universal curve and by $B_n(x, \ep)$ its intersection with $C_n$. When $\ep>0$ is sufficiently small (depending on $x$) the intersection of $B(x,\ep)$ with $C_0$ is either a disk or a union of two disks depending whether $x$ is a smooth point of $C_0$ or a node  (which lifts to two smooth points of the resolution, indexing the two local branches of $C_0$ meeting at $x$). For each lift $\wt x$ of $x$ to the resolution $\wt C_0$, denote by 
$\gamma_\ep(\wt x)$ the oriented boundary of its corresponding disk, and by 
$\gamma_{n, \ep}(\wt x)$ the corresponding loop in the boundary of $B_n(x, \ep)$ for $n$ sufficiently large (depending on $\ep$).

For each contact point $x\in R$ and each small radius $\ep(x)>0$ depending on $x$, consider the disjoint union  $\ma\sqcup_{x\in R} B(x, \ep)$. For $n$ sufficiently large (depending on $\ep$), $C_n \setminus \ma \sqcup_ {x\in R} B(x,\ep) = \ma\sqcup_k C_n^k(\ep)'$ is a disjoint union of compact pieces $C_n^k(\ep)'$, one for each $k$, 
which limit as  $n\ra \infty$ to $C_0^k \setminus \ma \sqcup_ {x\in R} B(x,\ep) $ and then as $\ep\ra 0$ to $C^k_0 \setminus R$. Denote by 
\bear\label{D.approx.depth.k}
C_n^k(\ep) =  C_n^k(\ep)' \cup \ma \sqcup_{x\in R^k} B_n(x,\ep) 
\eear
the subset of $C_n$ obtained this way. 
Here $n$ could also be 0, corresponding to the limit case $C_0$. Note that 
\bear\label{D.c.n.k.converge}
C_n^k(\ep)\ma\longra_{n\ra \infty} C_0^k(\ep)= \l(C_0^k \setminus \ma \sqcup_{x\in R^{>k}} B_0(x,\ep)\r) \cup  \ma \sqcup_{x\in R^k} B_0(x,\ep) \ma\longra_{\ep\ra 0} C_0^k(0)=(C_0^k \setminus R) \cup R^k
\eear 
The boundary of $C^k_n(\ep)$ decomposes as $\bd_+\sqcup - \bd_-$ where 
\bear
\bd_+ C^k_n(\ep)=\ma\sqcup_{x\in R^{>k} \cap \wt C_0^{k}}  \ga_{n,\ep}(x) \quad\mbox{ while } \quad  
\bd_- C^k_n(\ep) =\ma\sqcup_{x\in R^k\cap \wt C_0^{<k}}  \ga_{n,\ep}(x). 
\eear 
The curve $C_n$ is obtained by joining each piece $C_n^k(\ep)$ along $\bd_-C_n^k(\ep)$ and $\bd_+C_n^k(\ep)$ to the lower and respectively higher index $j$ pieces $C_n^{j}(\ep)$.  

As $\de\ra 0$, the inverse images $f^{-1}_0(U^k_\de)$ shrink to $f^{-1}_0(V^k)$. By construction, 
\bear\label{def.depth.ep.str}
C^{\ge k}_0(\ep)= f_0^{-1}(V^k) \cup  \ma \cup_{x\in R^{\ge k}} B_0(x,\ep) 
\eear
also shrink to  $f^{-1}_0(V^k)$ as the radii $\ep(x)\ra 0$ for each $x\in R^{\ge k}$ (independent of choices at points in $R^{<k}$). Using the local expansion of $f_0$ at the contact points, for each $\de>0$ small, we can choose radii $\ep(x)>0$ at each $x\in R^{>k}$ (converging to 0 as $\de\ra 0$) so that $f_0^{-1}(U^{k+1}_{2\de})\subset C_0^{\ge k+1}(\ep)$ so  
\bear\label{est.f0.2}
f_0(C^{\le k}_0(\ep)) \subset X\setminus U_{2\de}^{k+1} 
\eear
independent of the choice of radii at the rest of the contact points ($x\in R^{\le k}$). Conversely, for any small choice $\ep(x)>0$ of radii at $x\in R^{>k}$, we can find a $\de>0$ so that \eqref{est.f0.2} holds, with $\de\ra 0$ as $\ep(x)\ra 0$. 

On the other hand, by construction  
\bear\label{est.f.33}
C^{k}_0(\ep)\subset C^k_0\cup  \ma \sqcup_{x\in R^{k}} B_0(x,\ep) 
\eear
and $f_0(C_0^k)\subset V^k$.  So for each $\de>0$ small, we can choose radii $\ep(x)>0$ at each $x\in R^k$ (converging to 0 as $\de\ra 0$) so that 
\bear\label{est.f0.3}
f_0(C^{k}_0(\ep)) \subset U_{\de/2}^{k}
\eear
independent of the choice of radii at the rest of the contact points ($x\in R\setminus R^k$), and conversely, for any small choice $\ep(x)>0$ at $x\in R^k$ we can find a $\de>0$ so that \eqref{est.f0.3} holds. 

\begin{rem} In effect, starting from the stratification $f_0^*\mathsf{depth}$ of $C_0$ pulled back from $X$, we constructed a decreasing family of stratifications $\mathsf{depth}_{\ep}$ of a tubular neighborhood of $C_0$ in the universal curve, which restrict to stratifications  of $C_n$ with closed strata \eqref{def.depth.ep.str}. Of course, $X$ also has a decreasing family \eqref{D.str.X.de} of stratifications  $\mathsf{depth}_{\de}$ whose closed strata are $U_\de^k$, and the estimates above compare these stratifications. \end{rem}

Fix  $0<\de<<1$ small enough so that for all $k\ge 0$ the depth $k$ piece $f^k$ of the original limit $f_0$ has very small energy in the $\de$-neighborhood of the higher depth stratum $V^{k+1}$: 
\bear\label{fix.eta}
\mbox{ the energy of $f^k$ in the region $U^{k+1}_\de$ is less than $\al_V/100N$}
\eear
where $N=\dim X$ is an upper bound for the depth $k$. Note that $\de$ can be chosen as small as desired (because the energy density is nonnegative). But we have assumed that $f_0$ has at least one component in $V$ thus the energy of the restriction of $f_0$ to the complement of $C^0_0$ is at least $\al_V$. 

\begin{lemma}\label{L.estimate.f.n} Fix $f_0$ and $\de>0$ as above. Then for any sufficiently small choice $\ep(x)>0$ of radii at the contact points $x\in R$, and all $n$ sufficiently large (depending on $\ep$)
\begin{enumerate}[(i)]
\item the restriction of $f_n$ to $C_n^k(\ep)$ has energy at most $\al_V/50N$ in the region  $U_\de^{k+1}$ for all $k\ge 0$ 
\item the restriction of $f_n$ to the complement of $C_n^0(\ep)$ has energy at least $2\al_V/3$
\end{enumerate}
For each $k\ge 0$ fixed, and any sufficiently small choice $\ep(x)>0$ of radii at the points $x\in R^{>k}$, there exists $m>0$ (depending on $\ep$ but smaller than $\de$), a choice of small radii $\ep(x)>0$ at the points $x\in R^{k}$, and two other $m_\pm>0$, such that for all $n$ sufficiently large (depending on $\ep$): 
\begin{enumerate}[(i)]
\item[(iii)] $f_n(C_n^k(\ep))$ lies in $U^k_{m}\setminus U^{k+1}_{m}$ thus lifts to $\F_k$ 
\item[(iv)]  $f_n(\bd_- C_n^k)$ lies  $U^k_{m}\setminus U^{k}_{m_-}$ while  $f_n(\bd_+ C_n^k)$ lies in $U^{k+1}_{m_+}\setminus U^{k+1}_{m}$.
\end{enumerate}
As the choices $\ep(x)\ra 0$ for $x\in R^{>k}$, then $m(\ep)\ra 0$ and we can arrange that $m_\pm(\ep)\ra 0$.
\end{lemma}
\begin{proof} This follows from the estimates \eqref{D.c.n.k.converge}-\eqref{fix.eta} using  the fact that $f_n\ra f_0$ in the stable map compactification as maps into $X$. The first two conditions are clear. For (iii), fix any small choice $\ep(x)>0$ at all $x\in R^{>k}$. Use first \eqref{est.f0.2} to find an $m>0$  so that $f_n(C_n^k(\ep))\subset X\setminus U^{k+1}_m$ for any $n$ large. We can arrange that $m(\ep)$ decreases to zero as $\ep\ra 0$. 
Because as $n\ra \infty$ the restriction of $f_n$ to $\bd_+ C_n^k(\ep)$ uniformly converges to the restriction of $f_0$ which lands in a neighborhood of $f_0(R^{>k})\subset V^{k+1}$, we can find an $m_+$ such that the second inclusion in part (iv) holds. 

Next use \eqref{est.f0.3} for to find $\ep(x)>0$ at $x\in R^k$ so that $f_n(C_n^k(\ep))\subset U^{k}_m$, which then gives (iii). As the restriction of $f_n$ to $\bd_- C_n^k(\ep)$ uniformly converge to the restriction of $f_0$ to $\bd_- C_0^k(\ep)$ which is mapped off $V^k$, we can find an $m_-$ such that the first inclusion in part (iv) holds.  \end{proof}

\non {\bf Step 1.  Constructing a refined limit}. Next we rescale $X$ around $V$ by a certain amount $\la_n$ to catch a refined limit of the sequence $f_n$. To find $\la_n$, we consider the energy $E(f_n,t)$ of $f_n$ inside the annular region $A(t, \de)$ of (\ref{defn.an.r}) around $V$ in $X$.  For each $n$ fixed, the energy $E(t)=E(f_n, t)$ is a bounded variation decreasing function in $t$,  which is equal to 0 when $t=\de$ and is at least $2\al_V/3$ when $t=0$ (because for $\ep$ small, the restriction of $f_n$ to the complement of $C_n^0(\ep)$ lies in the $U_\de$ by (iii) and has energy at least $2\al_V/3$ by (ii)).

So  for $n$ sufficiently large, there exists a unique  $t_n=|\la_n|\ne 0$ such that 
 \bear\label{def.la.n}
 \liminf_{t\ra t_n} E(f_n, t) \le \frac{\al_V}2 \leq  \limsup_{t\ra t_n} E(f_n, t)
\eear
i.e. the energy of $f_n$ in the annular region $A(|\la_n|, \de)$ is essentially  $\al_V/2$. 
Then $\la_n\ra 0$, because if they were bounded below by $\mu>0$ then in the limit the energy in the annular region $A(\mu, \de)$ of $f_0$ and thus of $f^0$  would be  $\al_V/2$ which contradicts (\ref{fix.eta}).

We next show that $R_{\la_n} f_n$ regarded as a sequence of $R_{\la_n}^*J_n$-holomorphic maps into $\x$ have a Gromov convergent subsequence. For each $n$, the almost complex structure $J_n$ rescales to an almost complex structure $J_n^0$ on the level one building $X_1$ and deforms to $J_n'$ on the total space of $\x \ra B$ by the symplectic sum construction. For sufficiently large $n$, all these almost complex structures are tamed by any one of the symplectic forms $\om^{\al}$ of Remark \ref{R.ep.on.F} (when $0<\al\ll 1$), because being tamed is an open condition. In particular, the restriction of $\om^{\al}$ to $X_\la$ tames $R_\la ^*J_n$ for large $n$. Because $R_{\la_n} f_n:C_n \ra \x$ represent the same homology class, they have uniformly bounded energy measured with respect the associated metric on $\x$, and therefore a Gromov convergent subsequence to a map $g:B\ra X_1$ whose image lands in the central fiber $X_1$ of $\x$, and this result is independent of the choice of $0<\al \ll 1$, as the symplectic form is only used to tame the almost complex structures. 



Let $H\subset X_1\subset \x$ denote the limit as $\la_n\ra 0$ of the image in $\x$ under $R_{\la_n}$ of the region $A(|\la_n|, \de)$. Then $H$ is an open neighborhood of the singular locus $W_1$ of $X_1$. In fact, its restriction $H_k=H\cap \F_k$ to $\F_k$ is the union of the upper hemisphere of $\F_k$ and the $\de$-neighborhood of the fiber divisor $F_k$ of $\F_k$.  

By \eqref{def.la.n}, the limit $g$ does not have enough energy in the neck region to have a nontrivial component in the singular locus (as these would carry at least energy $\al_V$). It also has a resolution $\wt g: \wt B \ra \wt X_1$ without any components in the infinity section. 

Let $g_k$ be the restriction of $g$ to $B_k=g^{-1}(\F_k\setminus \F_{k+1})$, see \eqref{level1}. Each $B_k\subset B$ is a punctured curve, let $S\subset B$ be the union over $k$ of their punctures, and then remove $S$ from $B_k$. By induction on $k$, Lemma \ref{L.estimate.f.n}(iii) and (iv) implies that $g_k$ is the limit of the restrictions of the subsequence $R_{\la_n}f_n$ to $C_n^k(\ep)$ as $n\ra \infty$ and then the radii $\ep(x) \ra 0$ at the points $x\in R^{>k}$: condition (iii) implies that their limit is defined on a subset of $B_k$ while (iv) implies the limit of their boundaries converge to the punctures  of $B_k$. In particular $g_k: B_k \ra \F_k$ refines $f^k:C^k\ra V^k \setminus V^{k+1}$ (regarded as a map into the zero section of $\F_k$). 


 So by construction $g$ is a refinement of $f_0$; it remains to check that $g$ is a strict refinement of $f_0$, i.e. $g$ has a nontrivial component whose depth is strictly less than that of its projection $f_0$. Otherwise, for each $k$, $g_k$ can only differ from $f^k$ by trivial components (since we already know that $g_k$ has no nontrivial components in the infinity section). But then this contradicts \eqref{def.la.n}, as the trivial components do not have enough energy in $H_k$, nor do the $f^k$: by the $f^k$'s lie in the zero section where together have at most $\al_V/100$ energy in the $\de$-neighborhood of the fiber divisor $F_k$ by  \eqref{fix.eta}.  
\end{proof}
Intrinsically Proposition \ref{L.lim.level.1} describes a refinement of the Gromov topology, in which a sequence convergent in the Gromov topology may have several limit points in the once-refined topology, appearing as the limit of different subsequences, and a priori different ways to rescale the original sequence. The next Lemma essentially rules out the latter possibility. 
\begin{lemma}\label{L.unique.lim} For each once-refined-convergent subsequence, its limit constructed by Proposition \ref{L.lim.level.1} is unique up to rescaling the level one of the building by an overall factor $\la \in \cx^*$. 

More precisely, assume $g$ and $g'$ are two refined limits of the same subsequence $f_n$ but obtained using parameters  $\la_n, \la_n' \ra 0 $. Then  (i) $|\la_n/\la_n'|$ has a limit $\ne 0, \infty$ and (ii) $g$ is equal to $g'$ after a reparametrization of the domain and rescaling the level one of the building by an overall factor $\al \in \cx^*$. 
\end{lemma} 
\begin{proof} 
First of all, note that this notion of convergence is preserved under adding converging sequences of extra marked points to the domains or rescaling the parameters $\la_n$ by factors $\al_n$ as long as $\al_n\ra \al \ne 0, \infty$. In particular, if $g$ is the  refined limit of the sequence $f_n$ obtained by rescaling by $\la_n$ and 
$\al\ne 0, \infty$ then $\al g$ is the refined limit of the sequence $f_n$ obtained by rescaling by $\al_n\la_n$. 

Assume $g':C'\ra X$ is another refined limit of the same sequence $f_n$, obtained by rescaling by 
$\la_n'$ and adding a different collection $p_n'$ of $k'$ marked points to the domains to kill their automorphisms. By choosing a common lift $(C_n, p_n, p_n')\in \ov\M_{k+k'}$ and passing to the limit, we conclude that $(C, p, p')= (C', p,p')$ up to reparametrizations. So we can assume that both sequences had the same collection of marked points added to begin with to kill the automorphisms and drop that from the notation. Then $C=C'$ so the limits are defined on the same domain (e.g. the fiber of the universal curve), and their graphs in $\ov\U_k\ti X$ are embedded $J$ holomorphic curves. 

Next use energy considerations to show that the ratio $|\la_n'/ \la_n|$ has a limit $\al \ne 0, \infty$. First of all, since both $g$, $g'$ are strict refinements of $f_0$, they have a nontrivial component meeting the upper hemisphere, so by rescaling each one separately we can assume their energy in the upper hemisphere is essentially $\al_V/2$ in the sense  of (\ref{def.la.n}). Since both have very small energy near the singular divisor by (b) (coming only from trivial components) and the convergence is also in energy density, then $|\la'_n/\la_n|$ is uniformly bounded away from 0 and $\infty$ and therefore after passing to a further subsequence and rescaling one of the sequences by a fixed constant we can assume $\la_n'/\la_n\ra 1$. Now we can use the fact that the graphs of $R_{\la_n} f_n$ converge in 
Hausdorff distance to that of $g$ to conclude that the graphs of $g$ and $g'$ have Hausdorff 
distance zero. But their domains were equal, and the graphs are embedded $J$-holomorphic curves, thus they must be equal up to a reparametrization of the domain. \end{proof}
\non{\bf Step 2: Properties of the refined limit.}  We next want to understand the behavior of the refined limit $g:C\ra X_1$ constructed by Proposition \ref{L.lim.level.1} around the singular divisor $W_1$ where the pieces of the building are joined together; for that we restrict to the neck regions (both of the domain and of the target) where we use the local models (\ref{cood.nodes}) on the domain and (\ref{eq.xy=l}) in the target. 

We already know that the refined limit $g:C\ra X$ has a resolution $\wt g:\wt C\ra \wt X_1$ which has no components in the infinity divisor.  Lemma \ref{f.dec.depth.k} applies to this resolution to produce a further lift $\wh g$ recording the contact information of each depth $k$ piece of $\wt g$ to the higher depth stratum of the level one building $X_1$ with respect to its total divisor $D_1$.

A convenient way to record the contact information of the refined limit is via the semi-local models for the target $X_1$ described in \S \ref{rem.def.multilevel}. For each depth $k$ component $\Si$ of $C$, index the stratum of $\wt X_1$ containing the lift $\wt g(\Si)$ by the data $I_\Si$, $l_\Si$ and $\ep_\Si$ where (a) $I_\Si$ records the stratum of $X$ containing the projection of $\wt g(\Si)$ under the collapsing map $X_1\ra X$; (b) $l_\Si:I_\Si\ra \{0, 1\}$ is the multilevel map recording the level of $\wt g(\Si)$ in each direction and (c) $\ep_\Si:I_\Si \ra \{ +, 0\}$ is the multi-sign map recording the stratum of the total divisor containing $\wt g(\Si)$. Note that because $\wt g$ has no components in the infinity section, here $\ep_\Si(i) \ne -$ for all 
$i\in I_\Si$. Moreover, $I^\infty_\Si\ma=^{def} \ep_\Si^{-1}(\pm)$ indexes the directions in which 
$ \wt g(\Si)$ does not have a well defined contact information to the total divisor;  its cardinality records the depth of $\Si$ in $\wt g$, while the cardinality of $ I_\Si$ records the depth of the projection to $X$. 

Next, each contact point $x\in \Si$ similarly comes with an indexing set  $I(x)$, a multilevel map $l(x)$ and a multisign map $\ep(x)$ refining those associated to $\Si$, together with sequence of multiplicities $s(x): I(x)\setminus  I^\infty_\Si\ra \N$ recording the contact multiplicity 
of $\wh g$ at $x$ to the higher depth stratum of $\wt X_1$ indexed by $I(x)$, $l(x)$ and $\ep(x)$. Denote by $I^\infty(x)= I^\infty_\Si$ the directions in which the contact multiplicity is undefined (because $\wt g(\Si)\equiv  0$ in those directions), and by $I^\pm(x)$ and respectively $I^0(x)$ the partition of the remaining directions induced by $\ep(x)$.   

In local coordinates $z$ on $\Si$ at each contact point $x$, and normal coordinates $u_{i, l}\in N_i$  to the zero divisor in each level $l$ of the target, the resolution $\wh g(z)$ has an expansion:
\bear
\label{exp.u.v.x}
(u_{i, l_i(x)})^{\ep_i(x)}=\begin{cases}
a_i (x) z^{s_i(x)}+ O(|z|^{s_i(x)}) &  \text{ for all $i\notin I^\infty(x)$}
\\
0& \text{ for all $i\in I^\infty(x)$}
\end{cases}
\eear
where the leading coefficient  
\bear
0\ne  a_i(x)\in   N_i^{\ep_i(x)} \otimes (T_{x}^* \Si )^{s_i(x)} \quad \mbox{ for all } i\notin I^\infty(x)
\eear
With this notation, the conclusion of Proposition \ref{L.lim.level.1} can be strengthened as follows: \begin{prop}\label{L.lim.level.2} 
 Assume $J_n\ra J_0$ is a
sequence in $\JV(X, V)$, and consider a sequence  $f_n:C_n\ra X$   of maps in 
$\M(X,V)$ whose level zero limit $g_0:C_0\ra X$ has at least one nontrivial component in $V$. Then there exists a sequence $\la_n\ra 0$ of rescaling parameters such that,  after passing to a subsequence, $R_{\la_n} f_n$  converge to a limit $g_1:C_1\ra X_1$ into a level one building  with the following properties: 
\begin{enumerate}[(a)]
\item $g_1:C_1\ra X_1$ is a strict refinement of  $g_0:C_0\ra X_0$, without any nontrivial components in the singular locus $W_1$ of $X_1$;  
\item $g_1$ has a resolution $\wt g_1: \wt C_1 \ra \wt X_1$ whose contact information to the total divisor of $\wt X_1$ is described by (\ref{exp.u.v.x}); these maps fit in the diagram:
\bear\label{g.ref.f00} 
\xymatrix{
\wt C_1\ar[r]^{\xi}\ar[d]_{\wt g_1}&C_1\ar[r]^{\st}\ar[d]_{g_1} & C_0\ar[d]^{g_0}
\\
\wt X_1\ar[r]^{\xi}&X_1\ar[r]^p& X_0
}
\eear
\item for any intermediate curve $C_1\ra C'\ra C_0$ all the contact points of $\wt g_1$ descend to special points of $C'$, and moreover  
\begin{itemize}
\item $s(x_-)= s(x_+)$ and $\ep(x_-)= -\ep(x_+)$  for each node  $x_-=x_+$ of  $C'$; 
\item $s (x)= s(x_n)$ and $\ep(x)=\ep(x_n)$ for each marked point $x\in C'$ which is the limit of marked points $x_n$ of $C_n$; 
\end{itemize}
(whenever both sides are defined). Furthermore, $C$ is obtained from $C_0$ by inserting strings of trivial components $B_x$ (broken cylinders) either between two branches $x_\pm$ of a node $x$ of $C_0$ or else at a marked point $x$ of $C_0$; 
\end{enumerate}
For each once-refined-convergent subsequence, its limit $g_1:C_1\ra X_1$ is unique up to the $\cx^*$ action on the level one of the building.
\end{prop}
\begin{proof} Part (a) and (b) have already been proven, thus it remains to show (c).  First of all, we can decompose the domain $C$ into two pieces, $C^0$ and $B$ where $C^0$ consists of nontrivial components, and  $B$ consisting of components that get collapsed under the two maps $C\ra C_0$ and $X_1\ra X$.  Then each connected component $B_i$ of $B$ is an unstable genus zero curve (bubble tree) with either one or two marked points, corresponding to a special fiber of $C\ra C_0$ over either: 
\begin{enumerate}[(i)]
\item  a non special point in the case $B_i$ has only one marked point;  
\item a node in the case $B_i$ has two marked points and both belong to $C^0\cap B$ or
\item a marked point in the case $B_i$ has two marked points, but only one belongs to $C^0\cap B$ while the other is a marked point of $C$.  
\end{enumerate}
We have the same description for the special fibers of $C\ra C'$ and of $C'\ra C_0$ for any intermediate curve 
$C\ra C' \ra C_0$. 

For any point  $x\in C_0$, choose local coordinates at $x$ on the universal curve of the domains (containing $C_0$) and normal coordinates to $V$ in the target $X$ at $p$ as above (where $p=f_0(x)$ is a depth $k\ge 0$ point of $V$). Using the notations of  the proof of Proposition \ref{L.lim.level.1},  for $\ep, \de>0$ sufficiently small and $n$ sufficiently large, the image under $f_n$ of $B_n(x,\ep)$ is mapped in the $\de$ neighborhood of $p$ in $X$ but away from the $\de$ neighborhood of the depth $k+1$ stratum. Furthermore, since $f_n^{-1}(V)=P_n$ then $f_n$ maps $B_n(x,\ep)\setminus P_n$ into the annular region 
$O_\de$
\bear\label{u.s1}
0<|u_i|<\de \quad \mbox{ for all } i\in I
\eear
of $N_{V^k}$ around $p$ (away from the higher depth stratum). This region is homotopic to $(S^1)^k$, one factor 
for each branch $i$ of $V$ at $p$. 

If we also fix to begin with global coordinates on all bubble components in $B$,  we have a similar story for any point  
$x$ in any one of the intermediate curves $C'$: $f_n$ takes a sufficiently small punctured neighborhood $B_n(x,\ep)\setminus P_n$ into the annular region (\ref{u.s1}), where now $B_n(x, \ep)$ denotes the intersection of $C_n$ with the ball $B(x,\ep)$ about $x$ in the local model for domains containing their intermediate limit $C'$. Topologically, $B_n(x, \ep)$ is either a disk or an annulus, depending whether $x$ is a smooth point or a node of $C'$. In particular, for every lift 
$\wt x$ of $x$ to the resolution of $C'$ and thus to  $\wt C$, the corresponding boundary loop $f_n(\ga_{n,\ep}(\wt x))$ has a well defined winding number $s_{i, n}(\wt x)$ about the branch $i$ of $V$ at $p$.  

Furthermore, for $\ep>0$ sufficiently small, $B_n(x,\ep)$ contains no points from $P_n$ if $x$ is not a marked point of $C'$  and otherwise it contains precisely one point $x_n\in P_n$ which limits to $x$ as $n\ra \infty$.  This implies that, for $n$ sufficiently large,  and all $i\in I$, 
\begin{enumerate}[(i)]
\item if $x$ is not a special point of $C'$ then the winding numbers $s_{i, n}(\wt x)=0$;
\item if $x_\pm\in \wt C$ correspond to a node $x$ of $C'$ then $s_{i, n}(x_+)+ s_{i, n}(x_-)=0$. 
\item if $x\in C'$ is the limit of the marked points $x_n\in C_n$ then $s_{i, n}(\wt x)= s_i(x_n)$; 
\end{enumerate}
But the winding numbers $s_{i, n}(\wt x)$ of $f_n$ are related to those of the refined limit $g$. This is 
simply because the winding numbers of $f_n$ agree with those of the rescaled maps $R_{\la_n}f_n$, and these converge uniformly on compacts away from the nodes of $C$ to the refined limit $g:C\ra X_1$, which has well defined winding numbers about the zero section in all directions $i\notin I^\infty$. Since the loops $\ga_{n, \ep}(\wt x)$ stay away from all the nodes of $C$ (for $\ep$ sufficiently small) then for $n$ sufficiently large: 
\bear\label{match.sn=s}
\mbox{ $s_{i, n} (\wt x)=\ep_i(\wt x) s_i(\wt x)$ for all  $i\notin I^\infty$}
\eear
because of the expansion (\ref{exp.u.v.x}) of $g$ at $\wt x$.  Note that this give us no information about the winding numbers in the direction of $I^\infty$, where the winding numbers of $g$ are undefined. 

In particular, any contact point $\wt x$ of $g$ with $W_1\cup V_1$ or any of its strata has $s_i(\wt x)>0$ in  some 
direction $i$ which rules out case (i):  if $\wt x\in \wt C$ descends to a non special point $x$ of $C'$,  then $s_{i, n}(\wt x)=0$ by (i) which contradicts (\ref{match.sn=s}).  

In case (ii), for any node $x$ of $C'$  then  $s_{i, n}(x_-)+ s_{i, n}(x_+)=0$ and so  (\ref{match.sn=s}) implies that 
\best
s_i(x_-) = s_i(x_+)  \mbox { and  } \ep_i(x_-)= -\ep_i(x_+)\quad  \mbox { for all } i\notin I^\infty(x)\ma=^{def} 
I^\infty(x_-)\cup I^\infty(x_+)  
\eest
as both sides are well defined for such an $i$.  If $x$ is the limit of contact points $x_n$ of $f_n$ to $V$ then (\ref{match.sn=s}) implies that 
\best
s_i(x)= s_i(x_n)  \mbox { and  } \ep_i(x)= \ep_i(x_n)\ \quad \mbox { for all } i\notin I^\infty(x)
\eest  
So for example $x$ is a contact point of $g$ with the zero divisor $V_1$ if and only if $x_n$ was one for $f_n$ to $V$. 

Finally, this discussion implies that there are no components of $B$ with just one marked point (otherwise, contracting such component would give a curve $C'$ and a non special point $x$ on it which is impossible as case (i) is ruled out for all intermediate curves $C'$).  Therefore all the components of $B$ have precisely 2 special points, which must be the contact points with the zero and the infinity divisor in their fiber, and thus are indeed trivial components as in Definition \ref{def.triv.comp}. Furthermore, the only special fibers of $C\ra C'$ or $C'\ra C_0$ are strings of trivial components (broken cylinders). \end{proof}

\begin{rem}\label{R.intro.I} Proposition \ref{L.lim.level.2} and the discussion preceding it implies that each node $x$ of $C_0$, and more generally each node $y$ of any intermediate curve 
$C\ra C'\ra C_0$ which projects to $x$ comes with a partition of the original indexing set $I(x)$ of the stratum containing $f_0(x)$  into $I^\infty(y)$,  $I^0(y)$ and then $I^\pm(y)$ induced by the contact information of the refined limit $g:C\ra X_1$ to the total divisor.  $I^\infty(y)=I^\infty(y_+)\cup I^\infty(y_-)$ records those directions in which at least one of the local branches at $y_\pm$ of the resolution $\wt g$ lies in the total  divisor,  while $I^0(y)$ records the directions in which the $i$'th coordinates of both $\wt g(y_\pm)$ are nonzero, so both branches stay away from the total divisor in those directions. The  remaining directions $i\in I^\pm(y)$ come with a multiplicity $s_i(y)=s_i(x)>0$ and opposite signs $\ep_i(x_+)=\ep_i(y_+)=-\ep_i(y_-)=-\ep_i(x_-)\ne 0$ recording the two opposite sides of the level one building from which the two branches of $g$ come into the singular divisor, and also two leading coefficients $a_i(y_\pm)\ne 0$, intrinsically elements
\bear\label{a.v}
a_i(y_+) \in N_{f(x), i}^{\ep_i(x_+)} \otimes (T_{y_+} C )^{s_i(x)}\quad  \mbox{ and } \quad a_i(y_-) \in N_{f(x), i}^{\ep_i(x_-)} \otimes( T_{y_-} C)^{s_i(x)} 
\eear
for all $i\in I^\pm(y)$, where $N_{p, i}$ is the branch of the normal bundle to $V$ at the point $p$ in the direction $i$. 

Note that \eqref{l-l=ep} implies that for each node $y$ that is mapped to the singular divisor, and each fixed direction $i$ with $\ep_i(y)\ne 0$, one the two branches $y_\pm$ must be in level zero in that direction (the one for which $\ep_i$ evaluates to $1$) while the other one must be on level one.
\end{rem} 
\begin{ex} In the situation of Example \ref{ex.4dim} we could have two nodes of the domain mapped to $p$, one between the components $X$ and $\F_2$ while the other one between  the two $\F_1$ components, but no node at $p$ between say $X$ and $\F_1$ (as the branches of $g$ would not be on opposite sides of the singular divisor in all local directions).  
In the first case the node is between level 0 and level 1 (really local levels (0,0) and (1,1)), while in the second case it is between two level 1 floors,  or more accurately between a local level (0,1) and (1,0) floor. 
\end{ex} 
\subsection{Once-refined limits and the refined matching condition} Assume next that the once-refined limit constructed by  Proposition \ref{L.lim.level.2} has no components in the total divisor $D_1=W_1\cup V_1$. In this case, the rescaling process terminates in one step, and $g=g_1$ is the refined limit of the subsequence $f_n$. Then $g:C_1\ra X_1$  has a resolution $\wt g:\wt C_1\ra \wt X_1$ which is an element of the universal moduli space  
\bear\label{D.wt.g}
\wt g\in \M(\wt X_1, \wt W_1\cup \wt V_1)\ra \JV(X, V)
\eear
over the parameter space $\JV(X, V)$. The domain $\wt C_1$ is a resolution of $C_1$ obtained by resolving all its nodes $D$. The combined attaching map identifies pairs of marked points of $\wt C_1$ to produce the nodes of $C_1$, and simultaneously attaches the targets together to produce the singular locus $W_1$. The resolution $\wt g$ is essentially unique (up to reordering). 

The attaching map $\xi:\wt X_1\ra X_1$, when restricted to a depth $k\ge 1$ stratum is  a degree $2^k$ cover of  $W_1$. At each node $x_-=x_+$ of $C_1$ mapped into this stratum we also have a partition of the $2k$ normal directions  $N_{V^k}\oplus N^*_{V^k}$ into two length $k$  dual indexing sets $I_W(x_-)$ and $I_W(x_+)$  that record the two opposite local pieces of $\wt X_1$  containing the two local branches of $g$ at that node. According to our setup, we encode the indexing set of the branches of the divisor together with the contact multiplicities as decorations on the corresponding half edges of the dual graph cf Remarks \ref{top.info} and \ref{R.s.intrinsic}. 

\begin{rem} \label{match.naive} Let $\wt s$ denote the dual graph associated to the resolution $\wt g:\wt C_1\ra \wt X_1$, decorated by its (partial) contact information to the total divisor $D_1$. From it, we can read off the dual graphs of all the vertical arrows in \eqref{g.ref.f00}, as well as those of the original sequence $f_n$ if we assume their domains $C_n$ were smooth, see Remark \ref{R.dual.defm}. Denote by $s_\pm$ the contact information of $\wt s$ associated to the nodes $D$, i.e. $s_\pm(x)= s(x_\pm)$ for all $x\in D$, and by $s$ the dual graph obtained by contracting all the edges $D$ of $\wt s$. 

Then Proposition \ref{L.lim.level.2}(c) implies that  the image of the resolution \eqref{D.wt.g}  under the evaluation map (\ref{or.ev}) at the pairs $\{x_\pm\}_{x\in D}$ of marked points giving the nodes: 
\bear\label{ev.nodes}
\ev_{D}:\M_{\wt s} (\wt X_1, \wt W_1\cup \wt V_1)\longra  (W_1)_{s_+}\ti (W_1)_{s_-}
\eear
must land in the  diagonal $\Delta$. We call these the {\em naive matching conditions} because when $s_{\pm}$ contains depth $\ge 2$ points, the dimension of this stratum is in general bigger than the dimension of $\M_s(X, V)$!
\end{rem} 
\begin{lemma}\label{dim.diff1} The difference between the expected dimensions is 
\best
\dim \M_s(X, V) - \dim \ev_{s_\pm}^{-1}(\De) =  2 \sum_{x\in P(s_+)}  (1- k(x)) 
\eest 
where $P(s_+)$ denotes the collection of marked points associated with the sequence $s_+$. 
\end{lemma} 
\begin{proof} This is a simple adaptation of  the calculations of \cite{ip2}. The expected dimension of 
$\ev_{s_\pm}^{-1}(\De) $ is 
\best
\dim \ev_{s_\pm}^{-1}(\De)  &=& \dim \M_{\wt s} (\wt X_1, \wt W_1\cup \wt V_1) -\dim (W_1)_{s_+}
\eest
so the difference is
\best
\dim \M_s(X, V) - \dim \ev_{s_\pm}^{-1}(\De) = \dim \M_s(X, V) -  \dim \M_{\wt s} (\wt X_1, \wt W_1\cup \wt V_1) + \dim (W_1)_{s_+}
\eest
where 
\best
 \dim \M_s(X, V) &=&2c_1(T X)A_{s} + (\dim X-6) \frac{\chi} 2 +2\ell(s)-2A_{s} V
 \\
\dim \M_{\wt s} (\wt X_1, \wt W_1\cup \wt V_1)  &=&2c_1(T\wt X_1)A_{\wt s} + (\dim X-6) \frac{\wt \chi} 2 +2\ell(\wt s)-2A_{\wt s}\wt W_1-2A_{\wt s}V_1
\eest
But  $\chi=\wt \chi -2 \ell(s_+)$, $\ell(\wt s)= \ell (s) + 2\ell(s_+)$  and   $A_sV = |s|= A_{\wt s} V_1$,  while Lemma 2.4 of \cite{ip2} adapted to this context gives $
c_1(TX) A_s = c_1(T\wt X_1)A_{\wt s}-2A_{\wt s}\wt W_1 $ so
\best
&\dim \M_s(X, V) -  \dim \M_{\wt s} (\wt X_1, \wt W_1\cup \wt V_1) = (\dim X -2)\ell (s_+)
\eest
On the other hand, since the image under the evaluation map of each depth $k$ point lands  in a codimension $2k$ stratum of $X_1$  then  
\best
\dim (W_1)_{s_+}= \sum_{x\in P(s_+)} (\dim X-2k(x))
\eest
so the difference in dimensions is exactly as stated. \end{proof} 
\smallskip

The refined limit $g$ is  well defined only up to an overall rescaling parameter $\la\in \cx^*$ that acts on the level one of the building; 
this $\cx^*$ action induces an action on the restriction of the universal moduli space of maps into 
$\wt X_1$ to the parameter space $\JV(X, V)$, and the evaluation map \eqref{ev.nodes} descends to the quotient: 
\bear\label{ev.nodes.cx}
\ev_{D}: \M_{\wt s} (\wt X_1, \wt W_1\cup \wt V_1)/\cx^*\longra  (W_1)_{s_+}\ti (W_1)_{s_-}
\eear
Even after dividing by the $\cx^*$ action, Lemma \ref{dim.diff1} still implies that if we want to construct a relatively stable map compactification for a singular normal divisor (with depth $k\ge 2$ pieces), then we need some  refined matching condition, otherwise the boundary stratum is larger dimensional than the interior. Luckily, the existence of such refined compactification follows after a careful examination of the arguments in \cite{ip2}.   
\smallskip

It turns out that when $k\ge 2$, not all the maps $g:C\ra X_1$ (without components in $D_1$) whose resolution 
$\wt g$ satisfies the naive matching condition can occur as limits after rescaling of maps $f_n:C_n\ra X$ in $\M_s(X, V)$. 
To describe those that occur as limits,  we use the results of Section 5 of \cite{ip2}. For each node $x$ of $C$, we work in the local models (\ref{cood.nodes}) on the domain and (\ref{eq.xy=l}) in the target, using the local coordinates as described. 
The results of Proposition \ref{L.lim.level.m} can be strengthened as follows:
\begin{lemma}\label{L.lim.level.refined} Consider  $f_n:C_n\ra X$  a sequence of maps in $\M(X,V)$ as in Proposition \ref {L.lim.level.m}, and further assume that its refined limit $g:C\ra X_1$ constructed there has no components in the singular divisor $W_1$. Then for each node $x_-=x_+$ of $C$ which is mapped to $W_1$ we have the following relation: 
\bear\label{la=mu}
\lim_{n\ra \infty}\frac{ \la_n}{ \mu_n(x) ^{s_{i}(x)} }=a_i(x_-) a_i(x_+) \quad \mbox{ for each }  i\in I^\pm(x)
\eear
where $\mu_n(x)$ are the gluing parameters (\ref{cood.nodes})  describing $C_n$ at $x$ in terms of $C$, $\la_n$ is the sequence of rescaling parameters in the target. \end{lemma}
\begin{proof}  For each node $x$ of $C$ that is mapped to a depth $k(x)\ge 1$ stratum of $W_1$ with matching multiplicities $s(x)$, we work in the local models described above, where we can separately project our sequence  into each direction $i\in I(x)$, where the local model is that of rescaling a disk around the origin, as explained in more details in the next section. The projections are now maps into a rescaled family of disks, precisely the situation to which Lemma 5.3 of \cite{ip2} applies to give (\ref{la=mu}) in each direction, after using the expansions  (\ref{exp.u.v.x}).   \end{proof}

\begin{rem}\label{L.lim.level.model}  The local model of $\M_s(X, V)$ near a limit point $g:C\ra X_1$ which has no components in $W_1$ is described by tuples $(\wt f, \mu, \la)$ satisfying a {\em refined matching condition} (\ref{match.enh}). Here  $\wt f:\wt C\ra \wt X_1$ is an element of  $\ev_{D} ^{-1}(\Delta)\subset \M_{\wt s}(\wt X, \wt W_1)$ satisfying the naive matching condition  along $W_1$ described in Remark \ref{match.naive},  $\la \in \cx$ is the gluing (rescaling) parameter of the target, and $\mu \in \ma \bigoplus_{ x\in D} 
T^*_{x_-} C\otimes T^*_{x_+}C$ is the gluing parameter of the domain such that they also satisfy the condition
\bear\label{match.enh}
a_i(x_-) a_i(x_+) \mu(x) ^{s_i(x)} =\la  \qquad \mbox{ for all }  i\in I^\pm(x)
\eear
at each node $x\in D$ of  the domain, where $a_i(x_\pm)$ are the two leading coefficients (\ref{a.v}) of $\wt f$ in the $i$'th normal direction $i\in I^\pm(x)$ at the node $x\in D$.   
\end{rem}
Intrinsically, the gluing parameters $\mu$ in the domain are sections of the bundle  
\best
\ma\bigoplus_{x\in D} \L_{x_-} \otimes \L_{x_+}
\eest
while the gluing parameter $\la$ in the target  is naturally a section of the bundle $N\otimes N^*\cong \cx$. The condition (\ref{match.enh}) can be expressed as $k(x)$ conditions on the leading coefficients:
\best
a_i(x_-) a_i(x_+) = \la \cdot \mu(x) ^{-s_i(x)}
\eest
at each node. If we fix a small $\la \ne 0$, the existence of a  $\mu(x)\ne 0$ satisfying these relations imposes a $2k(x)-2$ dimensional condition on the leading coefficients, exactly what was missing in Lemma \ref{dim.diff1}. 

\begin{rem}\label{R.match.enh.log} The refined matching conditions become linear if we take their log: 
\bear\label{match.enh.log}
\log a_i^+(x)+\log  a_i^-(x) = \log \la - s_i(x) \log \mu(x)
\eear
which makes the transversality of this condition easier to prove, and also hints to the connection with log geometry. Here $\log$ is the appropriate extension of the map $\log: \cx^* \ra \R \ti S^1$ defined by $\log z = \log |z|+ i\arg z$  to the intrinsic bundles in question.
\end{rem} 

\begin{rem} \label{R.cx} There is a second  way to read (\ref{la=mu}).  Because the leading coefficients $a_i(x_\pm) \ne 0, \infty$, then for each node $x$ of $C$  its contact multiplicities must be independent of $i$, i.e.  
\best
s_i(x)=s_j(x) \mbox{ for all } i, j \in I^\pm(x). 
\eest
Moreover the resolution $\wt g_1:\wt C_1 \ra \wt X_1$ of the limit $g_1$ must satisfy the {\em refined matching condition} at each node $x\in D$ i.e. the image of $\wt g_1$ under  the refined evaluation map 
\bear\label{enh.ev.match}
\Ev_{x_\pm}:\M_{\wt s} (\wt X_1, \wt W_1) \ra  \P_{s(x)} (N  W_{J_-(x))}) \ti  \P_{s(x)} (N W_{ J_+(x)}) 
\eear
lands in the antidiagonal 
\best
\De^\pm= \{\;  ( [a_i], [a_i^{-1}]) \;   | \;  a_i\ne 0, \infty \mbox{ for } i \in I^\pm(x) \;  \} 
\eest 
\end{rem}
Unfortunately, in the presence of a depth $k\ge 2$ point, we cannot rescale the target such that the limit $f$ has no components in the singular locus $W_1$. The most we can do is to make sure it has no nontrivial components there, but at the price of getting several trivial components stuck in the singular divisor. Below are a couple of simple examples that illustrate this behavior.

\begin{ex}\label{neck.1} Consider the situation of Example \ref{ex.4dim.b}. Assume $f$ is a fixed stable map into $X$ which has two contact points $x_1$ and $x_2$ with $V$, both mapped into the singular locus $p$ of $X$, but such that $x_1$ has multiplicity (1,1) while $x_2$ has multiplicity (1, 2) to the two local branches  of $V$ at $p$. This means that in local coordinates $z_i$ around $x_i$ on the domain and  $u_1, u_2$ on the target around $p$, the map $f$ has the expansions 
\best
f(z_1)= (a_{11}z_1,\; a_{21} z_1)\qquad  f(z_2)= (a_{12}z_2, \;a_{22} z_2^2)
\eest
with finite, nonzero leading coefficients $a_{ij}$. Now add another marked point $x_0$ to the domain. As either $x_0\ra x_1$ or $x_0\ra x_2$ a constant component of $f$ is falling into $p$. 

Let's look at the case $x_0\ra x_1$.  Assume $x_0$ has coordinate $z_1=\ep$ so $f(x_0)= (a_{11}\ep, \;a_{21} \ep)$  and  $\ep \ra 0$.  Following the prescription of \cite{ip2}, we need to rescale the target by $\la= \ep$ to catch the constant component falling in. So in coordinates $u_{11}=u_1/\la$ and $u_{21}= u_2/\la$ we get 
\best
f_{1\la} (z_1)= (a_{11}z_1/\la, a_{21} z_1/\la)\qquad  f_{2\la}(z_2)= (a_{12}z_2/\la, a_{22} z_2^2/\la)\
\eest
In the domain, when we let $w_1= z_1/\ep=z_1/\la$ then $f_{1\la}$ converges to a level one nontrivial component which in the coordinates $u_{11}$ and $u_{21}$ has the expansion
\best
f_1(w_1) = (a_{11}w_1,\; a_{21} w_1)
\eest
This component lands in $\F_2$ and contains  the marked points $x_0$ and $x_1$ (with coordinates $w_1=1$ and respectively $w_1=0$), so it is the original component of $f$ that was falling into $p$ as $x_1\ra x_0$. 

\begin{figure}[h]
\centering
\includegraphics[height=4cm]{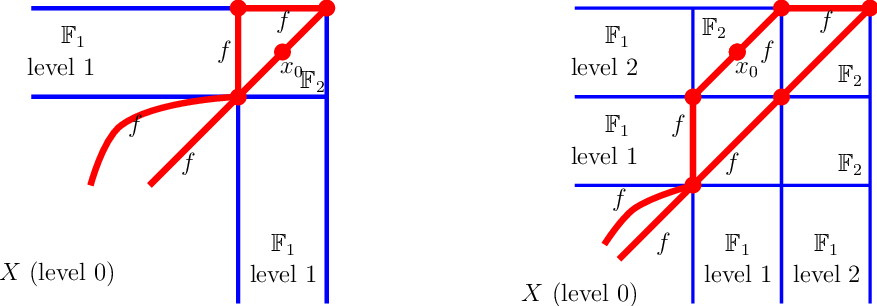}
\caption{(a) Limit as $x_0 \ra x_1$ \qquad  (b) Limit as $x_0 \ra x_2$}\label{F.4}
\end{figure}

But when we rescale the target by $\la=\ep$, the other piece of $f$ at $x_2$ also gets rescaled, and limits to trivial components in $\F_2$.  If we rescale the domain by $w_{21}=z_2/\sqrt \ep$ then $f_{2\la}(w_{21})= (a_{12}w_2/\sqrt \ep ,\; a_{22} w_2^2)$ converges to a trivial map in the neck
\best
f_{21} (w_{21})= (\infty,\;  a_{22} w_{21}^2)
\eest 
while if we rescale the domain by $w_{22}= z_2/\ep$ then  $f_{2\la}(w_{22})= (a_{12}w_{22} , a_{22} w_{22}^2/\ep )$ also 
converges to a trivial map in the zero divisor
\best
f_{21} (w_{22})= (a_{12}w_{22},\; 0)
\eest 
Putting these together, we see that the limit of $f$ as $x_0\ra x_1$ consists of a map into a level 1 building, which has one component $f$ on level zero and 3 components $f_1$, $f_{21}$ and $f_{22}$ on level one (all mapped into $\F_2$). But  only $f_1$ is a nontrivial component while the other two components are trivial, one of them mapped to the singular locus between $\F_2$ and $\F_1$ while the other one is mapped into the zero divisor of $\F_2$, see Figure \ref{F.4}(a).  

One can also see what happens as $x_0\ra x_2$. Then the limit is a map into a level 2 building, which now has 5 rescaled components, only one of them nontrivial (the one containing $x_0$). The piece of $f$ containing $x_1$ now gives rise to two trivial components $f_{11}$ and $f_{12}$ one on level (1,1) and the other on  level (2,2) piece $\F_2$. On the other hand, the piece of $f$ containing $x_2$ gives rise to  three components, the first one a trivial component in the neck between the level 1 piece 
$\F_1$ and $\F_2$, the next one a nontrivial component mapped to the level (1,2) piece $\F_2$ and the last piece is a trivial component mapped into the zero divisor of the level (2,2) piece $\F_2$, see Figure \ref{F.2}(b). 

Note that the only nontrivial component in this case lands in level one with respect to one of the directions, but also in level 
two with respect to the other direction, so we have a nontrivial component in each level. 
\end{ex} 
\par

\begin{ex} \label{neck.2}  If instead we  consider  Example \ref{ex.4dim}, where the two normal  directions at $p$ are globally independent, then the limit in the case $x_0\ra x_1$ would look just the same. However, the limit when $x_0\ra x_2$ would have fewer components, as we now can rescale independently  by two factors $\la_1=\ep$ and $\la_2=\ep^2$ getting a level (1,1) building.  The limit then has only 3 pieces on level one, all of them mapped to $\F_2$, but again only the piece containing $x_0$ is nontrivial. The other two pieces come from rescaling the piece of $f$ containing $x_1$ so they are both trivial, the first one in the neck between $\F_2$ and the level (1, 0) piece $\F_1$ while the other component lands in the zero divisor of the level (1,1) piece $\F_2$.  
\end{ex}
\begin{figure}[h]
\centering
\includegraphics[width=.7\textwidth]{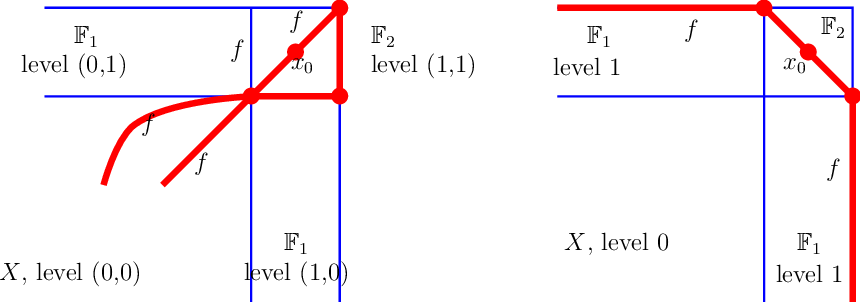}
\caption{Figure 4(a)  Example \ref{neck.2}\hskip1in  (b) Example \ref{neck.3} }
\end{figure}
\par 

\begin{ex} \label{neck.3} Finally, consider the case when $V$ is a union of the first two coordinate lines in $\P^2$ and $f_\ep :\P^1 \ra X$ are stable maps defined in homogenous coordinates by $f_\ep(z)=[\ep z, \ep z^{-1}, 1]$, all containing a marked point $x$  with coordinate $z=1$. As $\ep\ra 0$ the image of marked point $f(x)=[\ep,\ep, 1]$  falls into $p=[0,0,1]$. Rescaling the target around $p$ to prevent this gives rise to a level 1 building. The limit $f$ has now three components, all on level 1. Out of these,  only one is nontrivial and is  mapped into $\F_2$ (the one containing $x$) while the other two are trivial, each mapped into the zero section of a different $\F_1$ piece. 
\end{ex}
\par

The examples above show that we cannot avoid trivial components in the neck or in the zero divisor when $k(x)\ge 2$. The trivial components are uniquely determined by the behavior of the rest of the curve and the rescaling parameter, and are there only to make the limit continuous, such that the maps converge in Hausdorff distance to their limit.  The trivial components  satisfy only some partial version of the matching conditions, because some of their leading coefficients (but not all!) are either zero or infinity. 
\medskip

\subsection{Transversality for maps into level one buildings} 

If we assume that the refined limit $g_1:C_1\ra X_1$ has no components in $D_1$, its  resolution $\wt g: \wt C_1 \ra \wt X_1$ is an element of $\Ev^{-1}(\De_{\pm}) \subset \M_{\wt s}(\wt X_1, \wt D_1)$. 
Here $\Ev$ is the evaluation map corresponding to all the nodes of $C_1$, including its depth zero ones. Denote by $d_0$ their number, and use the notations of Remark \ref{match.naive}.     
\begin{lemma}\label{L.d=d-2} Each open stratum  $\M_{\wt s}(\wt X_1, \wt D_1)\ra \JV(X, V)$  is cut transversally and the refined evaluation map $\Ev$ is a submersion, at least at a point  $f$ with $\Aut \;C=1$ (or more generally if its graph $F:C\ra \ov\U\ti X$ is somewhere injective). 

So for generic $V$-compatible $(J, \nu)$, $\Ev ^{-1} (\De_\pm)/\cx^*$ is a smooth manifold of dimension 
\bear\label{dim.Ev.diag}
\dim \Ev ^{-1} (\De_\pm)/\cx^*= \dim \M_s(X, V) -2
\eear
at any $f:C\ra X_1$ with at least one nontrivial component in level one, at least as long as  its projection $f_0:C_0\ra X$ has $\Aut \;\ct( C) =1$.  
\end{lemma}
\begin{proof} This follows by a straightforward extension of the argument of Lemma 4.2 of \cite{ip1} to the context of normal crossing divisor. We briefly outline the main steps below. For simplicity of notation, we identify a singular space with its smooth resolution and so do not distinguish between a map $f:C\ra X_1$ and its resolution $\wt f:\wt C\ra \wt X_1$; in particular, the statements below become technically correct only after passing to appropriate resolutions. 

Assume first $f$ is an element of the absolute moduli space $\oM(\wt X_1)\ra \JV(X, V)$ without any components, nodes or marked points in $D_1$, and that $\Aut(C)=1$. We first show (a) this moduli space is cut transversally at $f$  and $\ev$ is a submersion (b) the leading order section $\si_x(f)$ is transverse to zero at $f$.   By induction on the multiplicity of the contact of $f$ at $x$ to $D_1$ this implies both that the moduli space  $\M(X_1, D_1)$ is cut transversally at $f$ as well as the transversality of the refined evaluation map $\Ev$. For both arguments, we separate the tangent and normal directions to $V$,  and use the fact that $(J, \nu)$-holomorphic map equation is linear in $\nu\in {\mathcal V}(X, V)$. Denote by $D_f$ the linearization in $f$ of this equation (keeping $C$, $J$, $\nu$ fixed). 

Consider the graph of $f$ inside $\ov\U\ti X_1\subset \ov\U\ti {\mathcal X}$ where 
$\x\ra B$ is the family of deformations (\ref{cal.X}) described in \S\ref{S.local.mod}.  Since $\Aut C=1$, the graph of $f$ is an embedded, possibly nodal curve: for each node $x_\pm$ of $C$ the graph has a node with transverse branches (because $C\subset \ov\U$ has transverse branches). As outlined in Remark \ref{R.nu.can.be} we can use the local models of both $\ov\U$ and $\mathcal X$ to construct specific variations in the parameter $\nu\in {\mathcal V}(X, V)$ supported in small balls around a somewhere injective, non-special point of the domain, which imply the required transversalities. 

Let  $g:C\ra X$ be the projection of $f$ to level zero, and regard each level one depth $k$ component $f:\Si \ra \F_k$ as a pair $(g, \xi)$ where $g:\Si\ra V^k$ is the projection of $f$ to the infinity divisor of $\F_k$ and $\xi\not\equiv 0, \infty$ a section of $g^* N_{V^k}$ with poles along the zero divisor of $\F_k$ (and zeros along the infinity divisor). The map $g$ can be regarded both as a map into $V^k$ without components in $V^{k+1}$ or as a map into the zero section of $\F^k$ without any components in the fiber divisor. So we have two linearizations at $g$, $D_g$ tangent to $V^k$ and the normal operator $D^N_g$ in the fiber direction of $\F_k$. Moreover, since $f$ is $(J, \nu)$ holomorphic then 
$g\in \M(V^k, V^{k+1})$ and $\xi$ is an element in the kernel of the $D^N_{g}$, see equation (7.1) of \cite{ip1}. 
\smallskip

This decomposes $f=(g, \xi)$. For part (a) of the transversality argument, the image of $g$ inside $\ov\U\ti X$ is embedded, so for each depth 
$k$ component $\Si$ of its domain there exists a small ball $B$ centered at $y\in\Si$ on which $g$ is somewhere injective and $g(B)\subset V^k\setminus V^{k+1}$. For any 
$\eta_y\in \Hom^{0,1}(T_y C, T_{g(y)}V^k)$  there exists a parameter $\nu\in{\mathcal V}(X, V)$ supported in a neighborhood of $(y, g(y))$ whose restriction to the graph of $g$ takes value $\eta_y$ at $y$, and thus each of the universal moduli spaces $\M(V^k)\ra \JV(X, V)$ are cut transversally at $g$.  

The normal component $\xi$ is a nonzero element of $\ker D^N_g$ (over each level one component).   Over the universal moduli space $\M(V^k)\ra \JV(X, V)$, the map 
$g\mapsto D^N_g$ defines a section of the bundle $\mathrm{Fred} \ra \M(V^k)$ which is 
transverse to each $\mathrm{Fred}^\ell$ stratum as in Lemma 6.4 of \cite{ip1},  using now the 1-jet of $\nu$ normal to $V$ around the point $(y, g(y))$. This proves that the universal 
moduli space $\M(X_1)\ra \J(X, V) $ is cut transversally at the point $f=(g, \xi)$. 
\smallskip

For part (b), for any map $f :C\ra X_1$ as above, pick any depth $k$ point $x\in  f^{-1}(D_1)$ and work separately in each direction $i\in I$. Consider the local expansion (\ref{f.exp}) of $f$ at $x$, and for each $0\le d\le s_i$ consider the section $\si_d(f, x)$ defined by the coefficient of $z^d$ in the expansion of $f^i$. By induction on $d$, one can show that $\si_d$ is transverse to the zero section  by first constructing a perturbation $\zeta$ in $f$ supported in a ball about $x$ and holomorphic in a smaller ball and then constructing a $\nu\in {\mathcal V}$ such that $D_f \zeta= -\nu$ as in Lemma 4.2 of \cite{ip1}. One could equivalently work instead with the $d$-jet of $f_i$ at $x$ and prove that is transverse to zero for any $d$.  Again, these arguments are done separately on tangent and normal components using the decomposition $f=(g, \xi)$ to make sure that such $\nu$ can be chosen compatible with $V$, proving that $\M(X_1, D_1)$ is cut transversally. At maps $f\in\M(X_1, D_1)$ whose leading coefficient 
$\si_x(f)\ne 0$ the same type of argument implies transversality of the map $\Ev_x=[\si_x(f)]$ from (\ref{enh.ev}). 
\smallskip

To prove the last statement of the Lemma, if $f$ has no trivial components, then the argument above applies to its resolution $\wt f$  to show $\Ev^{-1}(\De^{\pm})$ is smooth at $f$. If $f$ has any trivial components, the argument above applies only to its nontrivial components. On the trivial components, Gromov type perturbations vanish, so they are $J$-holomorphic. Their domain has only two marked points and are usually multiple covers. Nonetheless, because they project to a point in $V$, one can show they for {\em any} parameter $J$ they are cut transversally (with automorphism group $\Gamma_\al$), and the evaluation map at one of their marked points (but not at both) a submersion by simply calculating the cokernel of the linearization: in the normal direction, the linearization is the usual $\ov\bd$ operator on $S^2$ with values in a trivial bundle thus its cokernel vanishes. 
All together, we see that if $\Aut \ct f=1$ then $\Ev^{-1}(\De^\pm)$ is cut transversally. When $f$ has a nontrivial component in level one, the $\cx^*$ action on it is free, so the quotient is smooth. 

Finally, the calculations in Lemma \ref{dim.diff1}  together with the fact that the refined matching conditions impose an extra  $2k(x)-2$ dimensional condition for each {\em positive} depth node immediately imply (\ref{dim.Ev.diag}).
 \end{proof} 
\bigskip

\setcounter{equation}{0}
\section{The general limit of a sequence of maps}\label{s5}
\medskip

We are ready to describe what kind of maps may appear as refined limits of maps in $\M(X, V)$.  To construct the limit we inductively rescale the sequence $f_n$ to prevent (nontrivial) components from sinking into $V$, as described in the previous section. Starting with a stable map limit (level zero limit) $g_0:C_0\ra X_0$, as long as we have nontrivial components in the zero divisor we rescale again the sequence near the zero divisor to inductively obtain a  limit $g_m:C_m\ra X_m$ which is a strict refinement of the previous one  $g_{m-1}:C_{m-1} \ra X_{m-1}$ and which has no nontrivial components in the singular locus. This process terminates after finitely many steps with a final limit that has no {\em nontrivial} components in the total divisor, and satisfies a certain refined matching condition at depth $k\ge 2$ points. Note however that unlike the case of \cite{ip1}, in the limit there might be some trivial components lying in the total divisor, and this is something that cannot be avoided when $k\ge 2$, cf. Examples \ref{neck.1} and \ref{neck.2}. The matching conditions are much more involved, and are trickier to state because of the presence of these trivial curves.  

We start with the following generalization of Proposition \ref{L.lim.level.2}: 
\begin{prop}\label{L.lim.level.m}   Assume $J_n\ra J_0$ is a
sequence in $\JV(X, V)$ and consider a sequence  $f_n:C_n\ra X$   of maps in 
$\M_s(X,V)$ for a fixed refined dual graph $s$.  

Then there exists an $m\ge 0$ and a sequence of  rescaling parameters $\la_{n}=(\la_{n, l})_{l=1}^m$ such that, after passing to a subsequence, the rescaled sequence $R_{\la_n} f_n$ has a 
continuous $J_0$-holomorphic  limit $g_m:C_m \ra X_m$ into a level $m$ building with the following properties:
\begin{enumerate}[(a)] 
\item $g_m$ is obtained inductively as part of a tower 
\bear\label{f.ref.f-0} 
\xymatrix{
C_m\ar[d]_{g_m}\ar[r]^{\ct} &C_{m-1}\ar[d]_{g_{m-1}}\ar[r]^{\ct}&\dots \ar[r]^{\ct} &C_0 \ar[d]^{g_0} 
\\
X_m\ar[r]_{p}& X_{m-1} \ar[r]_{p}&\dots \ar[r]_{p}&X_0
}
\eear
of level $l$ limits $g_l:C_l\ra X_l$ starting from the stable map limit $g_0:C_0\ra X_0$; each $g_l$ is a strict refinement of $g_{l-1}$ and none of these maps has  any nontrivial components in the singular divisor;  
\item there is a tower of compatible resolutions 
\bear\label{gm.res.f} 
\xymatrix{
\wt C_m\ar[d]_{\wt g_m} \ar[r]&\wt C_{m-1}\ar[d]_{\wt g_{m-1}}\ar[r]&\dots\ar[r] &\wt C_0 \ar[d]^{g_0} 
\\
\wt X_m\ar[r]&\wt  X_{m-1} \ar[r]&\dots \ar[r]&\wt X_0
}
\eear
where the resolution $\wt g_i: \wt C_i \ra \wt X_i$ of $g_i:C_i\ra X_i$ has contact information to the total divisor of $\wt X_i$ described by an expansion of the form (\ref{exp.u.v.x}); 
\item for every positive level $l$  and any intermediate curve $C_l\ra C'\ra C_{l-1}$ all the contact points of $\wt g_l$ descend to special points of $C'$, and moreover  
\begin{itemize}
\item $s(x_-)= s(x_+)$ and $\ep(x_-)= -\ep(x_+)$  for each node  $x_-=x_+$ of  $C'$; 
\item $s (x)= s(x_n)$ and $\ep(x)=\ep(x_n)$ for each marked point $x\in C'$ which descends to the   limit of the marked points $x_n$ of $f_n$; 
\end{itemize}
(whenever both sides are defined). Furthermore, $C_l$ is obtained from $C_{l-1}$ by inserting strings of trivial components $B_x$ (broken cylinders) either between two branches $x_\pm$ of a node $x$ of $C_{l-1}$ or else at a marked point $x$ of $C_{l-1}$; 
\item[(d)] the final stage limit $g_m$, denoted $f:C\ra X_m$ and called the refined limit of the subsequence $f_n$, is a {\em relatively stable map} into $X_m$, i.e.  (i) $f$ has no nontrivial components in the total divisor of $X_m$ and (ii) $f$ has at least one nontrivial component in each positive level $l$.  
\end{enumerate}
For each refined-convergent subsequence of the original sequence, its refined limit $f$ is unique up to the $(\cx^*)^m$ action rescaling $X_m$ (described in Remark \ref{cx.m.act}). 
\end{prop}
\begin{proof} The existence of the tower of refined limits for the sequence $f_n$ follows immediately  by iterating the rescaling procedure of Proposition \ref{L.lim.level.1}. The limiting process takes place by induction on $l$ inside the deformation space  $\x_l \ra B^l$ of a level $l$ building constructed in \S\ref{S.local.mod}, with local trivializations $R_\la:X\ra X_{\la}$. At each step $l$ the convergence  $R_{\la_n(l)} f_n\ma\longra^G g_l$ is the sense of Definition \ref{D.Gr.top}. 

Denote by $J_0$ the induced parameter on $\x\ra B$ from $J_0$ by the rescaling process as in Remark \ref{R.ext.J.cal.X}. Since $J_n\ra J_0$ in $C^\infty$ in $\JV(X, V)$  then for any sequence $\la_n$ of rescaling parameters the sequence $(R_{\la_n}^{-1})^* J_n$ on $X_{\la_n}$ converges to  $J_0$ in the sense of part (iii) of Definition \ref{D.Gr.top}. 

The rescaling procedure of Proposition \ref{L.lim.level.1} starts with the sequence $(f_n, J_n)$ of elements in $\M_s(X, V)$ which Gromov converge to $(g_0, J_0)\in \ov\M(X)$ as maps into a level zero building $X$. When $g_0$ has nontrivial components in the zero divisor, it produced a sequence $\la_n(1)\in  \cx^*$ of rescaling parameters such that  
$(R_{\la_n(1)}f_n,  R_{\la_n} J_n) \ma\longra^G (g_1, J_0)$ as maps into $\x_1$. The construction of the refinement $g_1$ was semilocal, in the neighborhood of the zero divisor, so as long as the parameters converge, it extends to any sequence $f_n$ of maps into the fibers of $\x(l)$ to construct a refined limit inside $\x(l+1)$. This is because we can use another set of local trivializations to smoothly identify the tubular neighborhoods of zero sections of the smooth fibers of $\x_l$ with the standard 
local model of the neighborhood of $V$ in $X$. 

The rescaling process terminates in finitely many steps; there is an a priori 
uniform bound  $M$ on the number of times we need to rescale,  depending only the topological information $s$. More precisely, we first prove by induction on the level $m$ the existence of such a tower of refined limits which satisfies the properties (a)-(c).  Because both the topological type of the domain as well as the homology class of the image of $f_n$ is fixed (being part of the information encoded by the dual graph $s$) then by Gromov compactness,  after passing to a subsequence, we get a limit $g_0:C_0\ra X_0$ in the usual stable map compactification.  The number of components of $g_0$ is uniformly bounded by a constant $K$ which depends only on the dual graph  $s$, and in fact only on $(A_s, \chi_s, n_s)$,  but not on the sequence $f_n$. 
%

If the level zero limit $g_0$ has no components in $V$, then it already satisfies all conditions (a)-(d). So assume $g_0$ has at least one component in $V$; note that this component must be nontrivial as there are no trivial components in level 0. 

Next Proposition \ref{L.lim.level.2} implies the existence of a sequence of rescalings with a refined limit $g_1:C_1\ra X_1$ which satisfies the properties (a)-(c). The key idea was to use the 
rescaling technique normal to $V$ of Proposition \ref{L.lim.level.1} by finding a sequence of annular regions (\ref{defn.an.r}) around $V$ that carried essentially $\al_V/2$ energy see (\ref{def.la.n}); this region in effect acted as a buffer zone between level 0 and  level 1, see Figure \ref{F.3}. 
  

If $g_1$ has no components in the zero divisor, the process terminates; otherwise we use the rescaling technique of Proposition \ref{L.lim.level.1} again and keep inductively rescaling the sequence around the zero divisor by choosing at each step $l$ yet another sequence of same kind  of annular regions (with radii measured with respect to the rescaled metric of the previous step) on which the energy is essentially  $\al_V/2$; this is precisely what is needed to guarantee that the new limit is a strict refinement of the one obtained in the previous step, and also similarly satisfies all the other  properties (a)-(c). 
In the presence of depth $k\ge 2$ strata, the annular regions  (\ref{defn.an.r}) intersect nontrivially, see Figure \ref{F.3}.



Each level of this rescaling procedure may produce more and more trivial components, but $g_l$ and $g_{l-1}$ have the same nontrivial domain components (their domains differ only by trivial components by (c)), so this part of the domain stays constant in $l$. In particular each such nontrivial component $\Si$ is already part of the domain $C_0$ of the original limit 
$g_0:C_0\ra X_0$, where it has a depth $0\le k(\Si)\le \dim X$. Moreover, because $g_l$ is a strict refinement of $g_{l-1}$ then the depth of each nontrivial component $\Si$ in $g_l$ is less or equal to its depth in $g_{l-1}$ with at least one component with strict inequality. 


This means that the rescaling process must terminate in at most $M=K(\dim X+1)$ steps, with a limit $g_m$ which satisfies (a)-(c) but has no more components in the zero divisor. It remains to check that $g_m$ has a nontrivial component in each positive level $l\le m$. 
Assume by contradiction that $g_m$ has no nontrivial component in level $l$ whatsoever (it already has none in the total divisor); but $g_m$ agrees with $g_l$ up to level $l$, and therefore $g_l$ cannot have any nontrivial component in level $l$ either (away from the zero section). But this contradicts the fact that $g_l$ was a strict refinement of $g_{l-1}$, and completes the proof of property (d). \end{proof}   
\begin{rem} 
The number of levels may be bigger than the number of nontrivial components, as Figure 3(b) of Example \ref{ex.4dim.b} illustrates. In that case the limit is a map into a level 2 building, with only one nontrivial component which is simultaneously in both  level one and two (in different directions). It has the crucial property of not being fixed by the $\cx^*$ action in neither level. 
\end{rem}
Assume now  $f_n:C_n\ra X$ is a sequence of maps in $\M_s(X,V)$, and  let  $f:C\ra X_m$  be the relatively stable limit constructed by Proposition \ref{L.lim.level.m}; it fits as part of the diagram: 
\bear\label{g.ref.fmm} 
\xymatrix{
\wt C\ar[r]^{\xi}\ar[d]_{\wt f}&C\ar[r]^{\ct}\ar[d]_{f} & C_0\ar[d]^{f_0}
\\
\wt X_m\ar[r]^{\xi}&X_m\ar[r]^p& X
}
\eear
where $\wt f\in \ev^{-1}(\De)\subset \M(\wt X_m, \wt D_m)$ is its resolution, while $f_0\in \oM(X)$ is the usual stable map limit. We next describe the semi-local behavior of this limiting process in a neighborhood of the fiber $B_x$ of $C\ra C_0$ over a point $x\in C_0$ to extract further properties of the refined limit. By part (c), $B_x$ is either a point $x$ or else it is a string of trivial components with two end points $x_\pm$ (broken cylinder).  In the later case, which can happen only when $x$ is a special point of $C_0$,  we order the components $\Si_r$ of $B_x$ in increasing order as we move from one end $x_-$ to the other $x_+$, and make the following definition
\bear\label{D.stretch}
\mbox {the {\bf stretch} of a point $x\in C_0$  is $r(x)=$ the number of components of $B_x=\st^{-1}(x)$}
\eear
where by convention $r(x)=0$ whenever $B_x=x$. 
\smallskip

\subsection{Properties of the relatively stable  limit $f$ around $B_x$} Consider any node $x$ of $C_0$, and fix a choice $x^\pm $ of its two branches.  Order the components $\Si_r$ of $B_x$ and write
 \bear\label{B=ordered}
 B_x=(\Si_1, x_1^\pm)\cup(\Si_2, x_2^\pm)\cup \dots \cup(\Si_r, x^\pm_r)
 \eear
with nodes $x^+_r=x^-_{r+1}$ for each $r=1, \dots r(x)-1$. Set $x^+_0=x^+$ and 
$x^-_{r+1}=x^-$ to include the case when $B_x=x$. 

Consider next  $p_\pm= \wt f(x^\pm)\in \wt X_m$ the  two images in the resolution, while $p=f_0(x)$ is the common image in $X$.   By construction $f(B_x)$  lies in the fiber $F_p$ over  $p$ of the collapsing map $X_m\ra X$, where we can separately work one normal direction to $V$ at a time. Fix any of the directions $i\in I $ indexing the branches of $V$ at $p$, and let 
$\pi_i$ be the projection onto that direction,  defined on a neighborhood $U_p$ of $F_p$ in the semi-local model described in Remark \ref{rem.coord-iterate.m}.  The target of $\pi_i$ is nothing but the (global) model of the deformation of a disk $D^2$ in $N_i$ which is being rescaled $m$ times at 0; it is described in terms of the coordinates  $u_{i, l}$ and $v_{i, l}= u_{i, l}^{-1}$  by 
\bear\label{uv=i}
u_{i,l-1} v_{i,l}=\la_{i, l} \quad \mbox {for all } l=1, \dots,  m
\eear
for any collection $\la= (\la_{i, l})$ of small rescaling parameters. Choose also local coordinates $z, w$ at $x_\pm$ which then induce local coordinates on the universal curve of the domains at $x$ (the one containing $C_0$, which was assumed to be stable); the nearby curves are then described in the ball $B(x, \ep)$  by 
\bear\label{zw=mu}
zw= \mu(x)
\eear
where intrinsically the gluing parameters $\mu$ are local coordinates at $C_0$ on the moduli space of stable curves. 
Similarly choose (global) coordinates $z_j,  \; w_j=z_j^{-1}$ on the $j$'th component $\P^1$ of $B_x$, where $1\le j \le r(x)$, and where we set $z_0= z$ and $w_{r(x)+1}=w$. These provide global coordinates in the neighborhood $O_x$ of $B_x$ obtained as the inverse image of $B(x, \ep)$ under the collapsing map $C\ra C_0$, in which the nearby curves are described by
\bear\label{z.w.mu.j}
z_{r-1} w_{r} =\mu(x_r) \quad \mbox{for all } r=1, \dots,  r(x)+1
\eear  
where $\mu(x_r)$ is the gluing parameter at the $r$'th node $x_r$ of of $B_x'$, the intersection of $C$ with $O_x$.  In particular 
\bear\label{mu=mu.j}
\mu(x)= \prod_r \mu(x_r)
\eear
Note that $B'_x$ has $r(x)+2$ components $\Si_r$, the first and the last are the disks about $x^\pm$ while the remaining ones are the spherical components of $B_x$, see  (\ref{B=ordered}). 

For each component $\Si$  of $B_x$, $f$ may only have a partial contact information along the singular divisor at the two points  $0_\Si$ and $\infty_\Si$. In fact, in the coordinates on both the domain and target described above, $f|_{\Si}$ has an associated coefficient $a_i(0_\Si)=a_i^{-1}(\infty_\Si)\ne 0$ and a contact multiplicity $s_i(\Si)\ge 0$ for all $i\notin  I^\infty_\Si$, see Definition \ref{def.triv.comp} (with $ I^0_\Si$ indexing the directions in which $s_i(\Si)=0$, i.e. $\wt f(\Si)$ stays away from that stratum). Furthermore, because we already know that we have matching contact multiplicities at each node in all directions in which both sides are defined, then $s_i(\Si)=s_i(x)$ for all $i\notin I^\infty_\Si$. In the remaining directions $f$ still has a coefficient $a_i(\Si)$ which is 0 or $\infty$; the contact multiplicity $s_i(\Si)$ is technically undefined, but we can {\em define} it to be $s_i(x)$ for all $i\in I(x)$. 

As described in Remark \ref{R.intro.I}, each node $y$ of any intermediate curve $C\ra C'\ra C_0$ which projects to $x$ comes with a partition of the original indexing set $I(x)$ of the stratum containing $f_0(x)$ into $I^\infty(y)=I^\infty(y_-)\cup I^\infty(y_+)$,  $I^0(y)$ and then $I^\pm(y)$ induced by the contact information of the limit $f:C\ra X_m$ to the total divisor.  The directions $i\in I^\pm(y)$ come with a multiplicity $s_i(y)=s_i(x)>0$ and opposite signs 
\bear\label{ep.x=ep.y}
\ep_i(x_+)=\ep_i(y_+)=-\ep_i(y_-)=-\ep_i(x_-)\ne 0
\eear recording the two opposite sides of the level one building from which the two branches of $f$ come into the singular divisor, and also two leading coefficients $a_i(y_\pm)\ne 0$, intrinsically elements
\bear\label{a.v.y}
a_i(y_+) \in N_i^{\ep_i(x_+)} \otimes (T_{y_+} C )^{s_i(x)}\quad  \mbox{ and } \quad a_i(y_-) \in N_i^{\ep_i(x_-)} \otimes( T_{x_-} C)^{s_i(x)} 
\eear
where $N_i$ is the branch of the normal bundle to $V$  indexed by $i\in I^\pm(y)$. As we have seen in Example \ref{neck.3}, in general for a node $x$ of $C_0$, $l_i(x_+)$ could be bigger than $l_i(x_-)$ in some directions and smaller in others, so we  denote by  
\bear\label{def.min.max.l} 
l_i^-(y)=\min_\pm  \{ l_i(y_\pm) \}  \quad \mbox{ and } \quad  l_i^+(y)=\max_\pm \{ l_i(y_\pm)\}
\eear
for each node $y$ of any intermediate curve $C\ra C'\ra C_0$. 

Consider next  the restriction $f^i$ of $\pi_i\circ f$ to $B_x'$;  after collapsing all constant components (keeping only those for which $i\in I^\pm_\Si$) it has a stable map model $f^i:B'_i\ra (D^2)_m$ defined on a slightly bigger curve $B'_i$ containing $B_i$. We have a similar description of the nearby curves in terms of $B_i'$, but where now at each node $y$ of $B_i'$ the gluing parameter  (\ref{z.w.mu.j}) is  
\bear\label{mu.at.y}
\mu(y)=\prod _{z}\mu(z)
\eear 
where the product is over all nodes $z$ of $C$ in the fiber of the collapsing map $B_x\ra B_i$ at $y$. This formula extends (\ref{mu=mu.j}) which would correspond to the collapsing map $B_x\ra x$. 

\begin{lemma}\label{triv.comp.descr} Using the notations above and of (\ref{g.ref.fmm}), consider a sequence $f_n:C_n\ra X$ which after rescaling by the sequence $\la_n= (\la_{n, l})_{l=1}^ m$ has a  relatively stable limit $f:C\ra X_m$ as constructed by Proposition \ref{L.lim.level.m}.  

Fix a special point $x$ of $C_0$ and denote by $\mu_n(x_r)$ the corresponding parameters (\ref{z.w.mu.j}) describing the domains $C_n$ in a neighborhood of the fiber $B_x$ of $C\ra C_0$ over $x$. Finally, fix a direction $i\in I(x)$ of $V$ at $f_0(x)$. 

When $B_i\ne x$ then the restriction of $f^i$ to $B_i$ is a degree $s_i(x)$ chain of trivial components in $(D^2)_m$ stretching from $\pi_i(p_-)$ to $\pi_i(p_+)$, both of which must be on the total divisor of $(D^2)_m$. In particular for each node $y$ of $B_i'$, $f^i(y_-)=f^i(y_+)$ lands in the total divisor, with the two branches of $f^i$ landing on opposite sides of the divisor.  Furthermore, $f$ satisfies the following {\em refined matching condition at $y$}: 
\bear\label{enh.math.la.i.l}
\lim_{n\ra \infty}\frac{ \la_{n, l_i^+(y)}}{ \mu_n(y) ^{s_{i}(x)} }=a_i(y_+) a_i(y_-)
\eear
where $a_i(y_\pm)$  and  $s_i(y_\pm)=s_i(x)$ are the two leading coefficients of $f$, and respectively the contact multiplicity, while $l_i^+(y)$ is defined by (\ref{def.min.max.l}) and $\mu(y)$ by (\ref{mu.at.y}).

When $B_i=x$  is a node of $C_0$ then  $f^i(x_-)=f^i(x_+)$; if this lands in the total divisor of $(D^2)_m$, then  $f$ satisfies the corresponding refined matching condition at $x$ in the direction $i$:
\bear\label{enh.math.la.i.l.0}
\lim_{n\ra \infty}\frac{ \la_{n, l_i(x)}}{ \mu_n(x) ^{s_{i}(x)} }=a_i(x_-) a_i(x_+) 
\eear
\end{lemma} 
\begin{proof} This follows by refining the arguments in the proof of Lemma \ref{L.lim.level.2}, using also the information described above. Here we work locally in the neighborhoods $O_x$ and $U_x$ of $x$ and $p$ described above, where we can separately project onto the $i$'th direction. Denote by $C_n'$ the intersection of $C_n$ with $O_x$. 

Because we already know that $R_{\la_n}f_n$ converge to $f$, the projections $h^i_n$ of their restrictions to $C_n'$ will also converge, and the  limit will be precisely $f^i:B_i'\ra (D^2)_m$. In fact, $h^i_n$ are maps from  $C_n'$ into  nothing but 
$D^2_{\la_{n}}$, the $m$-times rescaled disk using the rescaling parameters $\la_{n} =(\la_{n,1}, \dots, \la_{n, m})$, and that is precisely the situation in which Lemma 5.3 of \cite{ip2} applies to give refined matching conditions above at each mode $y$ of $B_i'$. \end{proof}
\smallskip

First of all, Lemma \ref{triv.comp.descr} provides a concrete description of the trivial components in the refined limit $f:C\ra X_m$. The limit $f$ is continuous and its domain $C$   is obtained from $C_0$ by inserting strings of trivial spheres $B_x$ to stretch the image curve across the levels $l_i(x^\pm)$ in a zig-zagging fashion either (a) between $x_-$ and $x_+$ if $x$ is a node of $C_0$  or (b) at a contact point $x_-$ of $C_0$ with its respective zero section to stretch it all the way to a contact point $x_+\in C$ that is mapped to the zero section $V_m$.  The neck region of $C_n$ at $x$ is roughly equal to a trivial cylinder mapped in the fiber of the neck of the target over $f_0(p)$, which then gets further stretched, possibly several times to accommodate the rescaling done to catch all the nontrivial components of $f$. Once we fix an order $x_-$, $x_+$ of the two branches of $C_0$ at $x$, this orders the components of the chain $B_x$ as well as the two branches of each of its nodes. 

The first consequence of Lemma \ref{triv.comp.descr} is the following: 
\begin{cor} \label{C.triv.forget}  Consider a sequence $f_n:C_n\ra X$ which has a  relatively stable limit $f:C\ra X_m$ as constructed by Proposition \ref{L.lim.level.m}. For any node $y$ of $C$ that projects to $x$ in $C_0$ 
\bear\label{ep=l-l}
\ep_i(y_-)=-\ep_i(y_+)=l_i(y_+)-l_i(y_-)=\ep_i(x_-)= \mathrm{sign} (l_i(x_+)-l_i(x_-))&
\eear
in any direction $i\in I(x)\setminus I^{\infty}(y)$. 

Moreover, for any special point  $x$ of $C_0$, the image of the chain of trivial components $B_x$  connects $f(x_-)$ and $f(x_+)$ in the fiber $F_p$ of $X_m$ over $p=f_0(x)$ such that: (i)  at each step the levels change by either zero or $\ep_i(x_-)$ in  each direction $i\in I(x)$, and (ii) at each step the level changes in at least one direction. 
\end{cor} 
For each intermediate node $y$ that projects to $x$, we can therefore {\em define} a multi sign map by 
\bear
\ep_i(y_-)=\ep_i(x_-)= \mathrm{sign} (l_i(x_+)-l_i(x_-)) 
\eear
in all directions $i\in I(x)$. This definition agrees with the previously defined multi sign map in all directions $i\notin I^\infty(y)$, and is the only definition that makes the naive matching 
conditions that appear in part (c) of Proposition \ref{L.lim.level.2} true in all directions $i\in I$.  It also restricts the collection of possible refined dual graphs. 

Next, the asymptotics in Lemma \ref{triv.comp.descr}  impose further restrictions on both the rescaling parameters as well as on the refined limit $f$ at depth $\ge 2$ contact points with the total divisor. 
\begin{ex} \label{E.B.x=x}Take the simplest example when $B_x=x$, i.e. $x$ is a node between two nontrivial components and assume $x$ is a depth $|I^\pm(x)|\ge 2$ contact point to the singular divisor $W$. Lemma \ref{triv.comp.descr} implies that the limit $f$ is continuous at $x$; because the leading coefficients $a_i(x_\pm)$ of the limit are bounded away from zero and infinity then  (\ref{enh.math.la.i.l.0})  implies that the rate at which 
\bear
{\la_{n, l_i(x)}}^{1/{s_i(x)}} \ra 0
\eear
is the same in all directions  $i\in I^\pm(x)$, generalizing Remark \ref{R.cx}.  After multiplying the rescaling parameters by a constant factor in each level $l_i(x)$, we can arrange that the ratio
\bear\label{la.i=la.j}
{\la_{n, l_i(x)}}^{1/{s_i(x)}} / {\la_{n, l_j(x)}}^{1/{s_j(x)}}  \ra 1
\eear
for all $i, j\in I^\pm(x)$. Then   (\ref{enh.math.la.i.l.0})  becomes 
\best
[a_i(x^-)]_{i\in I^\pm(x)}= [a_i(x^+)^{-1}]_{i\in I^\pm(x)}\in \P_{s(x)}(N_{s(x)} W_{J(x^-)}) = 
\P_{s(x)}(N^*_{s(x)} W_{J(x^+)})
\eest
or equivalently
\bear
\Ev_{x_\pm}(\wt f)\in \De^{\pm}
\eear
generalizing (\ref{enh.ev.match}), and providing the statement of the refined matching conditions that the limit must satisfy at a node $x$ whose stretch $r(x)=0$.
\end{ex} 
If we denote the relative rescaling parameters between two levels $l_1$ and $l_2$ by 
\bear\label{La=min.max}
 \Lambda_{n}(l_1, l_2)= \prod_{k=\min (l_i)+1}^{\max(l_i)} \la_{n, k}
 \eear
 so that $\la_{n, l}=\Lambda_n( l-1, l)$, then the asymptotics in Lemma \ref{triv.comp.descr} extend across several levels to give the following
 \begin{cor}\label{cor.enh.lim.cond} Consider the situation of  Lemma \ref{triv.comp.descr}.  Then for any node $y$ of any intermediate curve $C\ra C'\ra C_0$: 
\bear\label{La=mu.asympt}
\lim_{n\ra \infty} \frac {\Lambda_n(l_i(y_\pm))}{\mu_n (y)^{s_i(y)}}=a_i(y_-) a_i(y_+)
\eear
for each $i\in I^\pm(y)$. 
\end{cor}

\begin{rem}\label{local.model} This implies that the local model of $\ov\M_s(X, V)$ near a refined limit $f:C\ra X_m$ into a level $m$ building can be described as a subset of  tuples 
$(\wt f,\la, \mu)$  near $(\mu, \la)=0$, up to the $(\cx^*)^m$ action on $(\wt f, \la)$. Here  $\wt f\in \M_{\wt s}(\wt X_m, \wt D_m)$, $\mu$ has one coordinate $\mu(z)\in \L_{z_-}\otimes \L_{z_-}$ for each node $z$ of $C$, while $\la$ has one coordinate $\la_l$ for each positive level $l$ of $X_m$.  The asymptotics (\ref{La=mu.asympt}) imply that the local model near $(\mu, \la)=0$ consists of such  
$(\wt f,\mu, \la)$ which satisfy the following full set of refined matching conditions:
\bear\label{a.mu=la.before.log}
a_i(y_-) a_i(y_+) \mu(y)^{s_i(x)} =\La(l_i(y_\pm))
\eear
for each $i\in I^\pm(y)$ and for each node $y$ of any intermediate curve $C\ra C'\ra C_0$. 

As we have already seen in Remark \ref{R.match.enh.log} these equations are best analyzed by taking $\log$ of both sides, obtaining the following linear system of equations 
\bear\label{full.enh.log}
\log a_i(y_-)+\log  a_i(y_+) + s_i(x)\ma\sum_{z\in D(y)}  \log \mu(z) =\ma\sum_{l=l_i^-(y)+1}^{l_i^+(y)}
\log \la_l
\eear
for each $i\in I^\pm(y)$ and for each node $y$ of any intermediate curve $C\ra C'\ra C_0$, where $D(y)$ denotes the collection of nodes $z$ of $C$ which collapse to $y$ under $C\ra C'$. This is typically an overdetermined system of equations, and the fact  that it  has solutions $(\mu,\la)\in (\cx^*)^{|D|}\ti (\cx^*)^m$ arbitrarily close to $(0, 0)$ is not  automatic, imposing restrictions both the rescaling parameters $\la$, the leading coefficients of $\wt f$ as well as on the topological type of the limit generalizing those in the Example \ref{E.B.x=x}. 

In particular, the asymptotics as $n\ra 0$ of the real part of parameters $\log\mu_n$ and 
$ \log \la_n$ imply that the following system in variables $\al(z)$ (one for each positive depth node $z$ of $C$)  and $\al_l$ (one for each positive level $l$ of $X_m$) must have strictly negative solutions:
\bear\label{full.enh.log.arg}
\ma\sum_{z\in D(y)}  \al (z) =\frac 1{s_i(x)}\ma\sum_{l=l_i^-(y)+1}^{l_i^+(y)}\al_l
\eear
with one equation for each direction $i\in I^\pm(y)$ and for each node $y$ of any intermediate curve $C\ra C'\ra C_0$.  One can eliminate the variables $\al(z)$ to get a set of linear conditions involving only $\al_l$ which then restricts the asymptotics of $|\la_{n, l}|$ as $n\ra \infty$.





 \end{rem}

\begin{rem}\label{R.triv.forget} 
Corollary  \ref{C.triv.forget} implies that for each direction $i\in I(x)$,  and each level $l_i^-(x) \le l\le l_i^+(x)$, there is precisely one node  $y_{i, l}(x)$ of $B'_i$ on level $l$ in direction $i$, which lifts to two points  $y^\pm_{i,l}(x)\in \wt C$  at which  $f$ has a well defined contact information in direction $i$, together with the string of trivial components $B_{i, l}(x)\subset B_x$ of $C$ on which $f$ is constant in direction $i$ and thus which are precisely all the level $l$ components of $B_x$ in direction $i$. 
With this notation, (\ref{enh.math.la.i.l}) becomes
\bear\label{enh.math.la.i.l.2}
\lim_{n\ra \infty}\frac{ \la_{n, l}}{ \mu_n(y_{i,l}(x)) ^{s_{i}(x)} }=a_i(y^+_{i,l}(x)) a_i(y^-_{i,l}(x))
\eear
for each  $i\in I(x)$,  and each level $l^-(x)<l\le l^+(x)$.  Because the right hand side of (\ref{enh.math.la.i.l.2}) is finite and nonzero,  the refined matching conditions {\em a fortiori} give  conditions on the relative rates of convergence of the rescaling parameters (involving the contact multiplicities), extending those of Remark \ref{R.cx}. 
\end{rem}
\begin{ex}\label{E.triv.forget} To see how this works in practice, assume for example that  $x\in C_0$ is a node such that $B_{i_1, l_1}(x)= B_{i_2, l_2}(x)$.  If the levels $l_1=l_2$ then the multiplicities $s_{i_1}(x)= s_{i_2}(x)$ must be the same, while if $s_{i_1}(x)\ne s_{i_2}(x)$ then  $l_1\ne l_2$ and the relative rate of convergence to zero of the two rescaling parameters in these two levels must be related, more precisely the two rates of convergence of   $\la_{n, l_j}^{1/s_j}$ as $n\ra 0$ are equal (as their limit is a  bounded, nonzero constant involving the leading coefficients of $f$).  This was precisely the case in Example \ref{neck.1}  (b). 
\end{ex} 
\smallskip


Denote by ${\mathcal N}(x)$ the collection of nodes of any intermediate curve $C\ra C'\ra C_0$ that project to $x\in C_0$ and by $L$ the collection of unordered pairs of distinct levels $(l_1, l_2)$ of the building $X_m$. 
Consider the following equivalence relation on  
\bear\label{def.T.part}
{\mathcal T}=\{\; y \; |\; y\in {\mathcal N}(x)\mbox{ for some } x\in C_0\} \sqcup L
\eear
generated by the relation $y\sim (l_1, l_2)$ if there exists a direction $i\in I^\pm(y)$ such that 
$\{l_i(y_\pm)\}=\{ l_1, l_2\}$. This partitions the set ${\mathcal T}$ into equivalence classes 
$\{[p]\;|\;p\in{\mathcal P}\} $. Note that this depends only on the refined dual graph $s$ of $f$. 

\begin{cor}\label{cor.enh.cond.2} Consider the situation of  Lemma \ref{triv.comp.descr}.  Then after multiplying the rescaling parameters by a constant factor in each level and  passing to a further subsequence, there exits (i) a sequence of $\cx^*$-parameters $t_{n, [p]}\ra 0$ indexed by equivalence classes $[p]\in \mathcal P$ and (ii) for each $p\in \mathcal T$ there exists a positive weight $w(p)>0$ and a coefficient $c(p)\in \cx^*$ such that  
\bear\label{asympt}
 \lim_{n\ra \infty} \frac {\La_{n}[ l^\pm]}{ t_{n,[l^\pm ]}^{w(l^\pm)}}= 1  \quad\mbox{ and }\quad
\lim_{n\ra \infty} \frac {\mu_n (y)}{ t_{n, [y]}^{w(y)}}= c(y). 
\eear
for each pair of distinct levels $l^\pm$ and each node $y\in {\mathcal N}(x)$. In particular, 
\bear\label{a.a.c=1}
a_i(y_-) a_i(y_+) c(y)^{s_i(x)} =1 
\eear
for all nodes $y\in  {\mathcal N}(x)$ and all directions $i\in I^\pm (y)$.
\end{cor} 
\begin{proof}  The conclusion of  Lemma \ref{triv.comp.descr} imposes conditions on the relative rates of convergence to zero of both the rescaling parameters $\la_{n, l}$ in the target, and also those of the domain $\mu_n(y)$. For each parameter $p\in {\mathcal T}$, let  $t_{n,  p} = \La_{n}(l_\pm)$ 
if $p=\{l_\pm\}$ or respectively $t_{n, p}=\mu_n(y)$ if $p=y$, where $\mu(y)$ is as in (\ref{mu.at.y}) and $\La_n(l_\pm)$ as in (\ref{La=min.max}). Equation  (\ref{La=mu.asympt}) implies that if 
$p_1\sim p_2$ then there exists some $w>0$ such that $t_{n, p_1}^w/ t_{n, p_2}$ is uniformly bounded away from zero and infinity for $n$ large. This partitions ${\mathcal T}$ and therefore $\mathcal P$ into equivalence classes each one corresponding to an independent (over  $\R_+$) direction of convergence to zero of these parameters. As ${\mathcal T}$ is finite, after passing to a further subsequence of $f_n$, we can arrange that all these quotients  have a finite, nonzero limit (within each equivalence class). In particular choosing a representative $ t_{n,[p]}$ for each equivalence class $[p]\in {\mathcal P}$ gives the asymptotics 
\bear\label{asympt.12}
 \lim_{n\ra \infty} \frac {\La_{n}[ l^\pm]}{ t_{n,[l^\pm ]}^{w(l^\pm)}}= c(l^\pm)  \quad\mbox{ and }\quad
\lim_{n\ra \infty} \frac {\mu_n (y)}{ t_{n, [y]}^{w(y)}}= c(y). 
\eear
where $w(p)>0$ and $c(p)\in \cx^*$ for each $p\in {\mathcal T}$. Multiplying the sequence $\la_{n, l}$ by the constant factor $c(l-1, l)^{-1}$ in each positive level $l$ gives us another equivalent sequence of rescaling parameters and another representative $f$ of the refined limit for which now $c(l^{\pm})=1$ for all distinct pairs of levels $l^{\pm}$. With these asymptotics, equation (\ref{La=mu.asympt}) immediately implies (\ref{a.a.c=1}). \end{proof}

Recall that a priori each trivial component $\Si$ that landed in the total divisor did not have a well defined contact information in the directions $i\in I^\infty(\Si)$. However, as we have already started seeing above, knowing the contact information of all the nontrivial components allows us to formally extend the contact information of the trivial components even in the directions $I^\infty$ in which their coefficients are zero or infinity, and thus the geometric contact information is technically undefined.  


\begin{defn}\label{def.triv.comp.decorated}  A {\em decorated trivial component} is a trivial component 
$f_\Si : (\Si, x^+, x^-) \ra (\F_V, V_0\cup V_\infty)$ that does not lie in the infinity section, together with the following extra data: 
\begin{enumerate}[(a)] 
\item a fixed isomorphism $T_{x^+} \Si\ma\cong_{\phi_\Si} T_{x^-}^*\Si$; 
\item the marked points $x_\pm$ have an attached multiplicity $s_i(x^+)=s_i(x^-)\ge 0$, two opposite signs $\ep_i(x_-)=-\ep_i(x_+)\ne 0$ and two dual elements $a_i(x_\pm) \in N_i^{\ep_i(x_\pm) }\otimes (T_{x_\pm}^*\Si)^{s_i(x)}$ for each direction $i\in I$, which agree with the usual contact information to the total divisor $D_m$ in the directions in which that can be geometrically defined. 
\end{enumerate} 
\end{defn} 
Denote by 
\bear\label{def.dec}
\M^{triv} (\wt X_m, \wt D_m)
\eear
the space of decorated trivial components (up to reparametrizations) and include them as part of 
$\M(\wt X_m, \wt D_m)$, as they now 
come with a well defined evaluation map (\ref{or.ev}), leading coefficients section (\ref{sect.lead}) and a refined evaluation map (\ref{enh.ev}) defined using the extra decorations. 

With this, the conclusion of Proposition \ref{L.lim.level.m} can strengthened as follows:
\begin{theorem}\label{thm.exist.lim}  Consider a sequence $f_n$  of $J_n$-holomorphic maps in $\M_s(X,V)$,  with $J_n\ra J$ in $\JV(X, V)$. 
Then there exists an $m\ge 0$ and a sequence of rescaling parameters $\la_k\ra 0$ in $(\cx^*)^m$ such that, after passing to a subsequence,  $R_{\la_n} f_{n}$ has a $J$-holomorphic limit $f:C\ra X_m$ with the following properties: 
\begin{enumerate}[(a)]
\item $f$ has a resolution $\wt f:\wt C\ra \wt X_m$ which is an element of 
$\M_{\wt s}(\wt X_m, \wt D_m)$ for some refinement of $s$ and a projection $f_0:C_0\ra X$ which is an element of $\ov\M(X)$; these fit in the diagram: 
\bear\label{g.ref.fmm.0} 
\xymatrix{
\wt C\ar[r]^{\xi}\ar[d]_{\wt f}&C\ar[r]^{\ct}\ar[d]_{f} & C_0\ar[d]^{f_0}
\\
\wt X_m\ar[r]_{\xi}&X_m\ar[r]_p& X}
\eear
\item the inverse image $f^{-1} (W_m)$ of the singular divisor consists only of nodes of $C$ or trivial components, while the inverse image of $f^{-1} (V_m)$ of the zero divisor consists only of marked points of $C$ or trivial components; moreover, for each special point $x\in C_0$, its inverse image $B_x$ in $C$ is a (possibly trivial) chain of decorated trivial components (broken cylinders) either between the  two branches $x_\pm$ of a node $x$ of $C_0$ or else at a marked point $x_0$ of $C_0$; 
\item  $f$  satisfies the {\em naive matching conditions}: for each node $y$ of $C$
\bear\label{naive}
f(y_-)=f(y_+), \quad s(y_-) = s(y_+), \quad \ep(y_-)=- \ep(y_+) 
\eear
 or equivalently $\ev_{y_\pm}(\wt f)\in \De$. 
\item $f$  satisfies the {\em refined matching conditions}: for each positive depth node $y$ of $C$, there exists  $c(y)\in \cx^*$ such that 
\bear\label{enh.a.y=d.full}
a_i(y_-) a_i(y_+) c(y)^{s_i(y)} = 1
\eear
for all $i\in I$, or equivalently $\Ev_{y_\pm}(\wt f)\in \De^\pm$. 
\item the sequences of parameters $\mu_n$ and $\la_n$ satisfy the 
asymptotics (\ref{asympt}) as $n\ra \infty$;  in particular, the linear system of equation (\ref{full.enh.log.arg}) has strictly negative solutions. 
\item $f$ is relatively stable, that is $f$ has at least one nontrivial component on each {\em positive} level $l$. 
\end{enumerate}
\non For each refined-convergent subsequence, its limit $f$ satisfying all these conditions is unique up to the action of a complex torus $T\le (\cx^*)^m$. \end{theorem}
\begin{proof} We first use Proposition \ref{L.lim.level.m} to obtain some limit $f:C\ra X_m$, defined up to the 
$(\cx^*)^m$ action on $X_m$, and which has all the properties described there. Fix such a representative $f:C\ra X_m$ of the limit; it will come with a resolution $\wt f$ which may have some trivial components, but no nontrivial components in the total divisor. The projection $f_0$ is the stable map limit of the original (unrescaled) sequence, and these fit in the diagram (\ref{g.ref.fmm.0}).   

However, so far the trivial components of the resolution $\wt f$ are yet undecorated, so in particular may have only a partial contact information to the total divisor. We need to show that there is an essentially unique way to decorate them in such a way compatible with the rest of the conditions. 

For each point $x$ of $C_0$ choose coordinates around $B_x$ in the domain and respectively around the fiber $F_p$ of $X_m \ra X$ over $p=f_0(x)$ as described before Lemma \ref{triv.comp.descr}. Note that when $B_x\ne x$ this involves a choice of dual coordinates $w_j=z_j^{-1}$ at the two end points of the trivial component $\Si_j$ which intrinsically corresponds to a choice of an isomorphism between the tangent space to $\Si_j$ at one of the points and its dual at the other point. Recall that each special fiber $B_x$ of $C\ra C_0$ was a string of trivial components with two end points 
$x_\pm$ (broken cylinder), which occurred only when $x$ was a special point of $C_0$.  We make the convention that 
if $x$ is a marked point of $C_0$ then the end $x_+$ of $B_x$ corresponds to the marked point $x\in C$ while the other end $x_-$ is where $B_x$ gets attached to the rest of the components of $C$.

Next, Lemma \ref{triv.comp.descr} implies that the limit $f$ satisfies the naive matching conditions (c) at all the special points  $y$ of any of the intermediate curves $C\ra C'\ra C_0$ in all the directions $i\notin I^\infty(y)$. For each trivial component $(\Si, y^-, y^+)$ of $C$ that is part of the trivial string $B_x$ with two end points $x_\pm$, and for all the directions 
$i\in I^\infty$ we  formally set $s_i(y)= s_i(x)$ and $\ep_i(y_\pm)=\ep_i(x_\pm)= -\ep_i(x_\mp)$ respectively. This choice is the unique one for which the naive matching conditions (c) are now satisfied in all directions  
$i\in I$, and therefore property (b) follows as well.  


Property (a) is still only partially satisfied, because the trivial components are not yet fully decorated, they are still missing a choice of leading coefficients  in the directions $i\in I^\infty$.
  
Next,  Corollary \ref{cor.enh.cond.2} implies property (e) as well as the fact that the refined matching conditions (\ref{enh.a.y=d.full}) are satisfied for each node $y$ of any intermediate curve $C\ra C' \ra C_0$ and any direction 
$i\notin I^\infty(y)$. In particular, at each node 
$x\in C_0$ 
\bear\label{match.a=c.d}
 a_i(x^-) a_i(x^+) c(x)^{s_i(x)}=1
\eear
for all directions $i\in I^\pm(x)$, and more generally for the node $y_r$ of the intermediate curve $C_r$ obtained by contracting the first $r$ components of $B_x$, using the notation (\ref{B=ordered})
\bear\label{match.a=a.r}
a_i(x^-) a_i(x_{r}^+)\cdot c(y_r)^{s_i(x)}=1
\eear
for all $i\notin I^\infty(x_{r}^+)$. If for $i\in I^\infty(x_{r}^+)$ we {\em define} $ a_i(x_{r}^+)$ by the same formula (\ref{match.a=a.r}), then this uniquely decorates the trivial components such that (\ref{match.a=a.r}) will now hold in all directions $i \in I$. Dividing two consecutive equations (\ref{match.a=a.r}) and using the fact that the leading coefficients are reciprocal of each other at the two end points of a trivial component proves (d) and also completes the proof of properties (a)-(f).  \end{proof} 
\begin{rem} The trivial components play a special role in the relative theory. Though typically multiple covers of degree $\al(\Si)$ defined by \eqref{multipl.s} and with a $\cx^*$ automorphism group, we can completely classify them and understand their moduli spaces as well as how they interact with the rest of the components. 

There are two fundamentally equivalent ways in which one can deal with the trivial components: one way is to completely eliminate them (once the topological implications of their presence on $\wt s$ are completely understood); the refined matching conditions reduce to equations \eqref{match.a=c.d}, indexed by each node $x$ of $C_0$, and involving only leading coefficients of the nontrivial components. These conditions are cut transversely, but are  cumbersome to work with, because they involve leading coefficients at two different points $f(x_\pm)\in X_m$, stretched across several levels (previously joined by the zig-zagging chain of trivial components that was forgotten). 

The second approach, which we took instead in Theorem \ref{thm.exist.lim} is to decorate the trivial components with an extra amount of information, and include them in the relative moduli space 
$\M_{\wt s}(\wt X_m, \wt D_m)$, as they now carry a full set of contact information to the total divisor 
$\wt D_m$ and have the expected dimension; the extra decorations then enter in a much simpler collection of refined matching conditions (\ref{enh.a.y=d.full}), now involving all the leading coefficients and indexed by each node $y$ of the refined limit $C$. These equations are also cut transversely, and as we have seen in the proof above essentially uniquely determine all the information about the trivial components, up to a residual finite group action, the group of roots of unity of order $\al(\Si)$. 

For example, a decorated trivial component $\Si$ is determined by $\mathrm{dim}_\cx X-1$ complex parameters, which matches the expected dimension (\ref{dim}) of the relative moduli space into 
$(\F_V, D_0\cup D_\infty\cup F_V)$. Moreover, at a depth $k\ge 2$ node $y$ of $C$, the existence of a solution $c(y)$ of the equation (\ref{enh.a.y=d.full}) imposes a  $k-1$ complex dimensional constraint on the leading coefficients at $y$ (equal to the expected dimension from the proof of Lemma \ref{L.d=d-2}). Should a solution $c(y)$ exist, it is unique up to a root of unity of order $\al(y)$.  
\end{rem}

\setcounter{equation}{0}
\section{The Relatively Stable Map Compactification}\label{s6}
\medskip


We are finally ready to define a relatively stable map $f:C\ra X_m$ into a level $m$ building which summarizes the topological restrictions on the types of limits  constructed in Theorem \ref{thm.exist.lim} for  sequences of maps in $\M(X,V)$. The refined limit $f$ constructed there was described in terms of a resolution $\wt f:\wt C\ra \wt X_m$, which was an element of the moduli space $\M_{\wt s}(\wt X_m, \wt D_m)\ra \JV(X, V)$ over the space of $V$-compatible parameters.  
The resolution $\wt f$ satisfied both the naive and refined matching conditions, which can be stated using of the two evaluation maps.


Recall that $\M_{\wt s}(\wt X_m, \wt D_m)$ now also includes trivial components, as long as they are decorated, see Definition \ref{def.triv.comp.decorated}. 
The notation can easily become unmanageable, but 
 for any fixed marked point $x$ of the domain $\wt C$ there is an evaluation map 
\bear\label{ev.or.m}
\ev_x:  \M_{\wt s} (\wt X_m, \wt D_m)\ra  D_{J(x)} 
\eear
 into the corresponding stratum of $(\wt X_m, \wt D_m)$ defined by \eqref{ev.naive}, a leading order section 
\bear\label{sect.lead.m}
\si_x:  \M_{\wt s} (\wt X_m, \wt D_m)\ra E_{x, s(x)}=\ev_x^*(ND_{J(x)})\otimes _{s(x)} \L_x
\eear
defined by \eqref{sect.lead}  and a refined evaluation map 
\bear\label{enh.ev.m}
\Ev_x:  \M_{\wt s} (\wt X_m, \wt D_m)\ra  \P_{s(x)}(N D_{J(x)}) 
\eear
defined as the weighted projectivization of (\ref{sect.lead.m}) and refining (\ref{ev.or.m}). Recording this information for all the contact points of $\wt C$ is summarized by the following diagram: 
\bear\label{ev.or.m.2}
\xymatrix{
E \ar[d]&&\P_{\wt s}(ND)\ar[d]
\\
\M_{\wt s} (\wt X_m, \wt D_m)\ar@/_/[u]_{\si}\ar[rr]^{\quad\ev}\ar[rru]^{\Ev}&& D_{\wt s}
}
\eear

According to our conventions \eqref{M.ft.V.strat}, the strata of $\M(\wt X_m, \wt D_m)$ are indexed by the refined dual graph of $\wt f$, but diagram \eqref{ev.or.m.2} only uses the contact information of $\wt f$ encoded as decorations on the half edges of the dual graph of $\wt f$. For the rest of this section we use the notation $\si=(\tau, s)$ for the decorated dual graph associated to a map, where $\tau$ is the usual dual graph and $s$ records only the contact information. 

With this convention, let $\wt \si=(\wt \tau, \wt s)$ be the (refined) dual graph of the resolution $\wt f$ in \eqref{g.ref.fmm.0}. The domain $\wt C$ is obtained by resolving the collection of nodes $D$ of $C$.  For each $y\in D$ let $s_\pm(y)= \wt s (y_\pm)$ denote the contact information at the nodes $y\in D$ and let $s$ denote the contact information at the remaining points (i.e. marked points of $C$). 
Similarly, the total divisor $D_m$ in a level $m$ building was the union of the singular divisor $W_m$ where the pieces of the building were joined together (in dual pairs) and the zero divisor $V_m$.    
\begin{defn}\label{D.map.m} A {\em naive map} from $C$ into a level $m$ building  
$(X_m, V_m)$ is a continuous function $f:C\ra X_m$ such that $f$ has a resolution 
$\wt f\in \M(\wt X_m, \wt D_m)$ which satisfies conditions (b)-(c) of Theorem  \ref{thm.exist.lim}.  
\end{defn}
The naive matching condition (c) is equivalent to the fact that  $\wt f$ belongs to the inverse image of the diagonal under the evaluation map at pairs of marked points giving the nodes $D$ of $C$:
\bear\label{ev.m.stable}
\ev_{D}: \M_{\wt s} (\wt X_m, \wt D_m)\longra W_{s_-} \ti W_{s_+}
\eear
extending (\ref{ev.nodes}). Recall that the normal direction to the singular divisor $W$ come in dual pairs, and here $s_-$ denotes the contact information associated to the branch $x_-$ of $\wt C$ which includes not just the multiplicities $s_i(x^-)$ but also the levels $l_i(x^-)$, the signs $\ep_i(x^-)$ and therefore the indexing set $J(x^-)$ of the branches of the total divisor $D_m$ which record the particular strata $W_{s_-}$  of the singular divisor $W$ that $f(x_-)$ belongs to; according to our conventions this includes the case when $x_-$ is an ordinary marked point, with empty contact multiplicity to $W$.  

Condition (b) implies that all the contact points of $\wt f$ descend to special points of $C$ and $C$ has no nodes on the zero divisor $V_m$ away from the singular divisor $W_m$; therefore the contact points $x$ of $\wt C$ which do not come from nodes of $C$ must be mapped to a stratum of the zero divisor $V_m$ (away from the singular divisor $W_m$) and record the contact information $s$ of $f$ along the zero divisor. Condition (b) also implies that $C$ is obtained from $C_0$ by possibly inserting strings of trivial components $B_x$ (broken cylinders) either between two branches $x_\pm$ of a node $x$ of $C_0$ or else at a marked point $x_0$ of $C_0$; the fact that the multi-signs $\ep$ are opposite at each node by (c) implies that each chain $B_x$  moves in a monotone zig-zagging fashion in the fiber of $X_m$ over $f_0(x)$, exactly as described in Corollary \ref{C.triv.forget}  (note that the level changes only in those directions in which the contact information is geometric). 

These describe restrictions on the types of refined dual graphs a naive map into $X_m$ could have. The asymptotics in part (e) of Theorem \ref{thm.exist.lim} impose one more condition on the dual graph $\si$ of the limit $f$: the system \eqref{full.enh.log.arg} must have strictly negative solutions. Denote by ${\mathcal T}_{V}$, perhaps more appropriately denoted ${\mathcal T}_{(X,V)}$ the collection of refined dual graphs with these properties, and call them {\em balanced} graphs.  This is a well defined notion, because the system depends {\em only} on the dual graph $\si$, and not on the leading coefficients of $f$. It has one variable for each positive level $l$ and  for each positive depth node $y$ of $C$ (the equations are independent of the depth zero nodes). 
Let $\wt T_\si\le (\cx^*)^{m}\ti (\cx^*)^{|D_+|}$ denote the complex torus generated by the solutions of this system, and $T_\si\le  (\cx^*)^{m}$ its projection onto the first factor. 
When $m=0$ the balanced conditions are empty, so any decorated dual graph is by convention balanced. When $m>0$ and $\si \in \mathcal T_{V}$ the torus $T_\si$ is positive dimensional. 

Note that the resolution $\wt f:\wt C \ra \wt X_m$ of a naive map $f:C\ra X_m$ is unique only up to reordering of the extra decorations, i.e. of the {\em choice} of an order for extra points over the nodes of the domain and of an indexing of the branches of the divisor over them. Let $G_\si$ denote the corresponding  product of symmetric groups from Remark \ref{R.aut.s} acting on these decorations.

\begin{defn}\label{D.enh.cond.f} A {\em map} into a level one building is a naive map $f:C\ra X_m$ as in Definition \ref{D.map.m} which has a resolution $\wt f$ satisfying the {\em refined matching conditions}, i.e.  (i) $\wt f\in \Ev^{-1}(\De^\pm)$ is 
in the inverse image of the antidiagonal $\De^\pm$ under the refined evaluation map at pairs of marked points $y^\pm$ corresponding to the nodes 
$D$  of $C$:
\bear\label{ev.t.m}
\Ev _{D}: \M_{\wt s} (\wt X_m, \wt D_m) &\ra& \prod_{y\in D}  \P_{\wt s(y)} (N W_{J_-(y)})\ti  \P_{\wt s(y)} (N W_{J_+(y)}) 
\eear   
and (ii) their refined dual graph is balanced, i.e. $\si\in \mathcal T_V$. 

Denote by $\M_\si(X_m, D_m)$ the collection of such maps $f$ (up to automorphisms of the domain) or equivalently of their resolutions $\wt f\in \Ev^{-1}(\De^\pm)$ (up to automorphisms of the domain and reordering of the decorations). 
\end{defn} 
For each $\si$ we get a stratum $ \M_{\si}(X_m, D_m)= \Ev^{-1}(\De^\pm) /G_\si$ defined 
via its resolution 
\bear\label{D.M.m}
\M^{\si}(X_m, D_m)= \Ev^{-1}(\De^\pm) 
\eear 
Note that the level $m=m(\si)$ of the building $X_m$ can be read off the refined  dual graph. The refined matching condition extends (\ref{enh.ev.match}), and keeps track not only of the image of $\wt f(y_\pm)$ in the singular divisor $W$  but also on its leading coefficients $(a_i(y_\pm))_i$  as elements of two dual normal bundles $N W_{J(y_-)} \cong  N^* W_{J(y_+)}$. It simultaneously restricts the topological type of $f$  to a balanced graph $\si\in \mathcal T_V$. 

From now on, by a resolution of map into a level $m$ building we mean one that satisfies the refined matching condition. Note that the level zero  $\M(X_0, V_0)$ is precisely the space $\M(X, V)$ as defined in \S \ref{s2}. Even though the $(\cx^*)^m$ action rescaling the target $X_m$ induces an action on the collection of naive maps into $X_m$, only the subtorus $T_\si\le (\cx^*)^m$ acts on their resolutions (i.e. preserves the refined matching conditions). In effect, by imposing the refined matching conditions we constructed a partial slice to the $(\cx^*)^m$ action.

\begin{defn} \label{D.map.m.stable} A map $f:C\ra X_m$ into a level $m$ building (as in Definition \ref{D.enh.cond.f}) is called {\em relatively stable} if it has at least one nontrivial component in each positive level $l$. 

Let $\ov \M(X, V)$ denote the collection of relatively stable maps $f:C\ra X_m$, up  to reparametrizations of the domains and rescaling of the target by the $(\cx^*)^m$ action. Equivalently,  it is the collection of resolutions $\wt f$ of relatively stable maps, up  to reparametrizations of the domains, rescaling the target by $T_\si$ and reordering the decorations by $G_\si$. 
\end{defn}
A relatively stable map $f$ and its resolutions $\wt f$ have a finite automorphism group. As before, the moduli space $\oM(X, V)$ comes with a stratification 
\bear\label{def.m.tau}
\oM(X, V)\ma \longra^\tau \mathcal T_V \ma \longra^m\N
\eear
that sends $f:C\ra X_m$ to its refined dual graph $\tau(f)=\si$ and then to its level $m$. Each open stratum $\M_\si(X,V)$ similarly comes with a finite resolution
\bear\label{m.s.tilde.}
\M^\si(X, V)=\M(X_m, D_m)/T_\si= \Ev^{-1}(\De^\pm)/T_\si
\eear defined as the collection of (decorated) resolutions $\wt f:\wt C\ra \wt X_m$ up to automorphisms of the domains and the complex torus $T_\si$ action on the target. The finite group $G_\si$ acts on $\M^\si(X, V)$ reordering the decorations and the quotient is $\M_\si(X,V)$. When $m=0$, $T_\si=G_\si=1$ so the level zero part of the moduli space $\oM(X, V)$ is again $\M(X, V)$. Note that diagram 
\eqref{M.ft.V.strat} is now fully extended to the compactification $\oM(X, V)$. 
 
\begin{rem} In the case $V$ has several components and we decide to rescale in $c$ independent directions, then we have a $(\cx^*)^{m_1} \ti \dots\ti (\cx^*)^{m_c}$ action on a multi-level  $m=(m_1, \dots m_c)$ building $(X_{m}, V_{m})$ that we take the quotient by. A map $f:C\ra X_m$ into a multi-level $m$ building is then called {\em relatively stable} if each multilevel $l=(l_1, \dots l_c)$ different from (0, \dots 0) contains at least one nontrivial component (or equivalently $f$ has finite automorphism group). 

The notion of stability therefore depends in how many independent directions we rescaled the target in, and so on  the particular group action we are taking the quotient by. For example, a  nontrivial component in the level (1,2) of a multi-directional building counts  as a nontrivial component in both level 1 and also in level 2 if we regard it as part of a uni-directional building, see Example \ref{neck.1}(b). If there are several independent directions, there is always a projection (stabilization) map from the unidirectional relatively stable map compactification $\ov \M(X, V)$ constructed in this paper to a smaller multi-directional relatively stable map compactification obtained by collapsing some multi-levels containing only trivial components (when regarded as independent multi-levels). 

There is also a forgetful map that forgets some of the components of $V$, and in particular
\bear\label{ft.v.final}
\mathrm{ft}_V:\ov\M(X, V)\ra \ov\M(X)
\eear 
maps a relatively stable map $f:C\ra (X,V)_m$ into its contraction $f_0:C_0\ra X$. Equivalently $f_0$ is the pushforward $p_* f$ under the collapsing map $p:X_m\ra X$.\end{rem} 
With these definitions, Theorem \ref{thm.exist.lim}  implies: 
\begin{theorem}\label{thm.lim.ex.2}  Consider the universal moduli space $\ov\M(X,V)\ra \JV(X, V)$ over the space of $V$-compatible parameters. Assume $f_n:C_n\ra X_{m_n}$ is a sequence of $J_n$-holomorphic maps in $\ov\M(X,V)$ such that (i) $J_n\ra J_0$ in $\JV(X, V)$ and (ii) the projections $\ft_V(f_n)\ra f_0$ in $\oM(X)$.  

Then there exists a level $m$ building $X_m$, a family $\x\ra B$ of deformations of it, a  sequence of parameters $\la_{n}\ra 0$ in $B$ and identifications 
$R_{\la_n}: X _{m_n}\ra X_{\la_n}$ such that after passing to a subsequence, $R_{\la_n} f_n$ converges (in the sense of Definition \ref{D.Gr.top}) to a relatively stable $J_0$-holomorphic map $f:C\ra X_m$ which is a refinement of $f_0$. 
\end{theorem}
\begin{proof} We first show that the number of levels $m_n$ as well as the rest of the topological type of the resolutions $\wt f_n\in \M_{\si_n}(\wt X_{m_n}, \wt D_{m_n})$ is uniformly bounded, thus can be assumed constant after passing to a subsequence. 

The number of nontrivial components of $f_n$ is the same as the number of components of its projection $p_* f_n$, thus bounded by the number $K$ of components of its limit $f_0$.  Since $f_n$ are relatively stable maps, they have a nontrivial component in each positive level. Each component can be in at most $\dim X$ different levels so the number of levels $m_n \le M=K \dim X$.  The number of chains of trivial components of $f_n$ is bounded by the number of nodes of its domain (thus of $C_0$) and the stretch  (\ref{D.stretch}) of each trivial chain by $M\dim X$; thus the number of trivial components of $f_n$ is bounded. This proves that there are finitely many possibilities for the  dual graphs of $f_n$, so after passing to a subsequence we can assume $\si_n=\si$ is constant. So $m_n=m$ and all the domains of $f_n$ are homeomorphic to a fixed nodal Riemann surface $\Si$, with a fixed smooth resolution $\wt \Si$ and attaching map $\xi:\wt \Si\ra \Si$.  

Next Theorem \ref{thm.exist.lim} applies to the resolutions $\wt f_n$ of $f_n$, which are stable $J_n$-holomorphic map into the smooth manifold $(\wt X_m, \wt D_m)$, to produce a relatively stable limit $g$ into a level $k$ building over $(\wt X_m, \wt D_m)$. But the target of $g$ after attaching the zero and infinity divisors according to the map $\xi:(\wt X_m, \wt D_m)\ra  (X_m, D_m)$ is nothing but  a level $m+k$ building over $(X, V)$. Therefore $f_n$ converges to a relatively stable limit $f$ into an  $m+k$ building. \end{proof} 
As we have seen, this notion of convergence defines a Hausdorff topology on the universal moduli space $\ov\M(X,V)\ra \JV(X, V)$ in which $\st\ti \Ev$ and $\ft_V$ are continuous. The maps in \eqref{def.m.tau} are stratifications (USC), while the composition 
\bear\label{D.h.ti.s}
\oM(X, V)\ma \longra^\tau \mathcal T_V \xrightarrow{\mathsf{h}\ti \mathsf{s}} H_2(X)\ti \Z\ti \Z \ti \mathcal S
\eear
is locally constant (continuous). Here $\mathsf{h}(\si) = (A_\si, \chi_\si, \ell(\si))$ and 
$\mathsf{s}(\si)= (\si(x))_{x\in P_\si}$ is the ordered partition of $A_\si \cdot V$ associated to the $\ell(\si)$ marked points $P_\si$, see Remark \ref{top.info}. This gives a decomposition 
\bear
\oM(X, V)=\ma \sqcup_{s\in \mathcal S} \oM_{s}(X, V)= \ma\sqcup_{A, \chi, s}  \oM_{A, \chi, s}(X, V)
\eear 
Each $[f]\in \oM(X, V)$ is represented by maps $f: C\ra X_m$ which in turn are described via their resolutions $\wt f:\wt C \ra \wt X_m$ and projection $f_0:C_0\ra X$, see diagram \eqref{g.ref.fmm}. 

With these preliminaries, the transversality argument from Lemma \ref{L.d=d-2} and the dimension counts of Lemma \ref{dim.diff1} extend to a level $m$ building:
\begin{prop}\label{L.dim.cuts2} For every $m\ge 0$, each stratum $\M_{\wt \si}(\wt X_m, \wt D_m)\ra \JV(X, V)$ is cut transversally at points $(\wt f, \wt C, J,\nu)$ with $\Aut(C_0)=1$ and the refined evaluation map $\Ev$ is a submersion. 

Therefore for generic $V$-compatible parameter $(J, \nu)$ each open stratum \eqref{m.s.tilde.} of $\oM(X, V)$ is a smooth manifold of dimension 
\bear\label{dim.Ev.diag.m}
\dim \M^{\si}(X, V)= 
\dim \ov\M_{\mathsf{h}(\si)} (X, V) - \dim T_{\si} - 2 d_0(\si)
\eear
at any point $(f, C)$ with $\Aut (\ct (C))=1$. Here $d_0(\si)$ is the number of depth zero nodes of $\si$, and $\dim \oM(X, V)$ is \eqref{dim.M.s}. In particular, the boundary strata are at least codimension 2.   \end{prop}
\begin{rem}
Note that in the presence of higher depth strata, the codimension of  the stratum into a level $m$ building (without any depth zero nodes) is not necessarily equal to $2m$, as again illustrated  by Example \ref{neck.1} (b). The stratum as $x_1\ra x_0$ is clearly only codimension 2, even though the limit is a map into a building with 2 levels. The reason is that the refined matching conditions impose further conditions on the rescaling parameters (in that case $\la_2=\la_1$), and so the rescaling parameters are no longer independent variables. The complex codimension of the stratum is the number of independent rescaling parameters (i.e. the dimension of the torus $T_{\si}$). 
\end{rem}

\begin{rem} If $J$ is integrable near $V$ then the weighted projective space 
$\P_{s(x)}(NV_{I(x)})$ can be regarded as an exceptional divisor in the blow up of the target at $V_{I(x)}$, and the refined evaluation map is the usual evaluation map into this exceptional divisor, relating it with the approach in Davis' thesis \cite{da} that worked very well in genus zero (but did not extend in higher genus). The only difference here is that this blow up is now a weighted blow up, which seems to keep better track of what happens in the limit when the multiplicities $s_i(x)$ are not equal, especially in higher genus. Of course, being a weighted blow up, it has singular strata, but with only orbifold singularities.
\end{rem}

\setcounter{equation}{0}
\section{The relative GW invariant}\label{s7} 
\bigskip

Theorem \ref{thm.lim.ex.2} implies that the relative moduli space  $\ov\M_{A, \chi, s}(X, V)$  for fixed parameter or over a fixed path of parameters  is compact. Proposition \ref{L.dim.cuts2} implies that for generic parameter each stratum of $\ov\M(X, V)$ has a smooth resolution, and all  boundary strata are codimension at least two (technically speaking, we only proved this under the simplifying assumptions of Remark \ref{nu.univ}). This is enough to imply that the image $(\st\ti \Ev)( \ov\M_{A,\chi,s}(X, V))$ represents a homology class, as described in  Sections 7 and 8 of \cite{ip1}.  In particular, we obtain the following generalization of Theorem 8.1 of \cite{ip1}:
\begin{theorem}\label{main.thm}  Assume $V$ is a normal crossing divisor in $X$ for some $V$-compatible pair $(J, \omega)$. For each fixed parameter $J\in \JV(X, V)$ the space of relatively stable maps $\ov\M_{A, \chi, s}(X,V)$ is compact and comes with a continuous map 
\bear
\st\ti \Ev: \ov \M_{A, \chi, s} (X, V) \ra \ov \M_{\chi, \ell(s)} \ti \P_{s} (N V)
\eear
Under the assumptions of Remark \ref{nu.univ}, for generic $V$-compatible $(J,\nu)$ the  image of $\ov \M_{A, \chi,s}(X,V)$ under $\st\ti \Ev$  defines a homology class $GW_{A, \chi,s}(X, V)$ in dimension \eqref{dim.M.s}. 

The class $GW_{A, \chi, s}(X, V)$ is independent of the perturbation $\nu$ and is invariant under smooth deformations of the pair $(X, V)$ and of $(\omega, J)$ through $V$-compatible structures;  it is called the  GW invariant of $X$ relative the normal crossing divisor $V$. 
\end{theorem} 

When $V$ is smooth, $\Ev$ is the usual evaluation map $\ev$ into $V_s$, so combined with Example \ref{V.sm} gives
\begin{cor} When $V$ is a smooth symplectic codimension 2 submanifold of $X$, the relative GW invariant constructed in Theorem \ref{main.thm} agrees with the usual relative GW invariant $GW(X,V)$ as defined in \cite{ip1}. 
\end{cor}

\begin{rem} \label{orbi.m.s} ({\bf Resolutions and Gluing}) Under the  assumptions of Remark \ref{nu.univ}, the gluing formula of \cite{ip2} can also be extended to this case to prove that for generic $J\in \JV(X, V)$ the local model of $\ov \M_s(X,V)$ normal to a boundary stratum is precisely described by the refined matching conditions  i.e. all solutions of the refined matching conditions glue to give actual $J$-holomorphic solutions in $X_\la$, providing the converse of the limiting argument of Theorem \ref{thm.exist.lim} for such generic parameter.  Below we outline the main ingredients involved, with the full details to appear in \cite{i-nor-sum}.  

First, the space of solutions to the refined matching conditions is typically not smooth at $\la=0$, but has smooth resolutions. This is already the case when the divisor $V$ is smooth. In that case the refined matching conditions 
\bear\label{enh.node.match}
a_i(x_-)a_i(x_+)\mu(x)^{s_i(x)}=\la
\eear
are automatically satisfied, i.e. given any naive map into a building, and any gluing parameter $\la$ of the target one could always find  gluing parameters $\mu(x)$ at each node $x$ of the domain satisfying \eqref{enh.node.match}. In fact, there are $s(x)$ different choices for each node, thus the multiplicity in the gluing formula, and the source of the branching in the moduli space. The easy fix is to include a choice of the roots of unity separating the different choices of the gluing parameters $\mu(x)$ at each node $x$, as explained in \cite{i2}. 


When $V$ is a normal crossing divisor, the story is similar. We saw that the local model of $\oM(X, V)$ near a point $f_0\in \M_\si(X, V)$ can be described by a collection of data $(f, \mu, \la)$, where 
$f\in \M_\si(X, V)$ is near $f_0$,  $(\la_1, \dots \la_m) \in (N\otimes N^*)^m\cong \cx ^m$ are gluing parameters of the target, and $\mu \in \ma\bigoplus_{x\in D} \L_{x_+} \otimes\L_{x_-}$ are the gluing parameters of the domain $C$, including those on the trivial components, which are all decorated, see (\ref{def.dec}).  There parameters must also  asymptotically satisfy for  $(\la, \mu)$ near $(0,0)$ the full set of refined matching conditions cf. \eqref{a.mu=la.before.log}:
\bear\label{match.enh.m.n}
a_i(y_-) a_i(y_+) \mu(y)^{s_i(x)} =\La(l_i(y_\pm))
\eear
for each $i\in I^\pm(y)$ and for each node $y$ of any intermediate curve $C\ra C'\ra C_0$. 

As we have seen, the existence of a family $\la_l\ra 0, \; \mu(z)\ra 0$ of solutions to the typically overdetermined system (\ref{match.enh.m.n}) is not automatic, and therefore imposes both  conditions on the leading coefficients of $\wt f$ as well as combinatorial conditions on the topological data $\wt \si\in \mathcal T_V$ associated to  $\wt f$.
The locus of the equations (\ref{match.enh.m.n}), thought as equations in the parameters $(\la, \mu)\in \cx^{m}\ti \cx^{|D_+|}$ is smooth for $(\la, \mu)\ne 0$, but may be singular at $(\la, \mu)=0$ (it is only a pseudo-manifold, or a branched manifold with several branches coming together at 0). We can instead use a refined compactification over $(\la, \mu)=0$ which has orbifold singularities, or equivalently work instead with the linear system of equations (\ref{full.enh.log}) obtained after taking $\log$ of these 
equations, i.e. consider a cylindrical (perhaps more appropriately called toroidal) compactification,  as in the last paragraph of \S 4 in \cite{ip2}. 

This is essentially what we did in the proof of Theorem \ref{thm.exist.lim} to get the refined matching conditions (\ref{enh.a.y=d.full}) in (d) for the gluing parameter  $\la=1$ together with the asymptotics in part (e). A brief inspection of equation (\ref{enh.a.y=d.full}) shows it is precisely the inverse image under $\Ev$ of the antidiagonal $\De^\pm$ in the weighted projective space $\P_{s(y)}NW_{J(y)}$, including its multiplicity $\al(x)$ in (\ref{multipl.s}) contributing to the order of the (finite) isotropy group $\Aut (f)\le T_\si$ of $f$ (whenever $\Aut C=1$).  

 As we have also seen in the proof of Theorem \ref{thm.exist.lim}, 
we could instead eliminate the trivial components, replacing the equations (\ref{match.enh.m.n})  with the following subset of equations: 
\bear\label{match.enh.m-trivial}
a_i(x^-) a_i(x^+)  \mu(x)^{s_i(x)}=\prod_{l=l_i^-(x)} ^ {l_i^+(x)} \la_{l}. 
\eear
with one equation in  each direction $i\in I(x)$ and for each node $x$ of $C_0$. Starting with a solution of (\ref{match.enh.m-trivial}) which involves only leading coefficients information at the nontrivial components of $f$, we can always uniquely solve the equations (\ref{match.enh.m.n}) to find all the leading coefficients of the (decorated) trivial components; the solutions will again have the correct asymptotics as  $\la,\mu\ra 0$  provided $\si\in \mathcal T_V$. 

No matter which of these models for the matching conditions we use, the approximate gluing map $\Gamma$ extending that in  Definition 6.2 of \cite{ip2} takes any triple $(f, \mu, \la)$ satisfying the refined matching conditions into an approximately $J$-holomorphic map $\wh f_{\mu, \la}: C_\mu\ra X_\la$ (without the refined matching condition the image would not land in $X_\la$). When the stratum containing $f_0$ is cut transversally at $f_0$ (i.e. the linearization is onto),   one can show that for any small fixed gluing parameters each approximate solution  $\wh f_{\mu, \la}$ can be uniquely corrected into an actual solution $f_{\mu, \la}: C_\mu'\ra X_\la$, providing the local model for the moduli space $\ov\M(X, V)$ near $f_0$. The key ingredients in the proof are the uniform estimates in Sections 6-9 of \cite{ip2} and especially the first eigenvalue estimate in the weighted Sobolev norms used there, all of which can be extended in a straightforward way to the case $V$ is a normal crossing divisor by working semi-locally in the necks of $C_\mu$ and $X_\la$ (as long as (\ref{match.enh.m.n}) are satisfied), as we have done in Sections 4 and 5 of this paper (see \cite{i-nor-sum} for more details).

\end{rem}
\subsection{Further directions} 
\smallskip

The next question is how the GW  invariants relative normal crossing divisors behave under degenerations. The degenerations we have in mind  come in several flavors. 

The first type of degeneration is one in which the target $X$ degenerates, the simplest case of that being the degeneration of a symplectic sum into its pieces. This comes down to the symplectic sum formula proved in \cite{ip2}, but where now we also have a divisor going through the neck. Consider for example the situation 
\best
(X,V)=(X_1,V_1)\#_U(X_2,V_2)
\eest
described in Remark  \ref{sym.sum.ex}, which means that $X=X_1\#_U X_2$, and simultaneously the divisor $V$ is the  symplectic sum of $V_1$ and $V_2$ along  their common intersection 
with $U$. The relative GW of the sum $(X,V)$ can be expressed in terms of the relative GW invariants of the pieces $(X_i, V_i\cup U)$; this type of formula allows one for example to compute the absolute GW invariants of a manifold obtained by iterating the symplectic sum construction. 
\begin{ex} Assume $V$ is a  normal crossing divisor in $X$, 
and let $X_1$ be the level one building associated to $(X,V)$, with its singular divisor $W_1$. Just as was the case in \cite{ip2}, the  rescaling process in this paper that constructs a family of deformations $X_\la\ma\longra_{\la\ra 0}X_1=X\cup_V \P_V$ can be reinterpreted as the inverse operation to the trivial symplectic sum 
\bear\label{triv.deg.dec} 
X=\ma \#_{W_1} \wt X_1= \ma\#_{V} ( X\sqcup \P_V)
\eear
of $X$ and $\P_V$  along $V$, or more precisely along the singular divisor $W_1$. The symplectic sum formula relates the GW invariants of the sum to the relative invariants of the pieces. Its complete proof will appear in \cite{i-nor-sum}, and follows by comparing the limit $\lim_{\la\ra 0}\oM(X_\la)$ to the moduli space $\oM(X_1, W_1)$ of the central fiber over $\la=0$. The later is defined just as in Definition \ref{D.enh.cond.f} via its resolutions $\Ev^{-1}(\De)\subset \oM(\wt X_1, \wt W_1)$ obtained by separating the nodes over the singular divisor. The trivial sum formula takes the form: 
\bear\label{sympl.sum.GW}
GW(X)=\sum_{s}  \frac{\al(s)} {|\Aut(s)|} \; GW_s(\wt X_1, \wt W_1)\cap \Ev^*\eta_{\De_s} \ma=^{\mathrm{def}}\;GW(X, V)\ma*_{\Ev} GW(\P_V, V_{\infty}\cup F_V)
\eear
where $\eta_{\De_s}$ denotes  the Poincare dual of the antidiagonal $\De_s^\pm$ in the weighted projective space $\P_s NW\ti \P_s N^*W$ (which satisfies Poincare duality over $\Q$),  and $s$ ranges over all balanced graphs $\mathcal T_{X_1, W_1}$. The multiplicity 
$\al(s)$ comes from the different choices of roots of unity in the refined matching condition \eqref{match.enh.m-trivial}, while $\Aut(s)$ comes from reordering the contact information associated to the nodes along $W$.  Note that as before, passing from $f:C\ra X_1$ to its resolution $\wt f:\wt C\ra \wt X_1$ requires choosing an {\em ordering} of the contact information $s$ along the singular divisor as in Remark \ref{R.aut.s}. With this, $s$ can be described via its resolution $\wt s\in \mathcal T_{\wt X_1, \wt W_1} = \mathcal T_{X, V}\ti  \mathcal T_{\P_V, V_{\infty}\cup F_V}$, up to the the (anti)-diagonal action reordering the extra choices.
\end{ex}

But there are other types of symplectic sums/smoothings of the target that these relative GW invariants should enter. The next simplest example is either the 3-fold sum or 4-fold sum defined by Symington in \cite{s3} (see also \cite{s4}). Both these constructions should have appropriate symplectic extensions to higher dimensions involving smoothings $X_\la$ of a symplectic manifold 
$X$ self intersecting itself along a symplectic normal crossing divisor $V$. The sum formula would then express the GW invariants of $X_\la$ in terms of the relative GW invariants of $(X,V)$. A special case of this is what is called a stable degeneration in algebraic geometry, in which case one has a smooth fibration over a disk with smooth fiber $X_\la$ for $\la \ne 0$ and whose central fiber $X_0$ has normal crossing singularities.  

\medskip

There is also a related question when the target $X$ is fixed, but now the divisor $V$  degenerates in $X$. The simplest case of that 
is the one in Example \ref{ex.cp2}, and serves as the local model of more general deformations. For example, a slightly more general case would be when we have a family of smooth divisors $V_\la$ degenerating to a normal crossing one $V_0$, which let's  assume has at most depth 2 points (i.e. its singular locus $W$ is smooth).  After blowing up $W$, this case can be reduced  to the case of a symplectic sum of the blow up of $X$  with a standard piece $\P_W$, constructed using the normal bundle of the singular locus $W$.  The divisors now go through the neck of the symplectic sum, but their degeneration happens only in $\P_W$ (therefore involves only local  information around $W$). So  if one can understand the degeneration locally near $W$,  one can again use the sum formula to relate the GW invariants of $(X, V_\la)$ to those of $(X, V_0)$.  
\smallskip

The discussion in this paper should also extend to the case when the target $X$ has orbifold singularities and the normal crossing divisor $V$ itself, as well as its normal bundle has an orbifold structure.  In this case the domains of the maps should also be allowed to have orbifold singularities. Again we have a very similar stratification of the domain and of the target  $V$ but now it has more strata depending also on the conjugacy classes of the isotropy groups; the orbifold evaluation maps take that into account as well. The corresponding refined matching condition will include that information, in the form of an additional balanced condition at each node as in  \cite{a}.



\begin{appendix}

\setcounter{equation}{0}
\renewcommand{\thesection}{\Alph{section}}
\renewcommand{\thesubsection}{\Alph{section}.\arabic{subsection}}
\renewcommand{\thetheorem}{\Alph{section}.\arabic{theorem}}
\renewcommand{\theequation}{\Alph{section}.\arabic{equation}}
  
\section{}

\subsection{Stratifications associated to a normal crossing divisor}\label{S.A.nc.strat}
   
Assume $V$ is a normal crossing divisor in $(X,\omega, J)$ and that $\iota: \wt V\ra V$ is its resolution and 
$\pi:N\ra \wt V$ its normal bundle. In particular this means that we have an immersion  $\iota:(U, \wt V)\ra (X, V)$ from some tubular neighborhood  $U$ of the  zero section $\wt V$ of $N$. 

 The divisor $V$ is stratified depending on how many local branches meet at a particular point (see \S \ref{S.A.stratif} for an overview of stratifications). Denote by $V^k$ the closed stratum of $V$ where at least $k$ local branches of $V$ meet, and let   $\op V^k=V^k\setminus V^{k+1}$ be the open stratum where precisely  $k$ local branches meet so 
 \bear\label{A.tower.V}
\dots \subseteq  V^{k+1} \subseteq V^k \subseteq\dots  \subseteq V^2\subseteq V^1=V\subseteq V^0=X
 \eear
 Then $\op V^k$ is smooth, both $\omega$-symplectic  and $J$-holomorphic and its normal bundle in $X$ is modeled locally on the direct sum of the normal bundles to each local branch of $V$: 
 \bear\label{nor.vk.open}
N_{{V^k},p}= \ma\bigoplus_{q\in \iota^{-1}(p)} N_q\quad  =\ma\bigoplus _{i\in I} N_{p_i}
\eear
where  $\iota^{-1}(p)= \{ p_i\; |\;i\in I\}$ indexes the $k$ local branches of $V$ meeting at $p\in \op V^k$. When $V$  does not have {\em simple} normal crossing these local branches may globally intertwine. The global monodromy of $N_{V^k}$ is determined by the monodromy  of the restriction
\bear\label{iota=k}
\iota: \iota^{-1}(\op V^k) \ra \op V^k
\eear
which describes a degree $k$ cover of $\op V^k$, its fibers indexing the $k$ independent directions of $N_{V^k}$ at $p$.  The original map $\iota$ is not a covering over the singular locus $V^{k+1}$ of  $V^k$, but it extends as a covering over each open stratum of the resolution $\wt {V^k}$ of $V^k$. 

The following lemma follows from the local model of a normal crossing divisor:
\begin{lemma}\label{L.ap.1} The closed stratum $V^k$ of $V$ has a resolution $\wt {V^k}$ which comes with a normal crossing divisor $W^{k+1}$ corresponding to the inverse image of the higher depth stratum $V^{k+1}$:  
\bear\label{def.wtV}
\iota_k: (\wt {V^k},  W^{k+1}) \ra (V^k, V^{k+1}).
\eear
The normal bundle to $V^k$ 
\bear\label{nor.vk}
\pi: N_{V^k} \ra \wt {V^k}
\eear
is obtained as in (\ref{nor.vk.open}) from the line bundle $N\ra \wt V$ and a degree $k$ cover $\iota$ of  $\wt {V^k}$ which extends  (\ref{iota=k}). By separately compactifying each normal direction we get a $(\P^1)^k$ bundle 
\bear\label{def-fk}
\pi_k: \F_k \ra \wt {V^k}
\eear
which comes with a normal crossing divisor 
$D_{k, 0} \cup D_{k,\infty}\cup F_k$  obtained by considering together its zero and infinity divisors plus the fiber $F_k$ over the divisor $W^{k+1}$ in the base. The resolutions of these divisors come naturally identified 
\bear\label{F=D}
\wt {F_{k}}\ma \longra^{\rho_k}_{\cong} \wt{ D_{k+1,\infty}}\ma \longra_{\cong} \wt{ D_{k+1,0}}
\eear
and their normal bundles are canonically dual to each other 
\bear\label{nor.op}
N_{F_k} \cong( N_{D_{k+1}, \infty})^* \cong N_{D_{k+1}, 0}.
\eear
\end{lemma}  
\begin{proof} The local model allows us to construct the resolution $\wt{V^k} $ of the closed stratum  ${V^k}$ as a smooth manifold, obtained by separating the branches that come together to form the next  stratum $V^{k+1}$ inside ${V^k}$, and simultaneously construct the resolution $\wt {W^{k+1}}$ of the corresponding divisor $W^{k+1}$ of (\ref{def.wtV}). The model for their normal bundles is induced from $N \ra \wt V$.  

There are several slight complications when $k\ge 2$. First, there is no direct map from the resolution $\wt {V^k}$ of the stratum $V^k$ of $V$ to the depth $k$ stratum $\wt V^k$ of  $\wt V$ (over which $N$ is defined). However the resolution $\wt {\wt V^k}$ of the depth $k$ stratum of  $\wt V$ is the degree $k$ cover of $\wt {V^k}$ whose fiber can still be thought as an indexing set for  the $k$ local branches of $V$ meeting at $p$:
\bear
\wt {\wt V^k}&\ma\longra^{\iota }& \wt{V^{k}}
\nonumber
\\
\label{iota+rho}
\downarrow^{\iota_{k}}&&\downarrow^{\iota_{k}}\\
\nonumber
\wt V \supset {\wt V}^k&\ma\longra^{\iota}&{V^{k}} \subset V
\eear
Here the vertical arrows are resolution maps. 
Therefore the pullback $\iota_k^* N $ of the normal bundle $N\ra \wt V$ still induces the same description of the  normal bundle of  
$V^k$: at each point $p\in \wt {V^k}$, 
 \bear\label{n.v.p}
N_{V^k,p}= \ma\bigoplus_{q\in \iota^{-1}(p)}(\iota_k^*N)_q\quad  =\ma\bigoplus _{i\in I} N_{p_i}
\eear
where  $\iota^{-1}(p)= \{ p_i |i\in I\}$ is the indexing set for  the $k$ local branches of $V$ meeting at $p$. Again, the cover $\iota$ may have nontrivial global monodromy which will induce a global monodromy in the normal bundle $N_{V^k}$ of $V^k$. 

Similarly, the resolution $\wt {W^{k}}$ of the normal crossing divisor $W^{k}$ of $\wt {V^{k-1}}$ is also a degree $k$ cover of the resolution $\wt{V^{k}}$ of the depth $k$ stratum of $V^k$:
\bear
\nonumber
\wt {W^{k}}&\ma\longra^{\iota}& \wt{V^{k}}\\
\label{iota+}
\downarrow^{\iota_{k}}&&\downarrow^{\iota_{k}}\\
\nonumber
 \wt {V^{k-1}} \supset W^{k} &\ma\longra^{\iota_{k-1}}&{V^{k}} \subset V^{k-1}
\eear
The vertical maps are resolution maps, while the fiber of the top map corresponds to  the indexing set of the  $k$ branches of $V$; the bottom map is the restriction of $\iota_{k-1}:\wt{V^{k-1}}\ra V^{k-1}$ to the corresponding divisor. This means in particular that the upper left corners of  (\ref{iota+rho}) and (\ref{iota+}) are the same, even though the lower left corners give two different factorizations: 
\bear\label{w=v}
\wt {W^{k}} = \wt {\wt V^k}  \ma\longra^{\iota} \wt V^k 
\eear
The fiber of this map $\iota$ indexes the $k$ local branches of $V$ coming together at a point $p\in  \wt V^k $. 

Next, the  $(\P^1)^k$ bundle 
\best
\pi_k:\F_k \longra \wt {V^k}
\eest
is obtained by separately compactifying each of the $k$ normal directions to $V$  along $V^k$, see  (\ref{n.v.p}). This means that its fiber at a point $p\in \wt{V^k} $ is 
\bear\label{fib.k}
\ma \ti_{i\in I} \P(N_{p_i}\oplus \cx)
\eear
 where $I$ is an indexing set of the $k$ local branches of $V$  meeting at $p$.  Globally,  the $\P^1$ factors of (\ref{fib.k}) may intertwine.  

Finally, the fiber divisor $F_{k}$ of $\F_k$ is by definition the inverse image of the divisor $W^{k+1}$ of $\wt{V^k}$.  Therefore its  resolution $\wt{F_{k}}$ is precisely the $(\P^1)^{k}$ bundle over the resolution  $\wt {W^{k+1}}$, whose fiber is (\ref{fib.k}). Moreover, the normal bundle to $F_k$ is the pull-back of the normal bundle of $W^{k+1}$ inside $\wt {V^k}$, which itself is the pullback  of $N\ra \wt V$ by 
$\iota_k$ of diagram (\ref{iota+rho}), see also (\ref{w=v}):  
\bear
N_{F_k}= \pi_k^* N_{W^{k+1}}= \pi_k^* \iota_k^* N
\eear

On the other hand,  the infinity divisor $D_{k+1, \infty}$ is by definition  the divisor in $\F_{k+1}$ where at least one of the $k+1$  fiber coordinates $(\P^1)^{k+1}$  is 
$\infty$, so  its resolution $\wt {D_{k+1, \infty}}$ is a $(\P^1)^k$ bundle; the base of this bundle is itself a bundle 
over $\wt {V^{k+1}}$, whose fiber consists of $k+1$ points, one for each of the $k+1$ directions of $V$ coming together.  Therefore by (\ref{w=v}),  $\wt {D_{k+1, \infty}}$ is the $(\P^1)^k$ bundle over the resolution  $\wt {W^{k+1}}$ whose fiber 
is (\ref{fib.k}). Furthermore, the normal bundle to the infinity divisor is dual to the normal bundle to the zero divisor and thus it is canonically identified to the corresponding pullback of $N^*\ra \wt V$.  

This means that we have a natural identification (\ref{F=D})  as $(\P^1)^k$ bundles over $\wt {W^{k+1}}$ and also the corresponding duality (\ref{nor.op}) of  their normal bundles.  \end{proof} 

\begin{ex} \label{ex0} {\bf (Local Structure)} One can see all these different stratifications and their resolutions in the local model, when $V$ is the union of the $n$ coordinate hyperplanes in $\cx^n$,  so its resolution $\wt V$ consists of $n$ disjoint planes 
$\cx^{n-1}$. The strata $V^k$ of $V$ are given by the vanishing of at least $k$ coordinates in $\cx^n$, while the strata 
$\iota^{-1} (V^k)$ of $\wt V$ are given by the vanishing of at least $k-1$ coordinates in each one of the $n$ disjoint planes 
$\cx^{n-1}$ of $\wt V$. Therefore the resolution $\wt{V^k}$ of $V^k$ consists of $\binom n k$ 
planes 
$\cx^{n-k}$, while the resolution $\wt{\wt V^k}$  of ${\wt V}^k$ consists of 
$n\binom{n-1}{k-1}=k\binom n k$ such planes, so the map 
$\iota:\wt{\wt V^k}\ra \wt {V^k}$  of (\ref{iota+rho}) is indeed a degree $k$ cover, whose fiber  labels the $k$ planes of $V$ coming together at a point $p\in \wt {V^k}$. Finally, the divisor $W^{k+1}$ inside $\wt {V^k}$ corresponds to the coordinate hyperplanes in $\wt{V^k}$, thus its resolution $\wt {W^{k+1}}$ consists of 
$n\binom{n-1}{ k}=(k+1)\binom  {n}{k+1}$ planes $\cx^{n-k-1}$, which is the same as the resolution of 
$\wt V^{k+1}$. This explains the diagram (\ref{iota+}) and the identification (\ref{w=v}). 
\end{ex}

The resolutions $\wt {V^k} \ra V^k$ constructed by Lemma \ref{L.ap.1} are part of a tower of compatible resolutions of the strata of $(X,V)$ that we describe next.  If $\iota:\wt V\ra V$ is the original resolution of $V$, for each finite (ordered) set $I$, consider the resolution 
\bear\label{A.D.V.I}
V_I=\{\;  (p,\;  \rho) \;|\; p\in X \text{ and } \rho: I \hookrightarrow \iota^{-1}(p) \}
\eear 
The fiber of  $V_I \ra X$ at $p$ consists of all injections $\rho:I\ra \iota^{-1}(p)$, i.e. ways of indexing by $I$ some of the local branches of $V$ at $p$, with $V_\emptyset=X$ and $V_{\{ i\}}=\wt V$ the original resolution of $V$. 

The resolutions $V_I$ are stratified by depth $k=|\iota^{-1}(p)|$,  with the top stratum corresponding to $k=|I|$ (i.e. $\rho$ is a bijection). The map 
\bear\label{A.ft.rho}
\ft_\rho:V_I\ra X
\eear 
forgetting the marking $\rho$ restricts over each open stratum $V^j\setminus V^{j+1}$  of $X$ to a topological covering.  The symmetric group $S_I$ acts freely on $V_I$ by reordering $I$.  The map \eqref{A.ft.rho} is $S_I$-invariant and descends to a map $\ft: V_I/S_I\ra V^{|I|}$, inducing a factorization 
\bear\label{A.V.s.i}
\ft_{\rho}:V_I\ra V_I/S_I\ra V^{|I|}  \hookrightarrow X
\eear 
Each  $I\subseteq J$ (or any order preserving injection $I\ra J$) induces a forgetful map 
\bear\label{A.V.s.forget}
\ft: V_J \ra V_I
\eear
For each $i\in I$ the evaluation map 
\bear\label{A.ev.i}
\ev_i: V_I \ra\wt V
\eear
is defined by $\ev_i(p, \rho)=\rho(i)$. The pullback of the normal bundle  $\pi:N \ra \wt V$ defines 
\bear\label{A.nor.V.I}
N_{V_I}= \ma\oplus _{i\in I} \ev_i^* N \ra V_I
\eear 
the `normal bundle' of the resolution $V_I$. This is an $S_I$-equivariant bundle which descends to 
a bundle over the smooth quotient $V_I/S_I$ (up to reordering of the branches), which  is equal to the  bundle  \eqref{n.v.p} whenever $I$ has order $|I|=k$.

In the local model from  Example \ref{ex0}, $V_I$ consists of $\frac{n!}{(n-k)!}$ planes $\cx^{n-k}$, corresponding to an ordered choice of $k=|I|$ out of the $n$ coordinates $x_i$ whose vanishing cuts out that plane. Therefore:
 \begin{lemma}\label{A.L.V-I} To any pair $(X, V)$ consisting of a normal crossing divisor $V$ in $X$  we can associate 
\begin{itemize} 
\item a stratification $\mathsf{depth}:X \ra \N$ with closed strata $V^k$, for $k\ge 0$; 
\item a tower of stratified resolutions $\{ V_I\}_I$ of the closed strata, indexed by finite ordered sets $I$, together with forgetful maps  \eqref{A.V.s.i} and  \eqref{A.V.s.forget} and evaluation maps \eqref{A.ev.i}; 
\end{itemize}
In particular, the resolution $\wt {V^k} \ra V^k$ of Lemma \ref{L.ap.1} is  
$\wt {V^k}= V_I/S_I$ where $I=\{1, \dots, k\}$.
\end{lemma}
One can further refine the stratification of $X$ by also keeping track of the connected components of the closed strata together with the maps induced on  $\pi_0(\cdot)$. 
When $V$ has {\em simple} normal crossings, Lemma \ref{A.L.V-I}  has a simpler interpretation,  cf. Example \ref{A.ex.V.I}. 

In  general, the constructions above are functorial and extend in much more general settings. For example, from the normal model to $V$ we get even more structure:
\begin{lemma}\label{A.L.V-I.de} Assume $V$ is a normal crossing divisor in $X$ whose neighborhood is modeled by the immersion $\iota:N\dashrightarrow X$ defined on some disk bundle of $\pi:N\ra \wt V$.  For any $\de> 0$ small, this induces
\begin{itemize}
\item a stratification $\mathsf{depth}_\de:X \ra \N$ whose closed strata $U^k_\de$ are neighborhoods of $V^k$, and open strata $X_{\de}^k$;  
\item a tower of (stratified) resolutions $\iota_k:\wt {U^k_\de}\ra U^k_\de$ (equal to the identity over the top stratum, and a topological covering on the rest), embeddings $\eta_k:\wt {U^k_\de}\ma \hookrightarrow N_{V^k}$ and projections $\pi_k:N_{V^k}\ra \wt {V^k}$ giving
\bear\label{A.X.k.finally}\xymatrix{
U^k_\de & \ar[l]_{\iota_k} \ar@{^(->}[r]^{\eta_k} \wt {U^k_\de} &  N_{V^k} \ar[r]^{\pi_k} & \wt {V^k}\ar[r]^{\iota_k}& V^k
}
\eear
Diagram \eqref{A.X.k.finally} restricts over the top strata to $X_\de^k \hookrightarrow N_{V^k} \ra V^k\setminus V^{k+1}$. 
\item a tower of stratified resolutions $\{ U_{I}\}_{I}$ indexed by finite ordered sets $I$, together with compatible forgetful maps \eqref{A.ft.ij.U} and evaluation maps \eqref{A.Ev.i.N}  
\item a tower of $S_I$-equivariant embeddings $\eta_I:U_{I}^\de \hookrightarrow N_{V_I}$ and projections $\pi_I:U_I^\de\ra V_I$ which restrict to the identity on the zero section $V_I\subset U_I^\de$ and whose quotient induce \eqref{A.X.k.finally}. 
\end{itemize}
Furthermore $\mathsf{depth}_{\de'} \ge \mathsf{depth}_\de$ for each each $\de'\ge \de>0$, and  $\mathsf{depth}_{\de}$ converges to $\mathsf{depth}$ as $\de\ra 0$. 
\end{lemma}
\begin{proof} Let  $\iota_\de$ denote the restriction of $\iota$ to the  closed $\de$-disk bundle $\wt U_\de$, where $\de$ is sufficiently small. Then $\iota_\de$ defines a stratification $\mathsf{depth}_\de:X\ra \N$ whose closed strata $U_\de^k$ consist of points $x\in X$ for which $\iota_\de^{-1}(x)$ has at least $k$ points. For $\de'\ge \de$, $\iota_\de$ is the restriction of $\iota_{\de'}$, and the number of points in the fiber may drop under restriction, so  $\mathsf{depth}_{\de'} \ge \mathsf{depth}_\de$.

For each ordered set $I$, let $U_{I}^\de$ denote the corresponding set \eqref{A.D.V.I}, but with $\iota$ replaced by $\iota_\de$: 
\bear\label{A.U.I.de}
U_I^\de=\{\;  (x,\rho)\;|\;  \;x\in X \text{ and }  \rho: I\hookrightarrow \iota_\de^{-1}(x))\; \} 
\eear
As before, the sets $U_I^\de$ are stratified by the number of points $k \ge |I|$ in the fiber of $\iota_\de$; the open strata $X_\de^{I, k}$ correspond to $k$ constant and are topological covers of $X^k_\de$ (the open strata of $U_\de^k$). There are two types of forgetful maps: 
\bear\label{A.ft.ij.U}
\ft_\rho:U_I^\de \ra U_\de^{|I|} \quad \text{ and } \quad  \ft: U_J^\de \ra U_I^\de  \quad \text{ for each $I\subseteq J$ or more generally $I\hookrightarrow J$} 
\eear
For $\de$ sufficiently small, the cardinality of $\iota_\de^{-1}(x)$ is the same as that of its projection 
$\pi( \iota_\de^{-1}(x))$ to the zero section for all $x\in X$ (by a proof by contradiction, using the fact that $X$ is compact). This means we also get a projection   
\bear\label{A.proj.U.I}
\pi_I:U_I^\de \ra V_I \quad \mbox{ defined by } \quad (x,\;  \rho) \mapsto (x, \pi\circ \rho)
\eear
with a `zero' section $V_I\ra U_I^\de$ (corresponding to $\rho=\pi\circ \rho$, i.e. the image of $\rho$ lies in the zero section of $N$). We also get an `evaluation' map
\bear\label{A.Ev.i.N}
\Ev_i: U_I^\de \ra  \ev_i^*N
\eear
defined by the formula $\Ev_i(x, \rho)=(\pi(\rho(i)), \rho(i))$ where $\rho(i)\in \iota^{-1}_\de(x)\subset N_{\pi(\rho(i))} $ is an element of the fiber of $N$ over the point $\pi(\rho(i))\in \wt V$. Their product induces the map 
\bear\label{A.emb.U.I}
\eta_I: U_I^\de \ra N_{V_I}
\eear
defined by $\eta_I(x, \rho)= ((x, \pi\circ \rho),  \ma \oplus_{i\in I} \rho(i))$ and which clearly  
projects to \eqref{A.proj.U.I}. 

All these constructions are equivariant with respect to the symmetric group $S_I$ action reordering $I$, thus descend to the quotient:
\best
V_I/S_I \hookrightarrow U_I^\de /S_I\hookrightarrow N_{V_I}/S_I \ra V_I/S_I 
\eest
and restrict to the identity on the zero section; letting  $ \wt {U^{k}_\de}=U_I^\de /S_I$  for $I=\{ 1, \dots, k\}$ gives \eqref{A.X.k.finally}. \end{proof}
\subsection{Spaces of parameters for the holomorphic map equation.} \label{A.spaces.param}  Assume $(X, \om)$ is a symplectic manifold with a fixed background metric and denote by $\J(X)$ the collection of $\om$-tamed almost complex structures $J$ with the $C^\infty$ topology; occasionally it is convenient to allow $\om$ to vary, in  which case  $\J(X)$ is the space of tamed pairs $(\om, J)$ on $X$.  Each tamed pair  induces a metric $g$ by symmetrizing $\om(\cdot, J\cdot)$. For transversality purposes, one considers $\J^{l, p}(X)$  the Sobolev completion in $W^{l, p}$ with $p>2$ and $l\ge 1$. 

The graph construction of Remark \ref{transv} induces corresponding Gromov-type parameter spaces $\JV(X)\subset \J(X')$ of pairs $(J,\nu)$ indexing deformations of the $J$-holomorphic map equation, where $X'=\ov\U\ti X$. The map $t\mapsto(J, t\nu)$ is a deformation retraction 
of $\JV(X)$ onto $\J(X)$, see (\ref{forget.nu}).  

If $p:X \ra Y$ is smooth fibration, the  space 
\bear\label{A.D.J.cpt.p}
\J(X\ma \ra^p Y)\subseteq \J(X)\ti \J(Y)
\eear 
of parameters {\em compatible} with $p$ is the subset of tuples $(J_X, \om_X, J_Y, \om_Y)$ for which  $dp\circ J_X=J_Y\circ dp$ and the restriction of $\om_X$ to each fiber of $p$ is symplectic (note that this does not impose any condition on $\om_Y$ other than it  tames $J_Y$). 

 \smallskip

Assume next $V$ is a normal crossing divisor in $(X, \om)$. A pair $(J, \nu) \in \J(X)$ is  {\em $V$-compatible} if  the following three conditions on their 1-jet along $V$ are satisfied (cf Definition 3.2 of \cite{ip2}):
\begin{enumerate}
\item[(a)] $J$ preserves $TV$ and $\nu^N|_V=0$;
 \end{enumerate}
 and for all $\xi\in N_V$, $v\in  TV$ and $w\in TC$: 
 \begin{enumerate}
\item[(b)] $[(\nabla _\xi J + J \nabla _{J\xi} J )(v)]^N= [(\nabla_v J)\xi + J (\nabla_{Jv} J)\xi]^N$
\item[(c)] $[(\nabla_\xi \nu+ J \nabla_{J\xi} \nu)(w)]^N = [ (J \nabla _{\nu(w)} J)\xi]^N$
 \end{enumerate}
Here $\xi \ra \xi^N$ is the projection onto the normal bundle $N_V$ of $V$ using the  splitting 
$\iota^*TX= T\wt V \oplus N_V$ coming from the local model at $V$. Denoting by $\JV(X, V)\subseteq \JV(X)$ the space of such $V$-compatible pairs, diagram (\ref{forget.nu}) restricts to: 
\bear\label{forget.nu>V}
\xymatrix{
\JV(X, V)\ar[r]&\ar@ /_1pc/ [l]   \J(X, V)
}
\eear
For any $J\in \J(X, V)$, $V$ is a normal crossing divisor in $(X,\om, J)$ cf Definition \ref{D.ncd}. In particular the restriction of $J$ along $V$ is an almost complex structure on 
$TX|_V\ra V$, or more precisely an almost complex structure on the resolution $\wt V$ which matches along the singular locus, i.e.  it descends to $V$. The restriction to $TV$ is an element of $\J(V)$ by (a)  while the restriction to the normal direction $N_V$ is a complex structure on $N_V$. There are similar parameters spaces $\J^1(N_V, V)$ defined by the restriction of the 1-jet of $J$ along $V$ that enters in condition (b). With this notation, we have projections 
\bear\label{j.tx.tv}
{\mathcal J}(X, V) \longra  {\mathcal J}^1(N_V, V) \longra  {\mathcal J}^0(N_V, V)\ra {\mathcal J}(V) 
\eear
that send  a $V$-compatible $J$ on $X$ to its 1-jet normal to $V$, and then to its restriction to $V$. Intrinsically, all these maps are defined by pullback via the immersion 
$\iota:(N_V,\wt V) \ra (X, V)$ which globally models $V$ and its tubular neighborhood. Using this diagram we can prove:
\begin{lemma} Assume $V$ is a normal crossing divisor in $(X, \om, J)$ as in Definiton \ref{D.ncd}. Then the collection of $V$-compatible parameters $\J(X, V)\ma \sim_{h.e.}\JV(X, V)$ is nonempty. 

If the branches of $V$ are $\om$-orthogonal, then the collection 
$\J_{cpt}(X, V)\subset \J(X, V)$ of $\om$-compatible parameters is nonempty. 
\end{lemma}
\begin{proof} We need to show that the existence of a $J$ satisfying condition (a) is enough to guarantee the existence of a $V$-compatible $J$, i.e. one that also satisfies the condition (b) on its 1-jet along $V$. 
The first map  in (\ref{j.tx.tv}) is surjective, and the arguments in the  Appendix  of \cite{ip1} show that the fiber of the second map is contractible, thus nonempty. Its target is $\J^0(N_V, V)$ which is nonempty by definition since $V$ is a normal crossing divisor. 

When the branches of $V$ are symplectic and orthogonal wrt $\omega$, one can first construct 
locally, in each coordinate patch around a depth $k$ point $p$ of $V$ a metric $g$ for which  $(TV_i)^{\om}=(TV_i)^{\perp g}$ for each local branch $i\in I$ of $V$ at $p$. If we denote $N_i=(TV_i) ^{\om}$, then this intrinsically defines a bundle $N_V$ over the resolution of $V$, as the fiber of 
$\wt V \ra V$ at $p$ is precisely $I(p)$, the set indexing the $k$ local branches of $V$ at $p$.  The metrics can be patched together to give an $\om$ compatible metric on $TX|_V=TV\oplus N_V$ and an $\om$ compatible $J$ on $X$, and which preserves both the branches of $V$ and the fibers of $N_V$. This uses the standard homotopy argument (based on the fact the space of metrics is 
convex), performed in a manner respecting the normal crossing structure of $V$ to prove that $\J^0_{cpt}(N_V, V)\ne \emptyset$. Next, when starting with an $\om$-compatible structure around $V$,  the deformation argument in the  Appendix  of \cite{ip1} gives rise to a an 
$\om$-compatible $J$ in the tubular neighborhood of $V$ satisfying the 1-jet condition there, which can then be extended to an $\om$-compatible $J$ in $\J(X, V)$. 
\end{proof}
\begin{ex} Assume $L \ra V$ is a complex line bundle over a smooth symplectic manifold $V$. The projectivization $p:\P(L\oplus \cx) \ra V$ is a Hamiltonian $\P^1$ fibration over $V$, and it comes with a zero and an infinity section (both of them smooth divisors). Then the space of parameters compatible with both the fibration and the zero and infinity section is nonempty. In fact, we can even construct parameters which are also compatible with the $S^1$ action on the fibers:  the existence of an $\om_X$ follows from Thurston's theorem while a compatible $J_X$ can be constructed starting from any $\om_V$ compatible $J_V$ on $V$ using the complex structure of $L\ra V$. 
\end{ex}

The space of parameters on a level $m$ building $(X_m, V_m)$ compatible with the total divisor is defined as the subset 
\bear\label{A.D.comp.div}
\J(X_m, V_m)\subseteq \J(\wt X_m, \wt D_m)\subseteq \J(\wt X_m, \wt W_m)
\eear
of pairs $(\om, J)$ on the  resolution $(\wt X_m, \wt D_m)$ that match along the singular locus $W_m$. It has a subset consisting of parameters that are also compatible with all the collapsing maps, or those that are  $(\cx^*)^m$-equivariant, and all come with corresponding Gromov-type extensions $\JV$. The restriction to level zero induces a fibration
\bear\label{A.resc.map.m}
\xymatrix{
\JV(X_m, V_m)\ar[r]&\ar@ /_1pc/ [l]   \JV(X, V)
}
\eear
and the rescaling process defines a family of sections  $(\om_\ep, J)$ in it. So $\JV(X_m, V_m)$ can also be thought as an extension of the parameter space $\J(X, V)$ indexing deformations of the equation (\ref{eq.f=nu}) for a map $f:C\ra X_m$, giving rise to a universal moduli space 
$\ov\M(X_m, V_m)\ra \JV(X_m, V_m)$ consisting of maps into $X_m$ satisfying the refined matching conditions (\ref{ev.t.m}) along the singular locus.  As mentioned in Remark \ref{R.ext.J.X.m}, these parameter spaces come in various flavors, depending how much of the structure of 
$X_m\ra X$ they are required to preserve. When restricted to the subspace of parameters which are invariant under the $(\cx^*)^m$ action on the positive levels, we also get a  corresponding `rubber'  moduli space 
 ${\ov\M}^\tau(X_m, D_m)/ T_\tau$ that appears in \eqref{m.s.tilde.} and which describes the level $m$ stratum of the relative moduli space $\ov\M(X, V)$. 
\begin{rem}
 When using Banach norms on the parameter spaces $\JV$ (rather than the Frechet  topology $C^\infty$),  the process of rescaling loses $m$ derivatives: e.g. the `section' above maps $\JV^{l+m,p}(X,V)\ra \JV^{l, p}(X_m, V_m)$, while the restriction maps $\JV^{l,p}(X_m, V_m)\ra \JV^{l, p}(X, V)$. However, for each $s$ we have an a priori topological bound on the maximum number $M$ of levels entering in the compactification $\ov\M_s(X, V)$. So for the purpose of transversality we can start with  
 Sobolev norms $W^{l+M, p}$ on the parameter space  $\J(X, V)$ and then take the limit as $l\ra \infty$ to get a smooth model for  the resolution of each stratum of $\ov\M_s(X, V)$.  
\end{rem}

For a family $\x\ra B$ of deformations of a level $m$ building, there is a space $\J(\x\ra B)$ of parameters compatible with both the fibration and the total divisor. Since $p$ is a singular fibration, this means that their restriction over each open stratum of $B$ (over which $p$ is a fibration) is  compatible with $p$ in the sense of \eqref{A.D.J.cpt.p}. The compatibility with the total divisor means that they lift to a pair of parameters on the smooth resolution.  There is also a corresponding Gromov-type perturbation space: 
\bear\label{A.D.J.defm.bd}
\xymatrix{
\JV(\x\ra B)\ar[r]&\ar@ /_1pc/ [l]   \J(\x\ra B)
}
\eear
Restricting to the fiber $X_m$ and then to $X$ in level zero defines fibrations
\best
\xymatrix{
\JV(\x\ra B)\ar[r]&\ar@ /_1pc/ [l]   \JV(X_m, V_m)\ar[r]&\ar@ /_1pc/ [l] \JV(X, V)
}
\eest
while the process of rescaling $(X, V)$ and then the symplectic sum of $X_m$ along the singular locus defines (non canonical) sections.
\smallskip

We end this section with local considerations involved in the construction of the perturbations 
$\nu\in {\mathcal V}(X, V)$ that enter the proof of transversality in Lemma \ref{L.d=d-2}. All required perturbations can be constructed supported in small balls about a fixed point, so it suffices to work on the standard local model of a normal crossing divisor at a point.
\begin{rem}\label{R.nu.can.be} 

For any element $\nu\in {\mathcal V} (X,V)$ and any map $f:C \ra X$, we can consider its restriction $\nu\in \Hom^{0,1}(TC , f^*TX)$ to the graph of $f$ (which is assumed to be an embedding). Around  any  depth $k$ point  $x\in C$ let $p=f(x)\in V^k$ and decompose $\nu$ into components $\nu_i$ one in each direction $i\in I(x)$ of $V$ at $f(x)$, and then $\nu^T$ in the tangent direction to the stratum $V^k$ containing $f(p)$. Working independently in these different directions one can construct required perturbations supported away from the higher depth stratum $V^{k+1}$. 

In particular, for any $\eta\in \Hom^{0,1}(T_xC, T_pV^k)$ there exists a $\nu\in {\mathcal V}(X, V)$ 
supported near $x$ such that $\nu(x, f(x))= \eta$, which is then used to prove transversality of the universal moduli space in the standard way. More generally, for any $\eta\in  \Hom^{0,1}(TC, f^*TX)$ which is supported in annuli about the contact points to $V$, there exists a perturbation 
$\nu\in {\mathcal V}(X, V)$ whose restriction to the graph of $f$ is $\eta$, which is used to prove transversality of the leading order section and of the refined evaluation map. 
\end{rem} 

\subsection{Stratifications}\label{S.A.stratif} Here we collect a few basic facts about the stratifications used in this paper, which are ubiquitous in GW theory. In short, if $X$ is a topological space, a {\em stratification} on it is a upper semi-continuous (USC) function $\si:X\ra \s$ to a partially ordered set (poset) $(\s,\fge)$; this simply means that  for each $s\in \s$  
\bear\label{A.D.X..fle.s}
X^{\fge s}=\{ x\in X\;|\; \si(x) \fge s\}
\eear 
are closed subsets of $X$, which we call the {\em closed strata}, and think of them as keeping track of points of $X$ whose `singularity' is at least $s$. The {\em open strata} are $X^s=\si^{-1}(s)$ (in general {\em not }open subsets of $X$!),  while  $X^{\fg s}$ is called the {\em boundary stratum} of $X^s$. (Warning: to match the standard conventions in combinatorial topology one should instead use the dual poset, i.e reverse the order on $\s$). 

A {\em forgetful map} is an order preserving map $\rho:\s'\ra\s$ (i.e. continuous); a stratification $\si':X\ra \s'$ such that $\si=\rho\circ \si'$ is called a {\em lift} or {\em refinement} of $\si$, and $\si=\rho_*\si'$ the {\em pushforward} of $\si'$. If $\rho$ is only USC, i.e $\rho_*\si' \fge \si$, then $\si'$ is called {\em finer} than $\si$.  

If $f:X\ra Y$ is a continuous map of topological spaces then $f^* \si_Y=\si_Y\circ f$ is the {\em pull-back} stratification. If  $\si_X$ is a stratification on $X$ which is constant along each fiber $f^{-1}(y)$ of $f$ when it descends to $Y$ to give the {\em push-forward} stratification $f_*\si_X$. 

Two stratifications $\si_i:X_i \ra \s_i$ induce a disjoint union stratification $\si_{1}\sqcup \si_{2}$ on $X_1\sqcup X_2$  and a product stratification $\si_1\ti\si_2$ on $X_1\ti X_2$. There is an equivariant version of the theory, where the group $G$ acts continuously on both $X$ and $\s$, and $s$ is $G$-equivariant, as well as a notion of morphism and isomorphism between two stratifications. 

If $\s$ has an initial element 0, $X^0$ is the {\em top stratum} of $X$.  If $\s$ has a final element, we also get a bottom stratum (the most `singular'). To each poset $\s$ we can always formally add an initial and a final object denoted $\pm \infty$ and then the top stratum is $X$ and the bottom stratum is empty.  

Finally, each poset $\s$ has an associated simplicial complex $\De(\s)$: 
the vertices of $\De(\s)$ are the elements of $\s$ and the faces of $\De(\s)$ are the chains in $\s$ (i.e. totally ordered subsets of $\s$). Conversely, each simplicial complex $\De$ has an associated  poset 
$\s(\De)$ of its nonempty faces ordered by inclusion. For each stratification $\si:X\ra \s$, the complex  $\De(\s)$ can be used to keep track of how the strata of $X$ fit together.

\begin{ex}\label{E.S.nc} If $V$ is a normal crossing divisor in $X$, $\mathsf{depth}: X\ra \N$ is a stratification, with the top stratum $X\setminus V$, and closed strata $V^k = X^{\ge k}$ where {\em at least} $k$ branches of $V$ meet, cf. \eqref{A.tower.V} and \eqref{A.D.X..fle.s}. A level $m$ building $X_m$ and its resolutions come with many more stratifications, discussed in detail in \S \ref{rem.def.multilevel}. For example, there is a stratification by multilevel
\bear\label{A.level.s}
\mathsf{level}:  X_m \ra \mathcal P
\eear
For each partition $\mu\in \mathcal P$ with $\mu_i$ parts of multiplicity $i$ and total degree $k=\sum_i \mu_i$,    the open stratum $\mathsf{level}^{-1} (\mu)$ consists of points $y\in X_m$ which project to a depth $k$ point $p\in X$ under the collapsing map $X_m\ra X$ and which are on level $i$ in $\mu_i$ of the $k$ normal directions to $V$ at $p$,  see \eqref{def.level.map}. 

The deformations $\x\ra B$ of a level $m$ building come similarly stratified: the base $B$ is a product of disks $(D^2, 0)$, thus carrying a product stratification, which induces by pullback a stratification on $\x$. The restriction of $\x\ra B$ over each open stratum is a fibration, whose fiber is a building with a fixed number of levels (equal to the depth of the stratum in $B$), with the top stratum corresponding  to buildings with 0 levels (unrescaled). 

\end{ex} 
\begin{ex}\label{E.S.X} If $X$ is any closed symplectic manifold, the moduli space $\oM(X)$ comes with many stratifications (cf Remark \ref{R.top.T.strat}), e.g.
\bear\label{D.stratif.X}
\mathsf{graph}: \oM(X) \ra \mathcal{T}, \quad \mathsf{nodes}:\oM(X) \ra \N,\quad \mathsf{n}:\oM(X) \ra \N, \quad 
\mathsf{degree}:\oM(X) \ra H_2(X, \Z) 
\eear
that associate to each stable map $f:C\ra X$ its dual graph, the number of nodes, marked  points of $C$ and respectively the homology class $f_*[C]$ it represents. There are forgetful maps that associate to each graph the collection of its edges, vertices, half edges or connected components, or just their cardinality: 
\bear\label{D.forget}
\mathsf{e}: \mathcal{T}\ra \N,   \quad \mathsf{v}:\mathcal{T}\ra \N,  \quad \mathsf{n}:\mathcal{T}\ra \N  \quad \text{ and }
\quad \mathsf{c}: \mathcal{T}\ra \N\quad
\eear
Since $\mathsf{nodes}= \mathsf{e}_*\mathsf{graph}$, the first stratification of \eqref{D.stratif.X} is a refinement of the second (and of the rest).

When $X=pt$ is a point, $\oM(X)=\oM$ is the Deligne-Mumford moduli space which comes with a natural normal crossing divisor (its nodal stratum); here the $\mathsf{nodes}$ and $\mathsf{depth}$ stratifications of $\oM$ are equal. If $\st:\oM(X)\ra \oM$ is the stabilization map, then on $\oM(X)$ the $\mathsf{ nodes}$ stratification is finer than it pullback from $\oM$, i.e.  $\mathsf{ nodes} \ge \st^*\mathsf{ nodes}$, recording the fact that the number of nodes of the domain $C$ is greater or equal than that of its image $\st(C)$. In general, without further information there is no other relation between the number of nodes of the domain and that of $\st(C)$ (i.e. $\mathsf{ nodes}$ is not a refinement of $\st^*\mathsf{nodes}$). Of course, in assuming that all domains are stable, we implicitly restricted to the case when these two stratifications are equal. On the other hand, the $\mathsf{graph}$ stratification of $\oM(X)$ is a refinement of $\st^*\mathsf{graph}$ because from the dual graph associated to $f$ we can read off not just the dual graph of its  domain $C$ but also that of $\st(C)$ (by inductively contracting all the unstable rational components). 

Another stratification of the moduli space 
\bear
\mathsf{reg}:\oM(X) \ra \N
\eear
associates to each stable map $f:C \ra X$ the dimension of the coker of the linearization $D_f$, with the top stratum consisting of regular points $f$ (i.e. for which $\cok D_f=0$). In fact this also comes in many flavors, one for each type of linearization $D_f$ considered e.g. depending whether we fix the domain or only its dual graph, whether we fix the image of the marked points in $X$ or not; for nodal domains we also have a choice of whether we consider only the linearization of the resolution $\wt f: \wt C\ra X$ or do we also linearize the matching conditions at the nodes, etc. Some of these stratifications are finer than others, but are not necessarily refinements of each other. One could go one step further and refine the dual graph stratification \eqref{D.stratif.X} by also keeping track of this information. 

Finally, one can also consider the stratification  
\bear\label{A.D.aut.X}
\mathsf{aut}:\oM(X) \ra \N
\eear
that associates to each stable map $f:C \ra X$ the order of its automorphism group $\Aut(f, C)$. It is in general finer than the pull back stratification $\st^*\mathsf{aut}$  from $\oM$. Of course, the stratification $\mathsf{aut}$ of the Deligne-Mumford moduli space can be refined to also keep track of the  automorphism groups (e.g. using the fact that $\oM$ is a global quotient orbifold). With more effort, \eqref{A.D.aut.X} could be refined to keep track of not just $\Aut (f, C)$ but also of the automorphism group $\Aut C$ of the (possibly unstable) domain $C$ together with its action on the components and special points of $C$ (i.e. on the dual graph) and also on the map $f$ (by reparameterization). We also avoided discussing this stratification by implicitly assuming all domains had trivial automorphisms. 
\end{ex} 
\begin{ex} If $V$ is a normal crossing divisor in $X$, all stratifications in Example \ref{E.S.X} extend to the  relative moduli space $\oM(X, V)$. The $\mathsf{depth}$ stratification of $(X,V)$ `lifts' to define a stratification 
\bear\label{D.stratif.cont}
\mathsf{contact}: \oM(X, V) \ra \mathcal{S}_{(X,V)}
\eear
associating to each relatively stable map $f:C\ra X_m$ its contact information $s$ to the zero divisor $V_m$. It can be refined to include the dual graph of $f$, now also decorated by the full contact information of $f$  to the total divisor $D_m$ (which includes the zero divisor $V_m$):
\bear\label{D.stratif.gr.V}
\mathsf{graph}: \oM(X, V) \ra \mathcal{T}_{(X,V)}
\eear
The evaluation/forgetful map that remembers only the contact information to the zero divisor relates these two stratifications by $\mathsf{ev}_{*} \mathsf{graph}= \mathsf{contact}$. In turn, these are by definition the pushforwards of the stratifications \eqref{M.ft.V.strat} associated to the resolution $\wt f:\wt C\ra \wt X_m$ regarded as an element of $\M(\wt X_m, \wt D_m)$. 
If $V\cup D$ is a normal crossing divisor where $D$ is a global branch,  the collapsing/forgetful map 
\bear\label{ft.d.v}
\ft_D:  \oM(X, V\cup D)\ra \oM(X, V)
\eear 
is compatible with the stratifications \eqref{D.stratif.cont}-\eqref{D.stratif.gr.V} using the corresponding forgetful maps on $\mathcal{T}$ or $\mathcal{S}$.  However, one has to be careful about the parameter spaces over which the maps  \eqref{ft.d.v} are defined, see \S \ref{A.spaces.param}.  For example, if we use only tamed pairs $\J$, then \eqref{ft.d.v} is defined over $\J(X, V\cup D)$,  but if we turn on Gromov type perturbations $\JV$, the map  \eqref{ft.d.v} is no longer defined on the entire parameter space $\JV(X, V\cup D)$, rather only on the subspace of perturbations which are independent of the extra contact points to $D$ (this causes extra complications in \cite{ip-don}). 

\end{ex} 
In all our examples above $\s$ is countable, and often equal to $\N$, or else comes with an (order preserving) map $\al: \s\ra \R$ (or a choice of `stability condition') with the following   {\em finiteness property}: for each $N\in \N$, there are finitely many elements $s\in \s$ with $\al(s)\le N$. The typical example is $\al:H_2(X, \Z) \ra \R$ defined by $\al(A)=\omega(A)$ the usual energy level, but that can also be refined to include the contribution from the domain (e.g. by keeping track of the total energy \eqref{energy}). 

For any stratification $\si:X\ra \mathcal S$ satisfying the finiteness property, the boundary strata are closed subsets of $X$, and the open stratum $X^s$ is an open subset of the closed stratum $X^{\fge s}$, which can be regarded as a compactification of $X^s$ (though it may contain more points than just limit points).  When $\si$ is continuous, the strata $X^s$ are both open and closed in $X$, thus $s$ indexes a {\em disjoint} union decomposition of $X$ into countably many pieces.  
\begin{ex} The $\mathsf{degree}$ \eqref{D.stratif.X}  is continuous, inducing $\oM(X)=\sqcup_{A\in H_2(X)}\oM_A(X)$.  When we allow disconnected domains, the number of connected components and their homology class are also continuous, so can be refined to $\mathsf{degrees}: \oM(X) \ra \mathrm{Sym}(H_2(X, \Z))$. The dual graph stratification \eqref{D.stratif.X} is a refinement of both, but is only USC in the Gromov topology, recording how the strata fit together.  
\end{ex}  

\begin{rem} Stratifications (of the type considered above) are closely related to filtrations. Any filtration $\dots \subseteq V^{k+1}\subseteq V^k\subseteq  \dots  \subseteq X$ of a topological space $X$ by closed sets $V^k$  induces a `depth' stratification $s:X\ra \Z$ by $s(x)=\max\{k\;|\; x\in V^k \}$, whose closed strata are $V^k$ (strictly speaking, the target of $s$ is $\Z$ if the union of all $V$'s is $X$ and their intersection is empty; otherwise the target  is $\Z\cup \{\pm \infty\}$). 

In fact, any collection of closed sets $\{ V_\al\}_{\al\in A}$ of $X$ is  part of a stratification  $s:X\ra \s_A$ whose closed strata are $V_I= \cap_{\al\in I}  V_{\al} $ with $V_\emptyset= X$, where $\s_A$ is the collection of subsets $I$ of $A$ (ordered by inclusion: $I\fle J$ iff $I\subseteq J$). When $A$ is finite, the forgetful map $o: \s_A\ra \N$ defined by $o(I)=|I|$ induces a coarser `depth' stratification $o_* s:X\ra \N$ whose closed strata $V^k$ keep track of at least how many of the original $V_\al$'s a point belongs to. 

Of course, any collection of open sets $\{ U_\al\}_{\al\in A}$ determines an `opposite' stratification, defined via their complements.  
\end{rem}
\begin{ex} \label{A.ex.V.I} Assume $V$ is a normal crossing divisor in $X$ with simple crossings. Then it has finitely many global branches $\{ V_\al\}_{\al \in A}$, each one a smooth compact submanifold of $X$. This induces a stratification $V_I$ of $X$  indexed by finite subsets of $A$ as described in the paragraph above. Because in this case the fiber $\iota^{-1}(p)$ of $\wt V=\sqcup_\al V_\al \ra V$ comes with an intrinsic injection $\iota^{-1}(p) \ra A$ to the indexing set $A$, the resolutions $V_I$ defined by 
\eqref{A.D.V.I} correspond (up to reordering of the branches) to $V_I=\cap_{\al\in I} V_\al$. 

 The coarser stratification $V^k$ by depth that forgets about the indexing of the branches of $V$ is nothing but the stratification/filtration \eqref{A.tower.V}. 
\end{ex} 
\begin{ex}\label{A.ex.V.tower} The tower \eqref{D.tower.build} of sub-buildings is another filtration by closed sets that induces a stratification $\mathsf{floor}: X_m \ra \Z$ of a level $m$ building, whose open stratum $\mathsf{floor}^{-1}(-i)=X_i\setminus X_{i-1}$ is the $i$'th floor of the building, with the top stratum $X$ (note the reversal in order). The tower of open sets 
\bear
X=\wt X_0\subset \wt X_1 \subset \wt X_2\subset \dots 
\eear
on the fundamental model \eqref{D2.rescaled} induces instead the `opposite' stratification \eqref{D.level.disk} on it by level, which gets promoted in a functorial fashion to the (multi)-level stratification in \S \ref{rem.def.multilevel}. The sign map \eqref{D.sign.disk} associated to the fundamental model is not a stratification (it is not USC), but can be replaced by a {\em pair} of stratifications $\ep_0, \ep_\infty:\wt X_m \ra \{0, 1\}$ on the fundamental model keeping track whether the coordinate of the point is 0 or not and respectively $\infty$ or not. They can be promoted in a functorial fashion to a stratification $\ep_0\ti \ep_\infty$ on any level $m$ building $\wt X_m$ carrying the same information as the multisign map. 
\end{ex} 
\begin{ex} The spaces of parameters $\J$ for the holomorphic map equation also come with various natural stratifications. One of them comes from the Sobolev completions $\J^{k, p}$ (with $p> 2$ and $k\ge 1$, partially ordered by the Sobolev embedding), with the continuos parameters being the initial object and the smooth parameters the final one. There is even a stratification (with only two strata) depending whether the pair $(\om,J)$ is compatible or just tame, or in the presence of a normal crossing divisor whether these parameters are compatible with the divisor or not, and with how many of its global branches. 

Gromov-type perturbations $\JV(X, V)$ also inherit stratifications from those on $\ov\U\ti X$, where the universal curve $\ov\U_{g, n}=\oM_{g, n+1}$ could have any one of the stratifications discussed above. The perturbation spaces $\JV(X_m)$ and $\JV(\x\ra B)$ associated to a level $m$ building $X_m$ and its deformations $\x\ra B$ similarly come with induced stratifications. As we have seen, these also come in several flavors, depending how much of the extra structure they are compatible with. 
\end{ex}

\end{appendix}



\end{document}